\documentclass[10pt]{article}

\usepackage{amsfonts,graphicx,color,epsfig,verbatim,multicol}

\usepackage[letterpaper,asymmetric,left=1.1in,right=1.1in,top=1.25in,bottom=1.25in,
bindingoffset=0.0in]{geometry}

\title{Constructing hyperbolic polyhedra using Newton's Method}

\author{Roland K. W. Roeder \footnote{rroeder@fields.utoronto.ca}\footnote{Partially supported by a U.S. National Defense Science and Engineering
Fellowship and by the Fields Institute as a Jerrold E. Marsden postdoctoral fellow.}}

\input{psfig.sty}

\newtheorem{thm}{Theorem}[section]
\newtheorem{prop}[thm]{Proposition}
\newtheorem{cor}[thm]{Corollary}
\newtheorem{lem}[thm]{Lemma}
\newtheorem{move}{Move}

\newcommand{\Endproof}{$\Box$ \vspace{.15in}}

\begin{document}

\maketitle

\begin{center} \bf Abstract \end{center}
{\small

We demonstrate how to construct three-dimensional compact hyperbolic polyhedra
using Newton's Method.  Under the restriction that the dihedral angles are
non-obtuse, Andreev's Theorem \cite{AND,AND2} provides as necessary and
sufficient conditions five classes of linear inequalities for the dihedral
angles of a compact hyperbolic polyhedron realizing a given combinatorial
structure $C$.  Andreev's Theorem also shows that the resulting polyhedron is
unique, up to hyperbolic isometry.  Our construction uses Newton's method and a
homotopy to explicitly follow the existence proof presented by Andreev,
providing both a very clear illustration of proof of Andreev's Theorem as well
as a convenient way to construct three-dimensional compact hyperbolic
polyhedra having non-obtuse dihedral angles.

As an application, we construct compact hyperbolic polyhedra having dihedral
angles that are (proper) integer sub-multiples of $\pi$, so that the group
$\Gamma$ generated by reflections in the faces is a
discrete group of isometries of hyperbolic space.  The quotient
$\mathbb{H}^3/\Gamma$ is hence a compact hyperbolic 3-orbifold, of which we
study the hyperbolic volume and spectrum of closed geodesic lengths using
SnapPea \cite{SNAPPEA}.  One consequence is a volume estimate for a
``hyperelliptic'' manifold considered in \cite{HYPER_ELL}.  }

\section{Introduction}

Andreev's Theorem \cite{AND,AND2} provides a complete characterization of
compact hyperbolic polyhedra having non-obtuse dihedral angles.  See also
\cite{ROE2,ROE} for an alternative exposition on the classical proof.  Other
approaches to Andreev's Theorem can be found by Rivin and Hodgson \cite{RH,H},
Thurston \cite{TH}, Marden and Rodin \cite{MR}, and Bowers and Stephenson
\cite{STEVE}.  In this paper we show that the classical proof from
\cite{AND,AND2,ROE2,ROE} is constructive when combined with Newton's Method for
solving nonlinear equations.  

Combinatorial descriptions of hyperbolic polyhedra that are relevant to
Andreev's Theorem fall into three classes, {\em simple, truncated}, and {\em
compound}, all defined later in this section.  The proof in \cite{ROE2,ROE}
provides an explicit continuous path in the space of polyhedra deforming a
given simple polyhedron $P$ to one of two which are easily constructed by hand:
the $N$-faced prism $Pr_N$ and the $N$-faced split prism $D_N$.  We use
Newton's method to follow such a path backwards deforming a computer
realization of $Pr_N$ or $D_N$ to a computer realization of the desired
polyhedron $P$.  This technique, which has been well studied in the literature,
is known as the homotopy method \cite{ALLGOWER,BLUM,SHUB1,SHUB2,SHUB3,SHUB4,BELTRAN}.
We illustrate the construction of simple polyhedra in Sections
\ref{SEC_WH},\ref{SEC_ALGORITHM}, and \ref{SEC_DIFFICULT}.

A similar deformation, again using Newton's method, allows us to construct
truncated polyhedra from simple polyhedra.  We demonstrate this deformation in
Section \ref{SEC_TRUNCATION}.  In Section \ref{SEC_COMPOUND} we show how to
construct a compound polyhedron as a gluing of two appropriate truncated
polyhedra.

In this way, our program graphically illustrates Andreev's proof of existence
for explicit examples.  In fact, writing this program and working
through Andreev's proof for some specific examples led to the detection an error
in the proof of existence, which has been corrected in
\cite{ROE2,ROE}.

A further benefit of this program is the construction of polyhedra whose
dihedral angles are proper integer sub-multiples of $\pi$.   As a consequence
of Poincar\'es polyhedron theorem \cite{POINCARE}, the group $\Gamma$ generated
by reflections in the faces of such a polyhedron is a discrete group of isometries of hyperbolic
space.  The quotient $\mathbb{H}^3/\Gamma$ is hence a compact hyperbolic
3-orbifold in which we study the hyperbolic volume and  spectrum of
closed geodesic lengths using SnapPea \cite{SNAPPEA}.  Such orbifolds and their
covering manifolds have been studied extensively
\cite{LOBELL,VESNIN_LOB,HYPER_ELL,LOB_VOL,SMALL_COVERS,KELLERHALS,RENI}.  In
fact the first example of a closed hyperbolic 3-manifold was obtained in this
way in 1931 by L\"obell \cite{LOBELL}.  One consequence of our study is a volume estimate
for a ``hyperelliptic'' manifold considered in \cite{HYPER_ELL}.

The reader should note that there are already excellent computer programs for
experimentation with hyperbolic 3-manifolds.  The program SnapPea
\cite{SNAPPEA} constructs hyperbolic structures on knot and link compliments,
as well as the hyperbolic Dehn surgeries on these compliments.  SnapPea
provides for the computation of a variety of geometry invariants of the
computed hyperbolic structure.  (See also \cite{WEEKS}.)  The program Snap
\cite{SNAP,SNAP_PAPER} provides a way of computing arithmetic invariants of
hyperbolic manifolds.  Both of these programs are quite easy to use and have
allowed for vast levels of experimentation, including a nice census of
low-volume hyperbolic manifolds and orbifolds.  

The experimentation done in this paper with the hyperbolic orbifolds obtained
from polyhedral reflection groups is very modest in comparison.  However, it is
an alternative way to construct hyperbolic structures on certain orbifolds
(and, in the future, possibly on manifold covers of these orbifolds) in a way
that these structures can nicely be studied by SnapPea (as well as Snap, and
other software, in the future). 

\vspace{.1in}

Let $E^{3,1}$ be $\mathbb{R}^4$ with the indefinite
metric $\Vert {\bf x} \Vert^2 = -x_0^2+x_1^2+x_2^2+x_3^2$.
In this paper, we work in the hyperbolic space $\mathbb{H}^3$ given by the
component of the subset of $E^{3,1}$ given by 

$$\Vert {\bf x} \Vert^2 = -x_0^2+x_1^2+x_2^2+x_3^2 = -1$$

\noindent
having $x_0 > 0$, with the Riemannian metric induced by the indefinite
metric

$$-dx_0^2+dx_1^2+dx_2^2+dx_3^2.$$

The hyper-plane orthogonal to a vector ${\bf v} \in
E^{3,1}$ intersects $\mathbb{H}^3$ if and only if $\langle{\bf v},{\bf
v}\rangle> 0$.  Let ${\bf v} \in E^{3,1}$ be a vector with $\langle{\bf
v},{\bf v}\rangle > 0$, and define

\begin{eqnarray*}
P_{\bf v} = \{{\bf w} \in \mathbb{H}^3 | \langle{\bf w},{\bf v}\rangle
= 0\} \mbox{      and      } 
H_{\bf v} = \{{\bf w} \in \mathbb{H}^3 | \langle{\bf w},{\bf v}\rangle \leq
0 \}
\end{eqnarray*}

\noindent to be the hyperbolic plane orthogonal to ${\bf v}$ and the
corresponding closed half space.

If one normalizes $\langle{\bf v},{\bf v}\rangle = 1$ and $\langle{\bf w},{\bf
w}\rangle = 1$ the planes $P_{\bf v}$ and $P_{\bf w}$ in $\mathbb{H}^3$
intersect in a line if and only if $\langle{\bf v},{\bf w}\rangle^2 <
1$, in which case their dihedral angle is $\arccos(-\langle{\bf v},{\bf
w}\rangle)$.  They intersect in a single point at infinity if and only if
$\langle{\bf v},{\bf w}\rangle^2 = 1$; in this case their dihedral angle is
$0$.

A {\it hyperbolic polyhedron} is an intersection

$$P = \bigcap_{i=0}^n H_{\bf v_i}$$

\noindent
having non-empty interior.

We will often use the Poincar\'e ball model of hyperbolic space, given by the
unit ball in $\mathbb{R}^3$ with the metric 

$$4\frac{dx_1^2+dx_2^2+dx_3^2}{(1 -\Vert {\bf x}\Vert^2)^2}$$ 

\noindent
and the upper half-space model of hyperbolic space,
given by the subset of $\mathbb{R}^3$ with $x_3 > 0$ equipped with the
metric 
$$\frac{dx_1^2+dx_2^2+dx_3^2}{x_3^2}.$$

\noindent
Both of these models are isomorphic to $\mathbb{H}^3$.

Hyperbolic planes in these models correspond to Euclidean hemispheres and
Euclidean planes that intersect the boundary perpendicularly.  Furthermore,
these models are conformally correct, that is, the hyperbolic angle between a
pair of such intersecting hyperbolic planes is exactly the Euclidean angle
between the corresponding spheres or planes.

\subsection{Combinatorial polyhedra and Andreev's Theorem}
A compact hyperbolic polyhedron $P$ is topologically a 3-dimensional ball, and
its boundary a 2-sphere $\mathbb{S}^2$.  The face structure of $P$ gives
$\mathbb{S}^2$ the structure of a cell complex $C$ whose faces correspond to
the faces of $P$.

Considering only hyperbolic polyhedra with non-obtuse dihedral angles
simplifies the combinatorics of any such $C$: 

\begin{prop} \label{TRIVALENT}
(a)  A vertex of a non-obtuse hyperbolic polyhedron $P$ is the
intersection of exactly 3 faces. \newline (b) For such a $P$, we can
compute the angles of the faces in terms of the dihedral angles; these
angles are also $\leq \pi/2$. 

\end{prop}
\noindent
See \cite{ROE2,ROE}.

The fundamental axioms of incidence place the following, obvious, further
restrictions on the complex $C$: 

\begin{itemize}

\item Every edge of $C$ belongs to exactly two faces.

\item A non-empty intersection of two faces is either an edge or a vertex.

\item Every face contains not fewer than three edges.

\end{itemize}

Any trivalent cell complex $C$ on $\mathbb{S}^2$ that satisfies the three
conditions above is an {\it abstract polyhedron}.  Since $C$ must be a trivalent
cell complex on $\mathbb{S}^2$, its dual, $C^*$, has only triangular faces and
the three above conditions ensure that it is a simplicial complex on
$\mathbb{S}^2$.  
The figure below shows an abstract  polyhedron $C$ drawn in the
plane (i.e.  with one of the faces corresponding to the region outside of the
figure.) The dual complex is also shown, in dashed lines.

\vspace{.1in}
\begin{center}
\begin{picture}(0,0)%
\epsfig{file=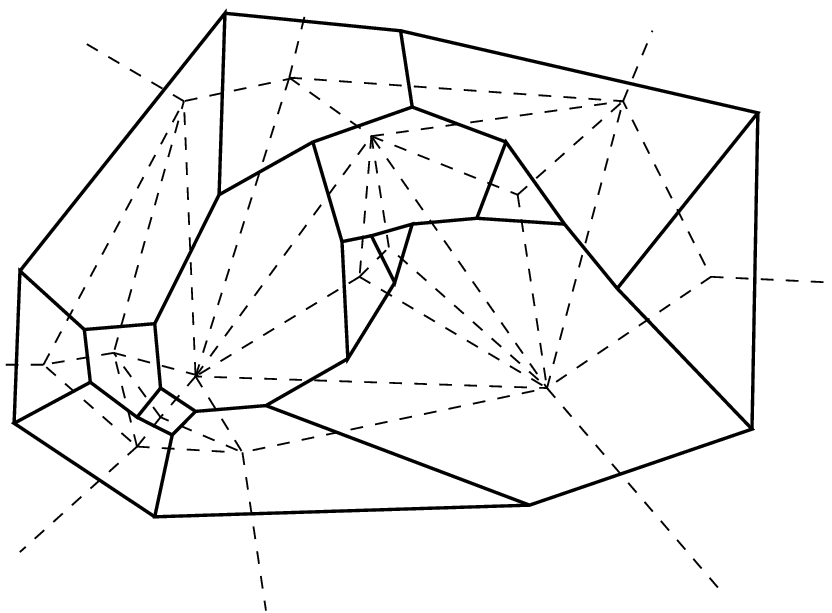}%
\end{picture}%
\setlength{\unitlength}{3947sp}%
\begingroup\makeatletter\ifx\SetFigFont\undefined%
\gdef\SetFigFont#1#2#3#4#5{%
  \reset@font\fontsize{#1}{#2pt}%
  \fontfamily{#3}\fontseries{#4}\fontshape{#5}%
  \selectfont}%
\fi\endgroup%
\begin{picture}(4015,2919)(139,-2143)
\end{picture}%

\end{center}
\vspace{.1in}

 We call a simple closed curve $\Gamma$ formed of $k$ edges of $C^*$ a
{\it k-circuit} and if all of the endpoints of the edges of $C$
intersected by $\Gamma$ are distinct, we call such a circuit a {\it
prismatic k-circuit}.  The figure below shows the same abstract polyhedron
as above, except this time the prismatic 3-circuits are dashed, the prismatic
4-circuits are dotted, and the dual complex is not shown.

\vspace{.1in}
\begin{center}
\begin{picture}(0,0)%
\epsfig{file=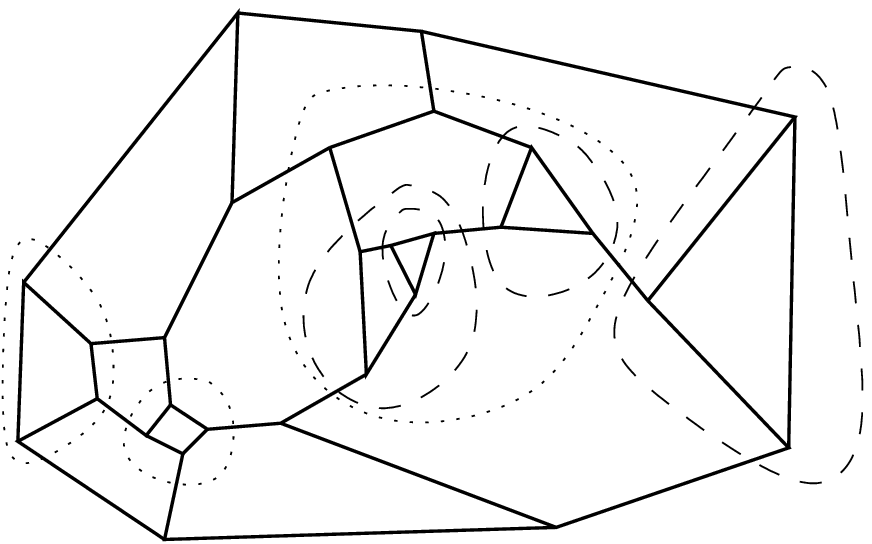}%
\end{picture}%
\setlength{\unitlength}{3947sp}%
\begingroup\makeatletter\ifx\SetFigFont\undefined%
\gdef\SetFigFont#1#2#3#4#5{%
  \reset@font\fontsize{#1}{#2pt}%
  \fontfamily{#3}\fontseries{#4}\fontshape{#5}%
  \selectfont}%
\fi\endgroup%
\begin{picture}(4152,2571)(80,-1853)
\end{picture}%

\end{center}
\vspace{.1in}

We say that a combinatorial polyhedron $C$ is {\em simple} if it has no prismatic $3$-circuits,
{\em truncated} if $C$ has prismatic $3$-circuits and each surrounds on one side a single
triangular face, and otherwise we call $C$ {\em compound}.  The combinatorial polyhedron
shown in the two above diagrams is compound.

\begin{thm} {\bf Andreev's Theorem}

Let $C$ be an abstract polyhedron with more than 4
faces and suppose that non-obtuse angles ${\bf a}_i$ are given corresponding
to edge $e_i$ of $C$.  There is a 
compact hyperbolic polyhedron $P$ whose faces realize $C$
with dihedral angle ${\bf a}_i$ at each edge $e_i$ if and only if
the following five conditions all hold:

\setcounter{enumi}{0}
\begin{enumerate}

\item For each edge $e_i$, ${\bf a}_i > 0$.

\item Whenever 3 distinct edges $e_i,e_j,e_k$ meet at a vertex, 
${\bf a}_i+{\bf a}_j+{\bf a}_k > \pi$.

\item Whenever $\Gamma$ is a prismatic 3-circuit intersecting edges
$e_i,e_j,e_k$, ${\bf a}_i+{\bf a}_j+{\bf a}_k < \pi$.

\item Whenever $\Gamma$ is a prismatic 4-circuit intersecting edges
$e_i,e_j,e_k,e_l$, then ${\bf a}_i+{\bf a}_j+{\bf a}_k+{\bf a}_l < 2\pi$.

\item Whenever there is a four sided face bounded by edges $e_1,$ $e_2,$
 $e_3,$ and $e_4$,
enumerated successively, with
edges $e_{12}, e_{23}, e_{34}, e_{41}$ entering the four vertices (edge
$e_{ij}$ connecting to the ends of $e_i$ and $e_j$), then:
$${\bf a}_1 + {\bf a}_3 + {\bf a}_{12} + {\bf a}_{23} + {\bf a}_{34} +
{\bf a}_{41} < 3\pi, \hspace{.2in} {\rm and}$$
$${\bf a}_2 + {\bf a}_4 + {\bf a}_{12} + {\bf a}_{23} + {\bf a}_{34} +
{\bf a}_{41} < 3\pi.$$
\end{enumerate}

Furthermore, this polyhedron is unique up to isometries of $\mathbb{H}^3$.
\end{thm}

\begin{cor}
\label{EX3APR} If $C$ is simple, i.e. has no prismatic $3$-circuits, there
exists a unique hyperbolic polyhedron realizing C with dihedral angles
$2\pi/5$.  \end{cor}

\vspace{.1in}

For a given $C$, let $E$ be the number of edges of $C$.  The subset
of $(0,\pi /2]^E$ satisfying these linear inequalities will be called the
{\it Andreev Polytope}, $A_C$.  Since $A_C$ is determined by linear
inequalities, it is convex.

Andreev's restriction to non-obtuse dihedral angles is emphatically necessary
to ensure that $A_C$ be convex.  Without this restriction, the corresponding
space of dihedral angles, $\Delta_C$, of compact (or finite volume) hyperbolic
polyhedra realizing a given $C$ is not convex \cite{DIAZ}.  In fact, the recent
work by D\'iaz \cite{DIAZ_ANDREEV} provides a detailed analysis of this space
of dihedral angles $\Delta_C$ for the class of abstract polyhedra $C$ obtained
from the tetrahedron by successively truncating vertices.  Her work nicely
illustrate the types of non-linear conditions that are necessary in a complete
analysis of the larger space of dihedral angles $\Delta_C$.

The work of Rivin \cite{RIV_IDEAL2, RIV_IDEAL1} shows that the space
of dihedral angles for ideal polyhedra forms a convex polytope, {\em without
the restriction to non-obtuse angles}.  (See also \cite{GUE}.)

Notice also that the hypothesis that the number of faces is greater than four
is also necessary because the space of non-obtuse dihedral angles for compact
tetrahedra is not convex \cite{ROE_TET}.  Conditions (1-5) remain necessary
conditions for compact tetrahedra, but they are no longer sufficient.

Bao and Bonahon \cite{BAO} prove a similar classification theorem for
hyperideal polyhedra.  Finally, the papers of Vinberg on discrete groups of
reflections in hyperbolic space \cite{AVS,VIN,VINREFL,VINVOL,VS} are also
closely related, as well as the work of Bennett and Luo \cite{LUO} and
Schlenker \cite{SCH2,SCH1,SCH3}.

\vspace{.1in}
Much attention has been focused on Andreev's Theorem from the viewpoint of
circle packings and circle patterns.  Given a polyhedron $P$ in the upper
half-space model of $\mathbb{H}^3$, the planes supporting the faces of $P$
intersect the boundary at infinity $x_3=0$ in a pattern of circles (and
straight lines) each with an orientation specifying ``on which side'' is the
polyhedron $P$.  Similarly, from such a pattern of circles and orientations one
can re-construct a polyhedron $P$.

The works of Thurston \cite{TH}, Marden and Rodin \cite{MR}, and Bowers and
Stephenson \cite{STEVE} all follow this approach to Andreev's Theorem.  In
fact, there is a beautiful computer program known as Circlepack
\cite{CIRCLE_PACK}, written by Ken Stephenson, that computes circle packings and
patterns of circles with specified angles of overlap.  All of the proofs from
this point of view use the conformal structure of the Riemann sphere
$\hat{\mathbb{C}} = \partial_\infty \mathbb{H}^3$ and use the correspondence
between conformal automorphisms of $\hat{\mathbb{C}}$ with isometries of
$\mathbb{H}^3$.

Instead of using the conformal structure on $\partial_\infty \mathbb{H}^3$, in
this paper we will work specifically with the metric structure of
$\mathbb{H}^3$.  (However, there is certainly some significant overlap with the
results in \cite{TH,MR,STEVE} and with the capabilities of the computer program
CirclePack \cite{CIRCLE_PACK}.)

\vspace{.1in}

We will now explain the implementation of a computer program whose input is
the combinatorial polyhedron $C$ and a dihedral angle vector ${\bf a} \in A_C$ and whose
output is a hyperbolic polyhedron realizing the pair $(C,{\bf a})$.

\subsection{An example}\label{SEC_EXAMPLE}

The following figure shows an explicit example of the data $(C,{\bf a})$
and the resulting polyhedron displayed in the conformal ball model using
the computer program Geomview \cite{GEO}.

\begin{figure}[htp]\label{BIG_FIG}
\begin{center}
\begin{picture}(0,0)%
\epsfig{file=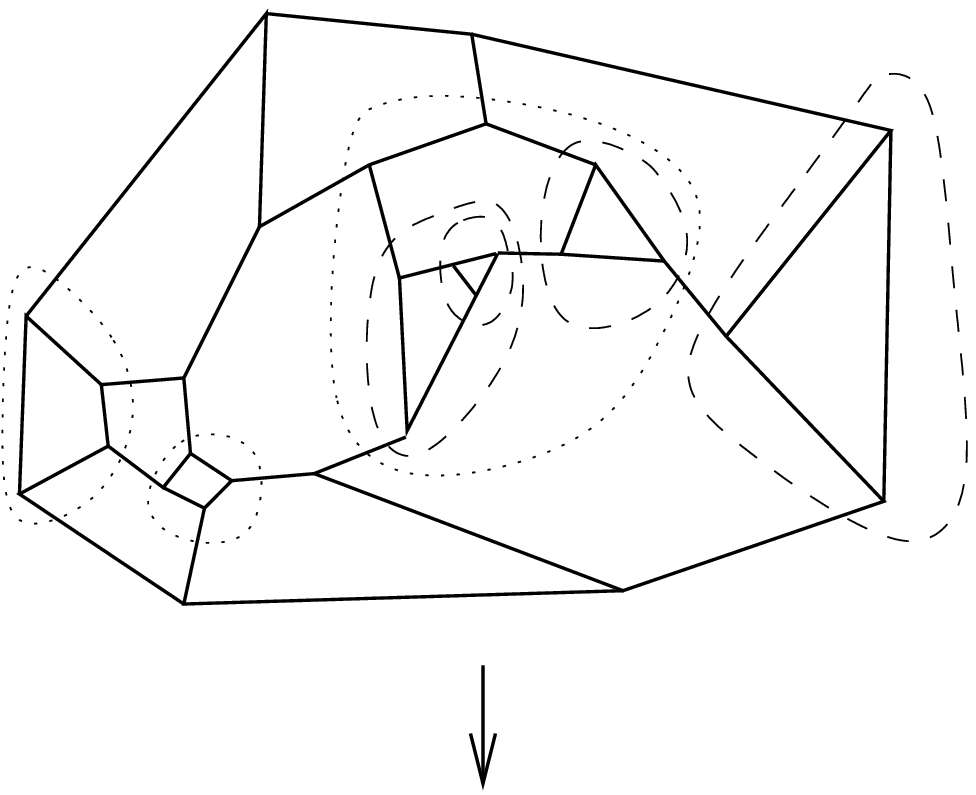}%
\end{picture}%
\setlength{\unitlength}{3947sp}%
\begingroup\makeatletter\ifx\SetFigFont\undefined%
\gdef\SetFigFont#1#2#3#4#5{%
  \reset@font\fontsize{#1}{#2pt}%
  \fontfamily{#3}\fontseries{#4}\fontshape{#5}%
  \selectfont}%
\fi\endgroup%
\begin{picture}(4654,3849)(82,-3083)
\put(878,-1751){\makebox(0,0)[lb]{\smash{{\SetFigFont{5}{6.0}{\familydefault}{\mddefault}{\updefault}{\color[rgb]{0,0,0}$\frac{2\pi}{5}$}%
}}}}
\put(2626,-2761){\makebox(0,0)[lb]{\smash{{\SetFigFont{12}{14.4}{\familydefault}{\mddefault}{\updefault}{\color[rgb]{0,0,0}Andreev's Theorem}%
}}}}
\put(581,-98){\makebox(0,0)[lb]{\smash{{\SetFigFont{5}{6.0}{\familydefault}{\mddefault}{\updefault}{\color[rgb]{0,0,0}$\frac{2\pi}{5}$}%
}}}}
\put(1101,-1498){\makebox(0,0)[lb]{\smash{{\SetFigFont{5}{6.0}{\familydefault}{\mddefault}{\updefault}{\color[rgb]{0,0,0}$\frac{2\pi}{5}$}%
}}}}
\put(720,-1049){\makebox(0,0)[lb]{\smash{{\SetFigFont{5}{6.0}{\familydefault}{\mddefault}{\updefault}{\color[rgb]{0,0,0}$\frac{2\pi}{5}$}%
}}}}
\put(1195,-729){\makebox(0,0)[lb]{\smash{{\SetFigFont{5}{6.0}{\familydefault}{\mddefault}{\updefault}{\color[rgb]{0,0,0}$\frac{2\pi}{5}$}%
}}}}
\put(1389,211){\makebox(0,0)[lb]{\smash{{\SetFigFont{5}{6.0}{\familydefault}{\mddefault}{\updefault}{\color[rgb]{0,0,0}$\frac{\pi}{2}$}%
}}}}
\put(3142,-308){\makebox(0,0)[lb]{\smash{{\SetFigFont{5}{6.0}{\familydefault}{\mddefault}{\updefault}{\color[rgb]{0,0,0}$\frac{\pi}{2}$}%
}}}}
\put(2716,-281){\makebox(0,0)[lb]{\smash{{\SetFigFont{5}{6.0}{\familydefault}{\mddefault}{\updefault}{\color[rgb]{0,0,0}$\frac{\pi}{2}$}%
}}}}
\put(3253,416){\makebox(0,0)[lb]{\smash{{\SetFigFont{5}{6.0}{\familydefault}{\mddefault}{\updefault}{\color[rgb]{0,0,0}$\frac{\pi}{4}$}%
}}}}
\put(2091,-81){\makebox(0,0)[lb]{\smash{{\SetFigFont{5}{6.0}{\familydefault}{\mddefault}{\updefault}{\color[rgb]{0,0,0}$\frac{2\pi}{5}$}%
}}}}
\put(3892,-560){\makebox(0,0)[lb]{\smash{{\SetFigFont{5}{6.0}{\familydefault}{\mddefault}{\updefault}{\color[rgb]{0,0,0}$\frac{\pi}{2}$}%
}}}}
\put(4381,-779){\makebox(0,0)[lb]{\smash{{\SetFigFont{5}{6.0}{\familydefault}{\mddefault}{\updefault}{\color[rgb]{0,0,0}$\frac{\pi}{2}$}%
}}}}
\put(3897,-1181){\makebox(0,0)[lb]{\smash{{\SetFigFont{5}{6.0}{\familydefault}{\mddefault}{\updefault}{\color[rgb]{0,0,0}$\frac{\pi}{2}$}%
}}}}
\put(3661,-1985){\makebox(0,0)[lb]{\smash{{\SetFigFont{5}{6.0}{\familydefault}{\mddefault}{\updefault}{\color[rgb]{0,0,0}$\frac{\pi}{4}$}%
}}}}
\put(2411,-1816){\makebox(0,0)[lb]{\smash{{\SetFigFont{5}{6.0}{\familydefault}{\mddefault}{\updefault}{\color[rgb]{0,0,0}$\frac{2\pi}{5}$}%
}}}}
\put(1786,-2248){\makebox(0,0)[lb]{\smash{{\SetFigFont{5}{6.0}{\familydefault}{\mddefault}{\updefault}{\color[rgb]{0,0,0}$\frac{2\pi}{5}$}%
}}}}
\put(1060,-1963){\makebox(0,0)[lb]{\smash{{\SetFigFont{5}{6.0}{\familydefault}{\mddefault}{\updefault}{\color[rgb]{0,0,0}$\frac{2\pi}{5}$}%
}}}}
\put(1141,-1696){\makebox(0,0)[lb]{\smash{{\SetFigFont{5}{6.0}{\familydefault}{\mddefault}{\updefault}{\color[rgb]{0,0,0}$\frac{2\pi}{5}$}%
}}}}
\put(945,-1578){\makebox(0,0)[lb]{\smash{{\SetFigFont{5}{6.0}{\familydefault}{\mddefault}{\updefault}{\color[rgb]{0,0,0}$\frac{2\pi}{5}$}%
}}}}
\put(1018,-1263){\makebox(0,0)[lb]{\smash{{\SetFigFont{5}{6.0}{\familydefault}{\mddefault}{\updefault}{\color[rgb]{0,0,0}$\frac{2\pi}{5}$}%
}}}}
\put(706,-1454){\makebox(0,0)[lb]{\smash{{\SetFigFont{5}{6.0}{\familydefault}{\mddefault}{\updefault}{\color[rgb]{0,0,0}$\frac{2\pi}{5}$}%
}}}}
\put(615,-1271){\makebox(0,0)[lb]{\smash{{\SetFigFont{5}{6.0}{\familydefault}{\mddefault}{\updefault}{\color[rgb]{0,0,0}$\frac{2\pi}{5}$}%
}}}}
\put(400,-1976){\makebox(0,0)[lb]{\smash{{\SetFigFont{5}{6.0}{\familydefault}{\mddefault}{\updefault}{\color[rgb]{0,0,0}$\frac{2\pi}{5}$}%
}}}}
\put(3387,-645){\makebox(0,0)[lb]{\smash{{\SetFigFont{5}{6.0}{\familydefault}{\mddefault}{\updefault}{\color[rgb]{0,0,0}$\frac{\pi}{4}$}%
}}}}
\put(401,-1604){\makebox(0,0)[lb]{\smash{{\SetFigFont{5}{6.0}{\familydefault}{\mddefault}{\updefault}{\color[rgb]{0,0,0}$\frac{2\pi}{5}$}%
}}}}
\put(1577,-309){\makebox(0,0)[lb]{\smash{{\SetFigFont{5}{6.0}{\familydefault}{\mddefault}{\updefault}{\color[rgb]{0,0,0}$\frac{2\pi}{5}$}%
}}}}
\put(1953,-292){\makebox(0,0)[lb]{\smash{{\SetFigFont{5}{6.0}{\familydefault}{\mddefault}{\updefault}{\color[rgb]{0,0,0}$\frac{\pi}{4}$}%
}}}}
\put(2436,-622){\makebox(0,0)[lb]{\smash{{\SetFigFont{5}{6.0}{\familydefault}{\mddefault}{\updefault}{\color[rgb]{0,0,0}$\frac{\pi}{2}$}%
}}}}
\put(2041,-862){\makebox(0,0)[lb]{\smash{{\SetFigFont{5}{6.0}{\familydefault}{\mddefault}{\updefault}{\color[rgb]{0,0,0}$\frac{\pi}{2}$}%
}}}}
\put(1374,-1660){\makebox(0,0)[lb]{\smash{{\SetFigFont{5}{6.0}{\familydefault}{\mddefault}{\updefault}{\color[rgb]{0,0,0}$\frac{2\pi}{5}$}%
}}}}
\put(225,-1275){\makebox(0,0)[lb]{\smash{{\SetFigFont{5}{6.0}{\familydefault}{\mddefault}{\updefault}{\color[rgb]{0,0,0}$\frac{2\pi}{5}$}%
}}}}
\put(1697,694){\makebox(0,0)[lb]{\smash{{\SetFigFont{5}{6.0}{\familydefault}{\mddefault}{\updefault}{\color[rgb]{0,0,0}$\frac{2\pi}{5}$}%
}}}}
\put(2661, 90){\makebox(0,0)[lb]{\smash{{\SetFigFont{5}{6.0}{\familydefault}{\mddefault}{\updefault}{\color[rgb]{0,0,0}$\frac{\pi}{4}$}%
}}}}
\put(2407,360){\makebox(0,0)[lb]{\smash{{\SetFigFont{5}{6.0}{\familydefault}{\mddefault}{\updefault}{\color[rgb]{0,0,0}$\frac{2\pi}{5}$}%
}}}}
\put(402,-932){\makebox(0,0)[lb]{\smash{{\SetFigFont{5}{6.0}{\familydefault}{\mddefault}{\updefault}{\color[rgb]{0,0,0}$\frac{2\pi}{5}$}%
}}}}
\put(2241,-1003){\makebox(0,0)[lb]{\smash{{\SetFigFont{5}{6.0}{\familydefault}{\mddefault}{\updefault}{\color[rgb]{0,0,0}$\frac{\pi}{3}$}%
}}}}
\put(2086,-504){\makebox(0,0)[lb]{\smash{{\SetFigFont{5}{6.0}{\familydefault}{\mddefault}{\updefault}{\color[rgb]{0,0,0}$\frac{\pi}{3}$}%
}}}}
\put(2972,-599){\makebox(0,0)[lb]{\smash{{\SetFigFont{5}{6.0}{\familydefault}{\mddefault}{\updefault}{\color[rgb]{0,0,0}$\frac{\pi}{2}$}%
}}}}
\put(1653,-1430){\makebox(0,0)[lb]{\smash{{\SetFigFont{5}{6.0}{\familydefault}{\mddefault}{\updefault}{\color[rgb]{0,0,0}$\frac{\pi}{3}$}%
}}}}
\put(2219,-670){\makebox(0,0)[lb]{\smash{{\SetFigFont{5}{6.0}{\familydefault}{\mddefault}{\updefault}{\color[rgb]{0,0,0}$\frac{\pi}{2}$}%
}}}}
\put(2276,-451){\makebox(0,0)[lb]{\smash{{\SetFigFont{5}{6.0}{\familydefault}{\mddefault}{\updefault}{\color[rgb]{0,0,0}$\frac{\pi}{2}$}%
}}}}
\put(2569,-406){\makebox(0,0)[lb]{\smash{{\SetFigFont{5}{6.0}{\familydefault}{\mddefault}{\updefault}{\color[rgb]{0,0,0}$\frac{\pi}{6}$}%
}}}}
\end{picture}%

\vspace{.1in}
\includegraphics[scale=.35]{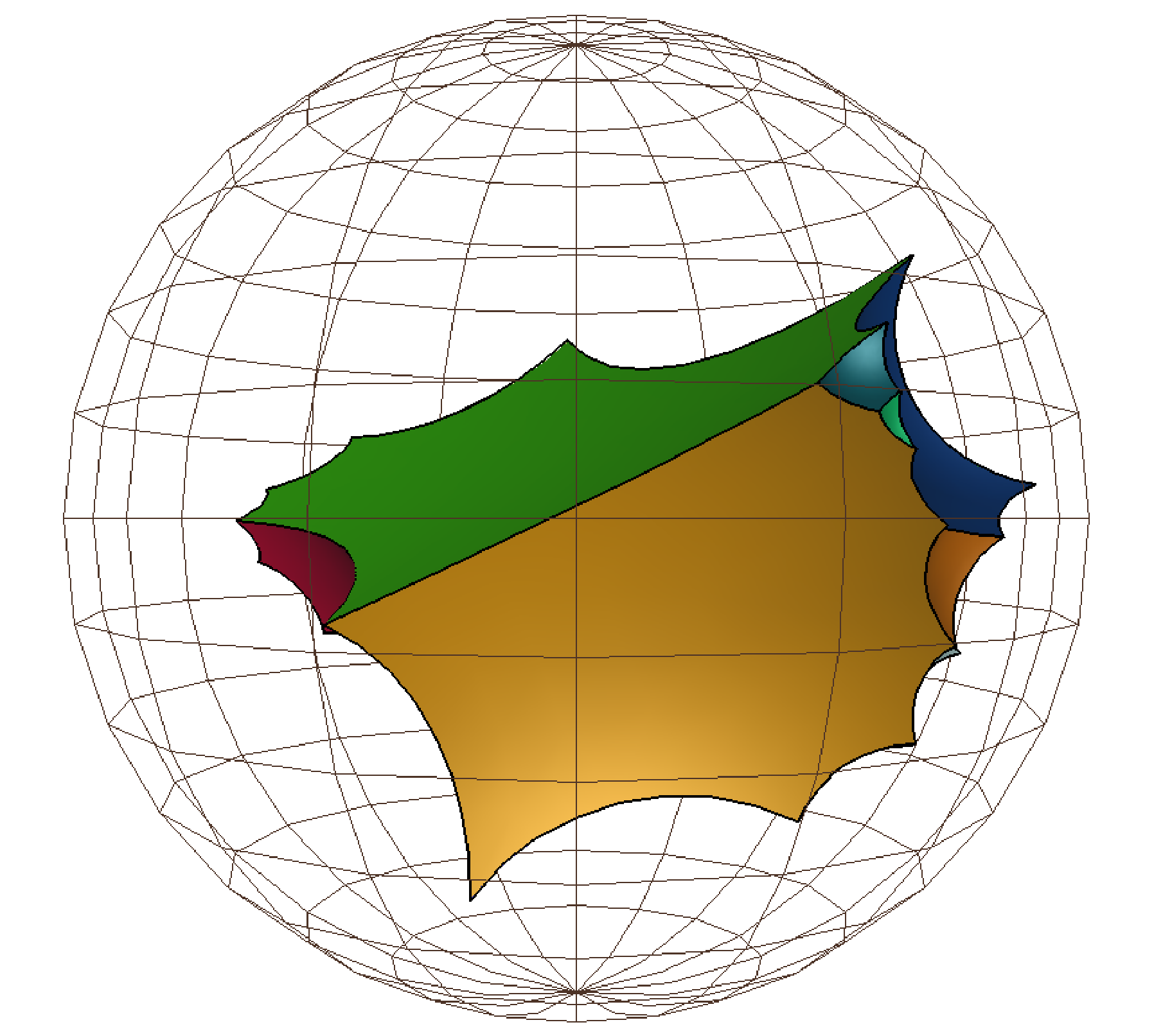}
\end{center}
\caption{\hbox{ }}
\end{figure}

\subsection{Outline of the proof of Andreev's Theorem}\label{SEC_OUTLINE}
In this section, we recall the major steps from the proof of Andreev's
Theorem that were presented in \cite{ROE2,ROE}.

Let $C$ be a trivalent abstract polyhedron with $N$ faces. We say that a
hyperbolic polyhedron $P \subset \mathbb{H}^3$ {\it realizes $C$} if there
is a cellular homeomorphism from $C$ to $\partial P$ (i.e., a homeomorphism
mapping faces of $C$ to faces of $P$, edges of $C$ to edges of
$P$, and vertices of $C$ to vertices of $P$). We will call each isotopy
class of cellular homeomorphisms $\phi : C \rightarrow \partial P$ a {\it
marking} on $P$.

We defined ${\cal P}_C$ to be the set of pairs $(P,\phi)$ so that
$\phi$ is a marking with the equivalence relation that $(P,\phi) \sim
(P',\phi ')$ if there exists an automorphism $\rho:  \mathbb{H}^3
\rightarrow \mathbb{H}^3$ such that $\rho(P) = P'$, and both $\phi '$
and $\rho \circ \phi$ represent the same marking on $P'$.

\begin{prop}
The space ${\cal P}_C$ is a manifold of dimension $3N-6$ (perhaps empty).
\end{prop}

\noindent
The proof is relatively standard and can be found in \cite{ROE2,ROE}.

Since the edge graph of $C$ is trivalent,
the number $E$ of edges of $C$ is the same as the dimension of
${\cal P}_C$.  
Given any $P \in {\cal P}_C$, let $\alpha(P) =
({\bf a}_1,{\bf a}_2,{\bf a}_3,...)$ be the $E$-tuple consisting of the dihedral
angles of $P$ at each edge (according to some fixed numbering of the edges of
$C$).  This map $\alpha$ is obviously continuous with respect to the topology
on ${\cal P}_C$, which it inherits from its manifold structure.

We let ${\cal P}_C^0$ be the subset of ${\cal P}_C$ consisting of polyhedra with
non-obtuse dihedral angles.  To establish Andreev's Theorem, we proved the
following statement:

\begin{thm}
For every abstract polyhedron $C$ having more than four faces, the mapping
$\alpha: {\cal P}_C^0 \rightarrow A_C$ is a homeomorphism.
\end{thm}

There were two major steps:

\begin{prop} \label{NONEMPTYIMPLIESAND}
If ${\cal P}^0_C \neq \emptyset$, then $\alpha : {\cal P}^0_C \rightarrow
A_C$ is a homeomorphism.
\end{prop}

We checked that $\alpha({\cal P}_C^0) \subset A_C$ by showing that conditions
(1)-(5) are necessary.  There is an open subset ${\cal P}^1_C \subset {\cal
P}_C$ containing ${\cal P}^0_C$ on which one can prove that $\alpha:{\cal P}^1_C \rightarrow
\mathbb{R}^E$ is injective, using a modification of Cauchy's rigidity for
Euclidean polyhedra.  This gives the uniqueness part of Andreev's Theorem.  Using invariance of domain,
it also gives that $\alpha : {\cal P}^1_C \rightarrow \mathbb{R}^E$ is a
local homeomorphism.  Because ${\cal P}^0_C \subset
{\cal P}^1_C$,  $\alpha$ restricted to ${\cal P}^0_C$ is a local homeomorphism, as
well.

We then showed that  $\alpha:{\cal P}_C^0 \rightarrow A_C$ is proper, which amounts
to showing that if a sequence of polyhedra $P_i$ in ${\cal P}_C^0$ degenerate
(i.e. leave  ${\cal P}_C^0$) then the sequence $\alpha(P_i)$ tends to $\partial A_C$. 
The fact that $\alpha:{\cal P}_C^0 \rightarrow A_C$ is a proper local homeomorphism
was sufficient to show that $\alpha({\cal P}_C^0)$ is open and closed in $A_C$.

\begin{prop}\label{EXIST}
If $A_C \neq \emptyset$, then ${\cal P}^0_C \neq \emptyset$.
\end{prop}

The second step was much more difficult because for each $C$ with non-empty
$A_C$ one needed to construct some polyhedron realizing $C$ (with non-obtuse
dihedral angles).  In fact the proof of Proposition \ref{EXIST} outlines a
scheme for how to construct a polyhedron realizing $C$.  The remainder of this
paper outlines how to follow this scheme explicitly on the computer using
Newton's Method and a homotopy method.

\section{A method for constructing Andreev polyhedra}

\subsection{Representing polyhedra on the computer}

All the constructions of polyhedra done in this paper are using Matlab
\cite{MAT} or the Free Software Foundation alternative Octave \cite{OCTAVE},
and all of the polyhedra are displayed in Geomview \cite{GEO}.  When doing
calculations, we represent a hyperbolic polyhedron $P$ having $N$ faces by
specifying $N$ outward pointing normal vectors vectors ${\bf v}_1,\cdots {\bf
v}_N$ each with $\langle {\bf v}_i,  {\bf v}_i \rangle = 1$, so that $P =
\bigcap_{i=1}^n H_{v_i}$.

Although such a list of $N$ vectors is sufficient to specify $P$, in order to
avoid repeated computation of the combinatorial structure of $P$ from these vectors
we additionally specify the adjacency matrix and a list
of all plane triples meeting at a vertex.  These three items are described
in a Matlab struct $P$, with $P.{\rm faces}$, $P.{\rm adjacency}$, and $P.{\rm
vert}$ holding the data mentioned above, respectively.

For example, the data for the polyhedron shown in Section \ref{SEC_EXAMPLE} is stored in
Matlab as:

\begin{multicols}{2}

{\tiny
\begin{verbatim}
New_poly = 
         vert: [28x3 double]
        faces: [16x4 double]
    adjacency: [16x16 double]

New_poly.faces = 
   36.5078  -10.7624   -0.3090  -34.8983
    4.9237   -1.5291   -1.2342   -4.6240
   -0.0000    0.8660   -0.5000   -0.0000
    4.5134   -2.1988   -2.3943   -3.2868
    2.7290   -1.9854   -0.3091   -2.1000
   13.3691   -5.0338   -0.3090  -12.4216
   65.0863  -19.6363   -2.7939  -61.9987
   35.9576   -9.6209   -1.5515  -34.6262
   51.5713  -13.3145   -0.3090  -49.8320
    5.8378   -0.5352   -0.3090   -5.8905
   -0.0000    0.0000    1.0000    0.0000
    1.2179   -1.4082   -0.7071    0.0000
   -7.8692    2.1329    3.6943    6.6879
    3.4039   -1.5744   -2.7269   -1.6344
   -1.0781    0.5773   -0.0000   -1.3524
   -2.1544    0.9964    1.7260   -1.2921

New_poly.adjacency =
  1   0   0   0   0   1   1   1   1   0   1   0   0   0   0   0
  0   1   1   1   1   1   1   1   0   1   0   0   1   0   0   0
  0   1   1   1   0   0   0   0   0   1   1   1   1   1   1   1
  0   1   1   1   1   0   0   0   0   0   0   1   1   1   0   0
  0   1   0   1   1   1   0   0   0   0   1   1   0   0   0   0
  1   1   0   0   1   1   1   0   0   0   1   0   0   0   0   0
  1   1   0   0   0   1   1   1   0   0   0   0   0   0   0   0
  1   1   0   0   0   0   1   1   1   1   0   0   0   0   0   0
  1   0   0   0   0   0   0   1   1   1   1   0   0   0   0   0
  0   1   1   0   0   0   0   1   1   1   1   0   0   0   0   0
  1   0   1   0   1   1   0   0   1   1   1   1   0   0   1   0
  0   0   1   1   1   0   0   0   0   0   1   1   0   1   1   1
  0   1   1   1   0   0   0   0   0   0   0   0   1   0   0   0
  0   0   1   1   0   0   0   0   0   0   0   1   0   1   0   0
  0   0   1   0   0   0   0   0   0   0   1   1   0   0   1   1
  0   0   1   0   0   0   0   0   0   0   0   1   0   0   1   1

New_poly.vert = 
     3     4    14
     6     5    11
     4    12     5
    12     5    11
     1     6     7
     1     7     8
     1     8     9
     1     9    11
    10     9    11
     1     6    11
     2     3    13
     2     4     5
     2     5     6
     2     6     7
     2     7     8
     9     8    10
     2     8    10
     2    10     3
    11    10     3
     2     4    13
     3     4    13
     3    12    14
     4    12    14
    11     3    15
    11    12    15
     3    12    16
     3    15    16
    12    15    16

\end{verbatim}
}

\end{multicols}

We display the polyhedra in Geomview using the hyperbolic mode and specifying
the conformal ball model.  The file format most convenient for our use is the
Object File Format, ``New\_poly.off.''  The first line of an Object format file
specifies the number of vertices, the number of faces, and the number of edges
of $P$ in that order: ${\rm num\_vert}$ ${\rm num\_faces}$ ${\rm num\_edges}.$ 

The next block of data is a list of the coordinates of vertices as points in
the unit ball.  (In fact, these are the coordinates of points in the Projective
Model for $\mathbb{H}^3$, not the Poincar\'e ball model which we describe in
the introduction.)  The last block of data is a list of the faces with each
face given by ${\rm vertex}_1$ ${\rm vertex}_2$ ... ${\rm vertex}_n$ ${\rm
colorspec}$, where the faces is spanned by  ${\rm vertex}_1$ ${\rm vertex}_2$
...  ${\rm vertex}_n$ and ${\rm colorspec}$ is an integer telling Geomview what
color to assign to the face.

\begin{multicols}{2}

{\tiny
\begin{verbatim}

28 16 42

0.093414 0.626297 -0.759378
0.668701 -0.423986 0.508660
0.480895 0.480927 -0.729094
0.533074 -0.046431 -0.835831
0.000602 -0.321164 0.944739
-0.109909 -0.298793 0.946284
-0.257482 -0.413626 0.871178
-0.241511 -0.517345 0.817366
-0.394737 -0.502851 0.762029
0.039860 -0.522084 0.841280
-0.198257 0.894806 -0.304826
0.473945 0.834475 -0.252208
0.632626 0.016126 0.705772
0.030462 -0.199457 0.975695
-0.112537 -0.193119 0.972142
-0.373893 -0.377061 0.844031
-0.376723 -0.130280 0.905450
-0.802208 0.544817 0.123196
-0.869841 -0.223571 -0.241933
0.160888 0.910483 -0.333052
0.007468 0.802257 -0.570580
0.104069 0.367146 -0.917226
0.301741 0.507113 -0.798909
-0.023848 -0.005394 -0.995609
0.158629 -0.006953 -0.985479
0.065661 0.192695 -0.976218
0.028183 0.094604 -0.992919
0.091162 0.094314 -0.989423

5 5 6 7 9 4 1
8 17 16 14 13 12 11 19 10 2
9 20 10 17 18 23 26 25 21 0 3
6 20 19 11 2 22 0 4
5 3 2 11 12 1 5
5 9 4 13 12 1 6
4 5 14 13 4 7
5 6 15 16 14 5 8
4 7 8 15 6 9
5 15 16 17 18 8 10
8 3 24 23 18 8 7 9 1 11
7 3 24 27 25 21 22 2 12
3 19 20 10 13
3 21 22 0 14
4 24 27 26 23 15
3 26 27 25 16
\end{verbatim}
}

\end{multicols}

Something is lacking when viewing the polyhedra displayed in the
two-dimensional images shown in this paper.  To alleviate this difficulty, the
Matlab and OFF files associated to each polyhedron that is constructed in this
paper are included as supplementary materials on the website of Experimental
Mathematics.  See the website \cite{GEO} for full details on the use of
Geomview.

\subsection{The desired polyhedron as a solution to $4N$ quadratic equations in $4N$ unknowns}
The proof of Andreev's Theorem gives that $\alpha_C: {\cal P}_C^0 \rightarrow
A_C$ is a homeomorphism, so the problem of constructing a polyhedron $P$
realizing $(C,{\bf a})$ can be expressed as the problem of finding a solution $P$ to
the equation $\alpha_C(P) = {\bf a}$.

Instead of working in ${\cal P}_C^0$, we write the desired polyhedron as a
solution of a system of $4N$ quadratic equations in $4N$ variables, where $N$ is the number of
faces of $C$.  Our solution is $N$ vectors $v_1,\cdots,v_N \in E^{3,1}$
satisfying
\begin{itemize}
\item{$\langle v_i,v_i \rangle = 1$ and}
\item{$\langle v_j,v_j \rangle = -\cos({\bf a}_{i,j})$ if faces $i$ and $j$ are
adjacent in $C$ and their common edge is assigned dihedral angle ${\bf a}_{i,j}$.}
\end{itemize}
These equations impose $E+N$ conditions on $4N$ variables, where $C$ has $N$
faces and $E$ edges.  

As mentioned in Section \ref{SEC_OUTLINE}, we have $E=3N-6$, so we have imposed
$4N-6$ conditions on $4N$ variables.  We impose $6$ additional conditions in
order to have the same number of equations and unknowns.  We normalize by
requiring that a chosen vector $v_i$ perpendicular to one of the faces agree with some
given $v$ (where $v$ is chosen so that $\langle v,v \rangle = 1$.)  We then
require that one of the vertices on the face perpendicular to $v_i$ is at a
given point $w$ in the plane $P_{v}$ and that a vertex adjacent to this vertex
be on a given line $l$ in $P_{v}$ through $w$.  One can check that these
normalizations provide $3$, $2$, and $1$ additional equations respectively.
(Notice that the six equations for this normalization are each linear.)  We denote
the normalization by a triple $(v,w,l)$.

We denote the resulting quadratic map by $F_{C,(v,w,l)}:\mathbb{R}^{4N}
\rightarrow \mathbb{R}^{4N}$.  Typically we will only mention the normalization
when necessary.  We denote the conditions described above for the right hand
side of the equations $F(x)=y$ as $({\bf a},0)$, where the ${\bf a}$ from this pair is
shorthand for the conditions $\langle v_i,v_i \rangle = 1$ and $\langle v_j,v_j
\rangle = -\cos(\alpha_{i,j})$ if faces $i$ and $j$ are adjacent in $C$, and
the $0$ represents the fact that the normalization $(v,w,l)$ is satisfied.

Andreev's Theorem asserts that if ${\bf a} \in A_C$, there is a real solution
to $F_{C,(v,w,l)}(x) = ({\bf a},0)$ corresponding to $N$ vectors ${\bf
v}_1,\cdots,{\bf v}_N$ in $E^{3,1}$ so that $P = \bigcap_{i=0}^n H_{\bf v_i}$
realizes the pair ($C$,${\bf a})$.  

\vspace{0.1in}
There are many sensible ways to numerically solve a system of quadratic
equations in the same number of equations as unknowns.  These include the
pre-packaged non-linear solvers in Matlab, Maple, and Mathematica, Newton's
Method, as well as Groebner basis techniques and fancier quadratically
constrained solvers.  

The difficulty is that with $4N$ quadratic equations in $4N$ unknowns, Bezout's
Theorem states that there will typically be $16N^2$ solutions.  On their own,
these solvers cannot easily be adapted to find the specific solution
corresponding to a convex polyhedron, without first finding all solutions (or
at least all real solutions) and then examining each solution to check if it
corresponds to the desired polyhedron.  Since some solutions may be much harder
to find than others, one could spend significant computation times pursuing
solutions that aren't of interest.

One way to ensure that the solution does correspond to a compact convex
polyhedron is to use an iterative method, like Newton's Method, for which an initial
condition that is sufficiently close to a given solution is garunteed to converge to
that root, in combination with a homotopy that garuntees that the nearest root
is always the root that corresponds to a compact convex polyhedron.  This is our
approach, which we describe in greater detail in the next few sections of the
paper.  We are not entirely sure that this method is faster than finding all of
the roots with ``brute force'' and then checking each solution to see if it is
the desired one, but our approach has the additional benefit that
it explicitly follows  Andreev's proof of existence, providing insight into
how this proof works for specific examples.

\subsection{Newton's Method and Homotopy methods}
Given two vector spaces $V$ and $W$ of the same dimension and a
mapping $F: V \rightarrow W$, the associated Newton map $N_F: V \rightarrow V$
is given by the formula \begin{eqnarray}\label{NEWTON} N_F({\bf x}) = {\bf x} -
[DF({\bf x})]^{-1}(F({\bf x})).  \end{eqnarray} \noindent 
\noindent
If the roots of $F$ are non-degenerate, i.e. $DF(r_i)$ is invertible for each
root $r_i$ of $F$, then the roots of $F$ corresponds bijectively to 
super-attracting fixed points of $N_F$.

Kantorovich's Theorem \cite{KANT} gives a precise lower bound
on the size of the basin of attraction for a root.

\begin{thm} {\rm ({\bf Kantorovich's Theorem})}\label{KANT}.
Let ${\bf a_0}$ be a point in $\mathbb{R}^n$, $U$ an open neighborhood of ${\bf a_0}$ in $\mathbb{R}^n$,
and $F:U \rightarrow \mathbb{R}^n$ a differentiable mapping with $[DF({\bf a_0})]$ invertible.

Let $U_0$ be the open ball of radius $\vert [DF({\bf a_0})]^{-1}F({\bf a_0}) \vert$ centered at ${\bf a_1} = N_F({\bf a_1})$.
If $U_0 \subset U$ and $[DF({\bf x})]$ satisfies the Lipshitz condition
$\parallel DF({\bf u_1}) - DF({\bf u_2})  \parallel \leq M \vert {\bf u_1} - {\bf u_2}  \vert$ for all ${\bf u_1, u_2} \in U_0$,
and if  the inequality
\begin{eqnarray}\label{KANTOROWICH}
\vert F({\bf a_0}) \vert \cdot \vert [DF({\bf a_0})]^{-1} \vert ^2 M \leq \frac{1}{2}
\end{eqnarray}

\noindent
is satisfied, then the equation $F({\bf x}) = 0$ has a unique solution in $U_0$, and Newton's Method
with initial guess ${\bf a_0}$ converges to it.
\end{thm}

For a proof of Kantorovich's Theorem see
\cite{HBH} or the original source \cite{KANT}.

While the dynamics near a fixed point can be easily understood by
Kantorovich's Theorem, the global dynamics of Newton's Method can be very
complicated, with loci of indeterminacy, and critical curves where $DN$ is not
injective.  In fact, the dynamics of Newton's method to solve for the common
roots of a pair of quadratic polynomials in $\mathbb{C}^2$ is a field of active
research \cite{HP,RNEWT}.  We expect that the global dynamics of the Newton map
to solve $F_{C,(v,w,l)}(x) = ({\bf a},0)$  is even significantly more
complicated than those in \cite{HP,RNEWT}.  In particular, we have no reason to
expect that a general initial condition in $R^{4N}$ will converge under
iteration of $N_F$ to any solution of $F_{C,(v,w,l)}(x) = ({\bf a},0)$ nor to
the specific solution representing a convex compact polyhedron $P$.

An approach that can sometimes be used to avoid the difficulties with the
global dynamics of Newton's Method is the homotopy method.   Suppose that you
want to solve $g(x)=y$.  The idea is to replace this equation by a family that
depends continuously on a single variable:
\begin{eqnarray*} g_t(x_t) = y_t
\end{eqnarray*} 
\noindent so that $g_1$ is the same function as $g$ and $y_1 =
y$, while $g_0(x) =y_0$ is an equation for which you already know a solution $x_0$.

Choose $k$ points $0=t_1,t_2,\ldots t_k=1$.  If $k$ is sufficiently large, then
$x_{t_1}$ may be in the basin of attraction of the Newton's method for
$f_{t_2}(x) = y_{t_2}$.  In this case, you can solve for $x_{t_2}$ and can
attempt to solve for $x_{t_3}$ using Newton's Method for $f_{t_3}(x) = y_{t_3}$
with initial condition $x_{t_2}$.  Repeating this procedure, if possible, leads
to the solution $x_1 = x_{t_k}$.

While this is obviously a very powerful method, there are many difficulties
choosing appropriate paths $g_t(x_t) = y_t$ and appropriate subdivisions
$0=t_1,t_2,\ldots,t_k=1$.  It is necessary to check that the conditions for
Kantorovich's Theorem are satisfied by $x_{t_j}$ for the equation
$f_{t_{j+1}}(x) = y_{t_{j+1}}$.  The biggest difficulty is to avoid the
situation where the derivative $\frac{\partial}{\partial x} f_t$ is singular
for some $t$.  Such points are described as being in the {\em discriminant
variety} and choosing paths that avoid the discriminant variety is a big
program of research.  These difficulties are discussed extensively by many
authors including Shub and Smale in
\cite{BLUM,ALLGOWER,SHUB1,SHUB2,SHUB3,SHUB4,BELTRAN}. 

\vspace{.1in}
The proof of Andreev's Theorem in \cite{ROE2,ROE} provides an explicit path
that we can use for a homotopy method to construct any simple polyhedron $P$ as
a continuous deformation of either the prism $Pr_N$ or the split prism $D_N$,
both of which can be easily constructed ``by hand''.  We will use this path for
our homotopy method: repeatedly using a polyhedron realizing a point on the
path as initial condition and solving for a polyhedron slightly further on the
path, chosen so that the dynamics of Newton's method converges to the correct
solution of $F$.

With a similar path we can use the homotopy method again to construct any
truncated polyhedron for which $A_C \neq \emptyset$.  We take a continuous
deformation of a simple polyhedron until the vertices to be truncated pass
$\partial_\infty \mathbb{H}^3$, and then add a finite number of additional
triangular faces intersecting the appropriate triples of faces perpendicularly.
Compound polyhedra are then constructed as gluings of a finite number of
truncated polyhedra.

\begin{prop}
The quadratic equation $F$ has a uniform Lipshitz constant on $R^{4N}$
depending only on the combinatorics $C$.
\end{prop}

\noindent
{\bf Proof:}
The proof is merely the observation that $F$ is quadratic so each of the second
derivatives are constant.
\Endproof

While we have checked that $F$ is Lipshitz, we make no effort to bound the
norm of the derivative $[DF]$ away from zero (hence avoiding the discriminant
variety).  In fact, for a typical problem this is very hard to do.  Instead, we
merely try the homotopy method with the path mentioned in the preceding
paragraph and we show that the method works for all of the constructions that
we attempt.  It may be interesting to provide a more rigorous basis for our use
of Newton's Method and the current choice of path.

\subsection{Deforming a given polyhedron using Newton's Method}\label{SEC_DEFORM}

Given a polyhedron $P$ realizing $C$ with dihedral angles ${\bf a} \in A_C$, it
is easy to use Newton's method to deform $P$ into a new polyhedron $P'$ having
any other angles ${\bf a}' \in A_C$.  Since $A_C$ is a convex polytope, 
choose the line segment between ${\bf a}$ and ${\bf a}'$ and subdivide
this segment into $K$ equally distributed points ${\bf a} = {\bf a}^0,{\bf
a}^1,{\bf a}^2,\cdots {\bf a}^{K-1}={\bf a}'$.  Then we use Newton's method
with initial condition corresponding to $P$ to solve for a polyhedron $P_1$
with dihedral angles ${\bf a}^1$.  We then repeats, using $P_1$ as initial
condition for Newton's method to solve for a polyhedron $P_2$ with dihedral
angles ${\bf a}^2$, and continue until reaching $P'$ realizing ${\bf a}'$.  If
the homotopy method has worked, then each step of Newton's method converges;
otherwise we can try a larger number of subdivisions $K$, or attempt to check
if the path has hit the discriminant variety.

In all of the calculations within this paper, when deforming the angles of a
given polyhedron $P$ within $A_C$, we use $K$ between $100$ and $300$
subdivisions, although this is sometimes a significant overkill.  

We consider it sufficient to show how to use Newton's method to construct some
polyhedron $P$ with non-obtuse dihedral angles for every $C$ that has $A_C \neq
\emptyset$.  From this $P$ one can construct any other $P' \in {\cal P}_C^0$
using the deformation described above. (This ease with which one can deform the
angles of a given polyhedron is an additional benefit of our homotopy method.)

In the next sections we will see how to connect individual paths
in $A_{C_1},\cdots,A_{C_k}$ so as to construct compact polyhedra realizing
$C_1$ a a sequence of deformations of a compact polyhedron realizing $C_k$.

\subsection{Simple polyhedra and Whitehead moves}\label{SEC_WH}

Recall that if $C$ is simple then
$\left(\frac{2\pi}{5},\cdots,\frac{2\pi}{5}\right) \in A_C$.  The goal of this
section and the following is to demonstrate the construction of a polyhedron
$P$ realizing any simple $C$ with these dihedral angles.

Andreev's Theorem provides a sequence of elementary changes (Whitehead moves) to the
reducing the combinatorics $C$ to one of two the
combinatorial polyhedra $D_N$ or $Pr_N$ depicted  below.

\vspace{.1in} \begin{center} 
\begin{picture}(0,0)%
\includegraphics{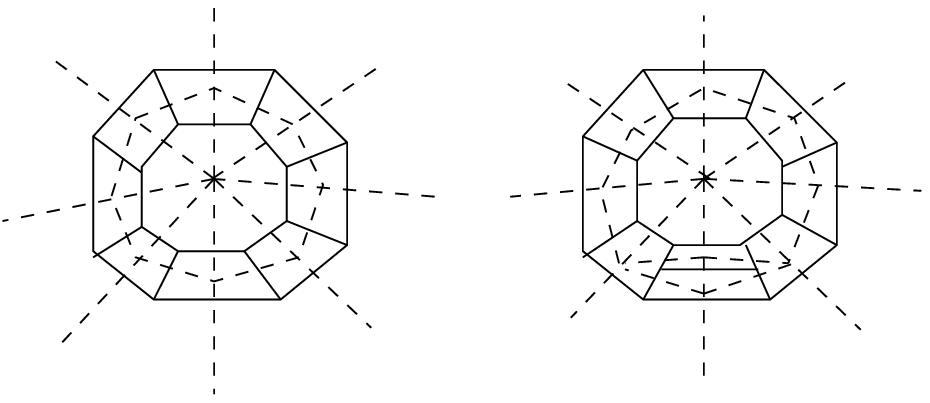}%
\end{picture}%
\setlength{\unitlength}{4144sp}%
\begingroup\makeatletter\ifx\SetFigFont\undefined%
\gdef\SetFigFont#1#2#3#4#5{%
  \reset@font\fontsize{#1}{#2pt}%
  \fontfamily{#3}\fontseries{#4}\fontshape{#5}%
  \selectfont}%
\fi\endgroup%
\begin{picture}(4224,2022)(1114,-1981)
\put(1669,-1968){\makebox(0,0)[lb]{\smash{\SetFigFont{7}{8.4}{\familydefault}{\mddefault}{\updefault}{\color[rgb]{0,0,0}Prism with 10 faces}%
}}}
\put(3769,-1950){\makebox(0,0)[lb]{\smash{\SetFigFont{7}{8.4}{\familydefault}{\mddefault}{\updefault}{\color[rgb]{0,0,0}Splitprism with 11 faces}%
}}}
\end{picture}
\end{center}

In this section we show how to create polyhedra realizing $D_N$ and $Pr_N$ and
how to do the Whitehead moves using Newton's method.
\vspace{.1in}

\begin{lem} \label{EXPRISM} Let $Pr_N$ and $D_N$ be the
abstract polyhedra corresponding to the $N$-faced prism and the $N$-faced
``split prism'', as illustrated below.  If $N > 4$, ${\cal P}^0_{Pr_N}$ is
nonempty and if $N > 7$, ${\cal P}^0_{D_N}$ is nonempty.
\end{lem}

\noindent
{\bf Construction:}
Construct a regular polygon with $N-2$ sides in the disc model for
$\mathbb{H}^2$.  ($N-2 \geq 3$, since $N \geq 5$.)  We can do this with
the angles arbitrarily small.  Now view $\mathbb{H}^2$ as the equatorial
plane of $\mathbb{H}^3$, and consider the hyperbolic planes perpendicular
to the equatorial plane containing the sides of the polygon.  In Euclidean
geometry these are hemispheres with centers on the boundary of the
equatorial disc.  The dihedral angles of these planes are the angles of
the polygon.

Consider two hyperbolic planes close to the equatorial plane, one
slightly above and one slightly beneath, both perpendicular to the
$z$-axis.  These will intersect the previous planes at angles slightly
smaller than $\pi/2$.  The region defined by these $N$ planes makes a
hyperbolic polyhedron realizing the cell structure of the prism.  Note that
our construction completes the proof of Proposition \ref{NONEMPTYIMPLIESAND},
for the special case $C = Pr_N$, $N \geq 5$.  

\begin{figure}[htp]
\begin{center}
\includegraphics[scale=1.0]{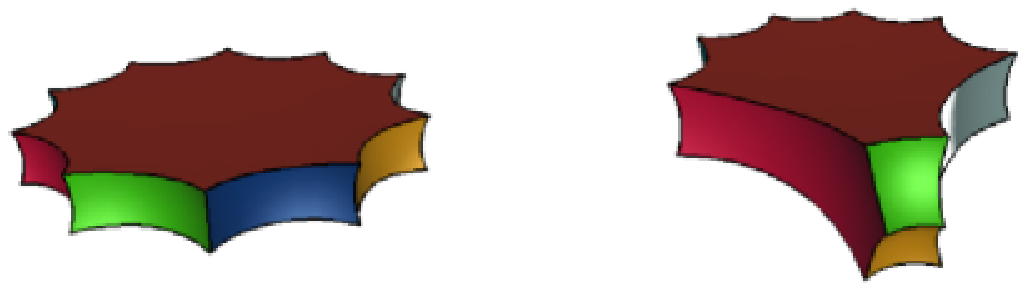}
\end{center}
\caption{\hbox{ }}
\end{figure}

For $N>7$ we will construct $D_N$ by cutting it into two
prisms each with $N-1$ faces and the dihedral angles shown below.

\begin{figure}[htp]
\begin{center} 
\begin{picture}(0,0)%
\epsfig{file=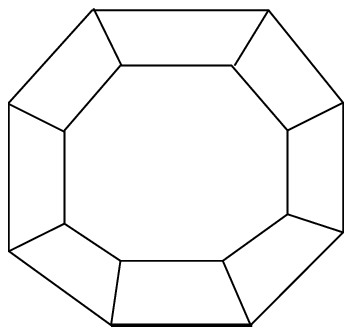}%
\end{picture}%
\setlength{\unitlength}{4144sp}%
\begingroup\makeatletter\ifx\SetFigFont\undefined%
\gdef\SetFigFont#1#2#3#4#5{%
  \reset@font\fontsize{#1}{#2pt}%
  \fontfamily{#3}\fontseries{#4}\fontshape{#5}%
  \selectfont}%
\fi\endgroup%
\begin{picture}(1824,1711)(830,-1539)
\put(1715,-1036){\makebox(0,0)[lb]{\smash{{\SetFigFont{7}{8.4}{\familydefault}{\mddefault}{\updefault}{\color[rgb]{0,0,0}$\pi/3$}%
}}}}
\put(2654,-710){\makebox(0,0)[lb]{\smash{{\SetFigFont{7}{8.4}{\familydefault}{\mddefault}{\updefault}{\color[rgb]{0,0,0}$\pi/2$}%
}}}}
\put(2445,-1223){\makebox(0,0)[lb]{\smash{{\SetFigFont{7}{8.4}{\familydefault}{\mddefault}{\updefault}{\color[rgb]{0,0,0}$\pi/2$}%
}}}}
\put(1172,-1356){\makebox(0,0)[lb]{\smash{{\SetFigFont{7}{8.4}{\familydefault}{\mddefault}{\updefault}{\color[rgb]{0,0,0}$\pi/2$}%
}}}}
\put(1039,-140){\makebox(0,0)[lb]{\smash{{\SetFigFont{7}{8.4}{\familydefault}{\mddefault}{\updefault}{\color[rgb]{0,0,0}$\pi/2$}%
}}}}
\put(1762, 88){\makebox(0,0)[lb]{\smash{{\SetFigFont{7}{8.4}{\familydefault}{\mddefault}{\updefault}{\color[rgb]{0,0,0}$\pi/2$}%
}}}}
\put(2464,-159){\makebox(0,0)[lb]{\smash{{\SetFigFont{7}{8.4}{\familydefault}{\mddefault}{\updefault}{\color[rgb]{0,0,0}$\pi/2$}%
}}}}
\put(1362,-729){\makebox(0,0)[lb]{\smash{{\SetFigFont{7}{8.4}{\familydefault}{\mddefault}{\updefault}{\color[rgb]{0,0,0}$\pi/3$}%
}}}}
\put(1457,-444){\makebox(0,0)[lb]{\smash{{\SetFigFont{7}{8.4}{\familydefault}{\mddefault}{\updefault}{\color[rgb]{0,0,0}$\pi/3$}%
}}}}
\put(2122,-710){\makebox(0,0)[lb]{\smash{{\SetFigFont{7}{8.4}{\familydefault}{\mddefault}{\updefault}{\color[rgb]{0,0,0}$\pi/3$}%
}}}}
\put(1457,-975){\makebox(0,0)[lb]{\smash{{\SetFigFont{7}{8.4}{\familydefault}{\mddefault}{\updefault}{\color[rgb]{0,0,0}$\pi/3$}%
}}}}
\put(1576,-109){\makebox(0,0)[lb]{\smash{{\SetFigFont{7}{8.4}{\familydefault}{\mddefault}{\updefault}{\color[rgb]{0,0,0}$\pi/2$}%
}}}}
\put(1593,-1256){\makebox(0,0)[lb]{\smash{{\SetFigFont{7}{8.4}{\familydefault}{\mddefault}{\updefault}{\color[rgb]{0,0,0}$\pi/2$}%
}}}}
\put(2122,-1211){\makebox(0,0)[lb]{\smash{{\SetFigFont{7}{8.4}{\familydefault}{\mddefault}{\updefault}{\color[rgb]{0,0,0}$\pi/2$}%
}}}}
\put(1172,-402){\makebox(0,0)[lb]{\smash{{\SetFigFont{7}{8.4}{\familydefault}{\mddefault}{\updefault}{\color[rgb]{0,0,0}$\pi/2$}%
}}}}
\put(830,-710){\makebox(0,0)[lb]{\smash{{\SetFigFont{7}{8.4}{\familydefault}{\mddefault}{\updefault}{\color[rgb]{0,0,0}$\pi/2$}%
}}}}
\put(1187,-1070){\makebox(0,0)[lb]{\smash{{\SetFigFont{7}{8.4}{\familydefault}{\mddefault}{\updefault}{\color[rgb]{0,0,0}$\pi/2$}%
}}}}
\put(2178,-187){\makebox(0,0)[lb]{\smash{{\SetFigFont{7}{8.4}{\familydefault}{\mddefault}{\updefault}{\color[rgb]{0,0,0}$\pi/2$}%
}}}}
\put(1731,-1508){\makebox(0,0)[lb]{\smash{{\SetFigFont{7}{8.4}{\familydefault}{\mddefault}{\updefault}{\color[rgb]{0,0,0}$\pi/4$}%
}}}}
\put(2024,-938){\makebox(0,0)[lb]{\smash{{\SetFigFont{7}{8.4}{\familydefault}{\mddefault}{\updefault}{\color[rgb]{0,0,0}$\pi/3$}%
}}}}
\put(1976,-413){\makebox(0,0)[lb]{\smash{{\SetFigFont{7}{8.4}{\familydefault}{\mddefault}{\updefault}{\color[rgb]{0,0,0}$\pi/3$}%
}}}}
\put(2374,-856){\makebox(0,0)[lb]{\smash{{\SetFigFont{7}{8.4}{\familydefault}{\mddefault}{\updefault}{\color[rgb]{0,0,0}$\pi/2$}%
}}}}
\put(2385,-562){\makebox(0,0)[lb]{\smash{{\SetFigFont{7}{8.4}{\familydefault}{\mddefault}{\updefault}{\color[rgb]{0,0,0}$\pi/2$}%
}}}}
\put(1740,-299){\makebox(0,0)[lb]{\smash{{\SetFigFont{7}{8.4}{\familydefault}{\mddefault}{\updefault}{\color[rgb]{0,0,0}$\pi/3$}%
}}}}
\end{picture}%

\end{center}
\caption{\hbox{ }}
\end{figure}

These angles satisfy Andreev's conditions (1) -- (5) so we can use Newton's method
to deform the prism constructed in the previous paragraph to have these angles.
Gluing this prism is to its mirror image, the edges labeled $\pi/2$ on the
outside disappear as edges, and the edges labeled on the outside by
$\pi/4$ glue together becoming an edge with dihedral angle $\pi/2.$ Hence,
we have constructed a polyhedron realizing $D_N$, assuming  $N > 7$.
Notice that when $N \leq 7$ the combinatorics of $D_N$ coincides with that of
$Pr_N$.
\Endproof

Assume that the two vertices incident at an edge $e$ are trivalent.  A
Whitehead move $Wh(e)$ on edge $e$ is given by the local change of the abstract
polyhedron described in the following diagram.  The Whitehead move in the dual
complex is dashed.  Often we will find it convenient to describe the Whitehead
move entirely in terms of the dual complex, in which case we write $Wh(f)$.

\begin{figure}[htp]
\begin{center} 
\begin{picture}(0,0)%
\includegraphics{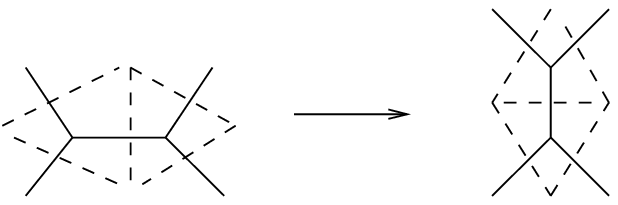}%
\end{picture}%
\setlength{\unitlength}{3947sp}%
\begingroup\makeatletter\ifx\SetFigFont\undefined%
\gdef\SetFigFont#1#2#3#4#5{%
  \reset@font\fontsize{#1}{#2pt}%
  \fontfamily{#3}\fontseries{#4}\fontshape{#5}%
  \selectfont}%
\fi\endgroup%
\begin{picture}(2975,1424)(1339,-1412)
\put(2157,-1381){\makebox(0,0)[lb]{\smash{\SetFigFont{8}{9.6}{\familydefault}{\mddefault}{\updefault}{\color[rgb]{0,0,0}Whitehead move on edge $e$}%
}}}
\put(3613,-72){\makebox(0,0)[lb]{\smash{\SetFigFont{8}{9.6}{\familydefault}{\mddefault}{\updefault}{\color[rgb]{0,0,0}$e_1$}%
}}}
\put(4290,-83){\makebox(0,0)[lb]{\smash{\SetFigFont{8}{9.6}{\familydefault}{\mddefault}{\updefault}{\color[rgb]{0,0,0}$e_2$}%
}}}
\put(4314,-1085){\makebox(0,0)[lb]{\smash{\SetFigFont{8}{9.6}{\familydefault}{\mddefault}{\updefault}{\color[rgb]{0,0,0}$e_3$}%
}}}
\put(3614,-1096){\makebox(0,0)[lb]{\smash{\SetFigFont{8}{9.6}{\familydefault}{\mddefault}{\updefault}{\color[rgb]{0,0,0}$e_4$}%
}}}
\put(4004,-687){\makebox(0,0)[lb]{\smash{\SetFigFont{8}{9.6}{\familydefault}{\mddefault}{\updefault}{\color[rgb]{0,0,0}$e'$}%
}}}
\put(1799,-671){\makebox(0,0)[lb]{\smash{\SetFigFont{8}{9.6}{\familydefault}{\mddefault}{\updefault}{\color[rgb]{0,0,0}$e$}%
}}}
\put(1374,-340){\makebox(0,0)[lb]{\smash{\SetFigFont{8}{9.6}{\familydefault}{\mddefault}{\updefault}{\color[rgb]{0,0,0}$e_1$}%
}}}
\put(2426,-351){\makebox(0,0)[lb]{\smash{\SetFigFont{8}{9.6}{\familydefault}{\mddefault}{\updefault}{\color[rgb]{0,0,0}$e_2$}%
}}}
\put(2483,-1051){\makebox(0,0)[lb]{\smash{\SetFigFont{8}{9.6}{\familydefault}{\mddefault}{\updefault}{\color[rgb]{0,0,0}$e_3$}%
}}}
\put(1357,-1096){\makebox(0,0)[lb]{\smash{\SetFigFont{8}{9.6}{\familydefault}{\mddefault}{\updefault}{\color[rgb]{0,0,0}$e_4$}%
}}}
\put(1992,-556){\makebox(0,0)[lb]{\smash{\SetFigFont{8}{9.6}{\familydefault}{\mddefault}{\updefault}{\color[rgb]{0,0,0}$f$}%
}}}
\put(3793,-497){\makebox(0,0)[lb]{\smash{\SetFigFont{8}{9.6}{\familydefault}{\mddefault}{\updefault}{\color[rgb]{0,0,0}$f'$}%
}}}
\put(2832,-476){\makebox(0,0)[lb]{\smash{\SetFigFont{8}{9.6}{\familydefault}{\mddefault}{\updefault}{\color[rgb]{0,0,0}$Wh(e)$}%
}}}
\end{picture}

\end{center}
\caption{\hbox{ }}
\end{figure}

The following lemma appears in \cite{ROE2}:

\begin{lem} \label{WHITEHEAD} Let the abstract polyhedron $C'$ be obtained
from the simple abstract polyhedron $C$ by a Whitehead move $Wh(e)$. 
Then if ${\cal P}^0_C$ is non-empty so is ${\cal P}^0_{C'}.$ \end{lem}

The proof constructs a sequence of polyhedra realizing $C$ with dihedral angles
chosen so that the edge $e$ converges to a single point at infinity.  A
carefully chosen small perturbation of this limiting configuration results in a
compact polyhedron realizing $C'$ with non-obtuse dihedral angles.  

Suppose that we have a polyhedron realizing $C$ with all dihedral angles equal
to $\frac{2\pi}{5}$, and  choose a small $\epsilon > 0$.  To implement a
Whitehead move using the computer, we assign dihedral angle $\epsilon$ to the
the edge $e$ and dihedral angle $\frac{\pi}{2}$ to the four edges sharing an
endpoint with $e$.  Leaving the dihedral angles of the remaining edges the
same, the resulting set of angles is in $A_C$ and hence we can use Newton's
Method to deform $P$ into a polyhedron $P_1$ realizing $C$ with these new
angles.

If $\epsilon$ was chosen small enough, $P_1$ will be in the basin of attraction
for a polyhedron realizing $C'$ with the edge $e'$ replacing $e$, the dihedral
angle at $e'$ equal to $\epsilon$, and all other dihedral angles as in $P_1$.
We call the resulting polyhedron $P_2$.  Since $C'$ is simple we can deform
$P_2$ to have all dihedral angles $\frac{2\pi}{5}$, obtaining $P'$.

The following diagram shows these four steps when doing a Whitehead move on
one of the edges of the dodecahedron.  Here and elsewhere in this paper we use $\epsilon$
approximately  $\frac{\pi}{45}$.  (A smaller $\epsilon$ may be necessary when constructing polyhedra
with a very large number of faces.)

\begin{figure}[htp]
\begin{center}
\begin{picture}(0,0)%
\epsfig{file=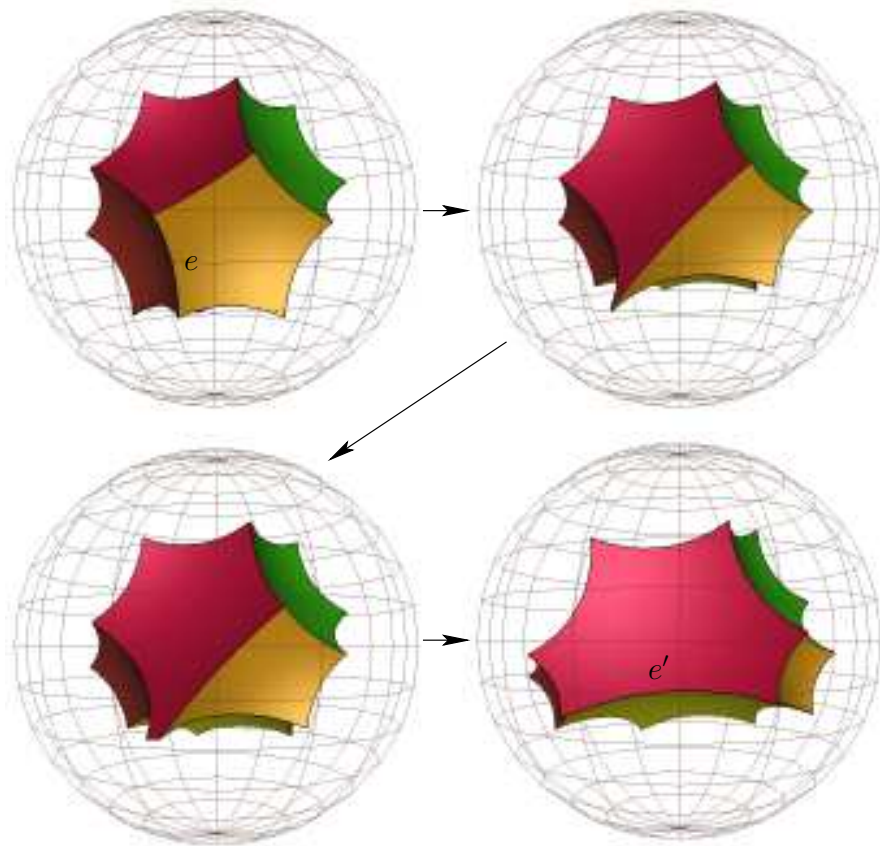}%
\end{picture}%
\setlength{\unitlength}{3947sp}%
\begingroup\makeatletter\ifx\SetFigFont\undefined%
\gdef\SetFigFont#1#2#3#4#5{%
  \reset@font\fontsize{#1}{#2pt}%
  \fontfamily{#3}\fontseries{#4}\fontshape{#5}%
  \selectfont}%
\fi\endgroup%
\begin{picture}(5550,5550)(1201,-5911)
\put(2326,-2086){\makebox(0,0)[lb]{\smash{{\SetFigFont{12}{14.4}{\familydefault}{\mddefault}{\updefault}{\color[rgb]{0,0,0}$e$}%
}}}}
\put(5239,-4653){\makebox(0,0)[lb]{\smash{{\SetFigFont{12}{14.4}{\familydefault}{\mddefault}{\updefault}{\color[rgb]{0,0,0}$e'$}%
}}}}
\end{picture}%

\end{center}
\caption{Whitehead move $Wh(e)$ on edge $e$ of the dodecahedron}
\end{figure}

\vspace{.1in}

If we can find a sequence of combinatorial Whitehead moves reducing a given
simple abstract polyhedron $C$ to either $Pr_N$ or $D_N$ via a sequence of
simple abstract polyhedra $C_1,\cdots,C_N$, we can use Newton's method to
perform this sequence of Whitehead moves in the reverse order, constructing
geometric polyhedra that realize $C_N, C_{N-1},\cdots,C_1$, and finally $C$.
Before explaining why such a sequence always exists, we demonstrate this
process for the dodecahedron.

\begin{figure}[htp]
\begin{center}

\scalebox{1.3}{
%
}
\end{center}
\caption{\hbox{ }}
\end{figure}

We can then use Newton's Method to do the reverse sequence of Whitehead moves $Wh(8,11)$, $Wh(4,11)$,
$Wh(1,2)$, $Wh(9,11)$, $Wh(2,4)$, $Wh(1,6)$, $Wh(7,11)$, $Wh(6,9)$, $Wh(1,5)$, and
$Wh(1,4)$ geometrically, constructing the dodecahedron
from $D_{12}$.  The following diagram shows this process.

\begin{figure}[htp]
\begin{center}
\includegraphics[scale=1.2]{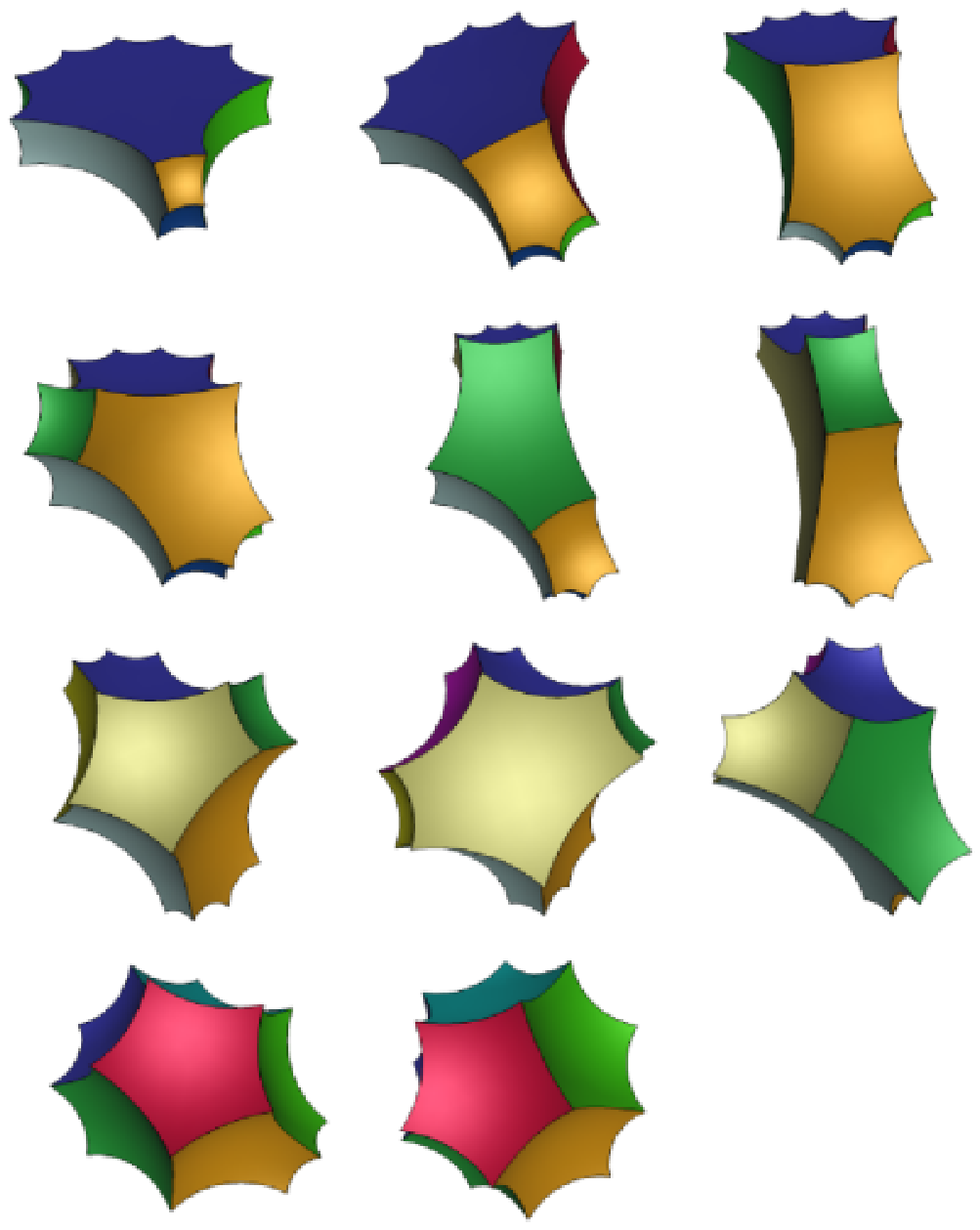}
\end{center}
\caption{Construction of the dodecahedron from $D_{12}$ using $10$ Whitehead Moves.}
\end{figure}

\subsection{A Lemma on Whitehead Moves}\label{SEC_ALGORITHM}

The following lemma from \cite{AND} and 
\cite{ROE2} is necessary to prove Andreev's Theorem and for our
construction of simple polyhedra doing Whitehead moves geometrically with Newton's Method. 

\begin{lem} {\bf Whitehead Sequence} \label{ROEDER}

Let $C$ be a simple abstract polyhedron on $\mathbb{S}^2$ which is
not a prism.  If $C$ has $N > 7$ faces, $C$ can be simplified to $D_N$ by a finite
sequence of Whitehead moves such that {\it all of the intermediate
abstract polyhedra are simple.}

\end{lem}

Theorem 6 in Andreev's original paper contains our Lemma \ref{ROEDER}.
Andreev's original proof of Theorem 6 provides an algorithm to produce the
Whitehead moves needed for this lemma, but the algorithm {\it contains a glitch
}.  An error was detected when the author tried to implement it
as a computer program.  The algorithm failed on the first test, when given a dodecahedron. 

In the instead using $Wh(6,9)$ for the fifth Whitehead move of the sequence
described in the previous section, Andreev's algorithm uses either $Wh(2,6)$ or
$Wh(2,5)$.  In both cases it produces an abstract polyhedron which had a
prismatic 3-circuit, see below:

\begin{figure}[htp]
\begin{center}
\scalebox{1.2}{
\begin{picture}(0,0)%
\epsfig{file=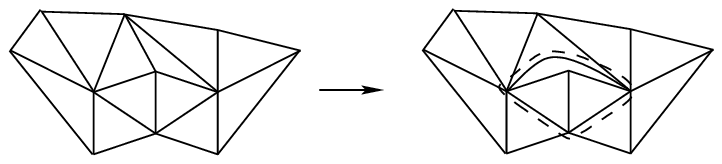}%
\end{picture}%
\setlength{\unitlength}{3947sp}%
\begingroup\makeatletter\ifx\SetFigFont\undefined%
\gdef\SetFigFont#1#2#3#4#5{%
  \reset@font\fontsize{#1}{#2pt}%
  \fontfamily{#3}\fontseries{#4}\fontshape{#5}%
  \selectfont}%
\fi\endgroup%
\begin{picture}(3556,879)(1074,-3976)
\put(3792,-3867){\makebox(0,0)[lb]{\smash{{\SetFigFont{5}{6.0}{\familydefault}{\mddefault}{\updefault}{\color[rgb]{0,0,0}$5$}%
}}}}
\put(1568,-3954){\makebox(0,0)[lb]{\smash{{\SetFigFont{5}{6.0}{\familydefault}{\mddefault}{\updefault}{\color[rgb]{0,0,0}$10$}%
}}}}
\put(2525,-3355){\makebox(0,0)[lb]{\smash{{\SetFigFont{5}{6.0}{\familydefault}{\mddefault}{\updefault}{\color[rgb]{0,0,0}$12$}%
}}}}
\put(2169,-3518){\makebox(0,0)[lb]{\smash{{\SetFigFont{5}{6.0}{\familydefault}{\mddefault}{\updefault}{\color[rgb]{0,0,0}$9$}%
}}}}
\put(2223,-3912){\makebox(0,0)[lb]{\smash{{\SetFigFont{5}{6.0}{\familydefault}{\mddefault}{\updefault}{\color[rgb]{0,0,0}$8$}%
}}}}
\put(1545,-3178){\makebox(0,0)[lb]{\smash{{\SetFigFont{5}{6.0}{\familydefault}{\mddefault}{\updefault}{\color[rgb]{0,0,0}$6$}%
}}}}
\put(2104,-3248){\makebox(0,0)[lb]{\smash{{\SetFigFont{5}{6.0}{\familydefault}{\mddefault}{\updefault}{\color[rgb]{0,0,0}$4$}%
}}}}
\put(1247,-3161){\makebox(0,0)[lb]{\smash{{\SetFigFont{5}{6.0}{\familydefault}{\mddefault}{\updefault}{\color[rgb]{0,0,0}$7$}%
}}}}
\put(1074,-3351){\makebox(0,0)[lb]{\smash{{\SetFigFont{5}{6.0}{\familydefault}{\mddefault}{\updefault}{\color[rgb]{0,0,0}$3$}%
}}}}
\put(1393,-3642){\makebox(0,0)[lb]{\smash{{\SetFigFont{5}{6.0}{\familydefault}{\mddefault}{\updefault}{\color[rgb]{0,0,0}$11$}%
}}}}
\put(1865,-3582){\makebox(0,0)[lb]{\smash{{\SetFigFont{5}{6.0}{\familydefault}{\mddefault}{\updefault}{\color[rgb]{0,0,0}$2$}%
}}}}
\put(1810,-3871){\makebox(0,0)[lb]{\smash{{\SetFigFont{5}{6.0}{\familydefault}{\mddefault}{\updefault}{\color[rgb]{0,0,0}$5$}%
}}}}
\put(2616,-3675){\makebox(0,0)[lb]{\smash{{\SetFigFont{5}{6.0}{\familydefault}{\mddefault}{\updefault}{\color[rgb]{0,0,0}{\tiny $Wh(2,6)$}}%
}}}}
\put(3550,-3950){\makebox(0,0)[lb]{\smash{{\SetFigFont{5}{6.0}{\familydefault}{\mddefault}{\updefault}{\color[rgb]{0,0,0}$10$}%
}}}}
\put(4507,-3351){\makebox(0,0)[lb]{\smash{{\SetFigFont{5}{6.0}{\familydefault}{\mddefault}{\updefault}{\color[rgb]{0,0,0}$12$}%
}}}}
\put(4151,-3514){\makebox(0,0)[lb]{\smash{{\SetFigFont{5}{6.0}{\familydefault}{\mddefault}{\updefault}{\color[rgb]{0,0,0}$9$}%
}}}}
\put(4205,-3908){\makebox(0,0)[lb]{\smash{{\SetFigFont{5}{6.0}{\familydefault}{\mddefault}{\updefault}{\color[rgb]{0,0,0}$8$}%
}}}}
\put(3527,-3174){\makebox(0,0)[lb]{\smash{{\SetFigFont{5}{6.0}{\familydefault}{\mddefault}{\updefault}{\color[rgb]{0,0,0}$6$}%
}}}}
\put(4086,-3244){\makebox(0,0)[lb]{\smash{{\SetFigFont{5}{6.0}{\familydefault}{\mddefault}{\updefault}{\color[rgb]{0,0,0}$4$}%
}}}}
\put(3229,-3157){\makebox(0,0)[lb]{\smash{{\SetFigFont{5}{6.0}{\familydefault}{\mddefault}{\updefault}{\color[rgb]{0,0,0}$7$}%
}}}}
\put(3056,-3347){\makebox(0,0)[lb]{\smash{{\SetFigFont{5}{6.0}{\familydefault}{\mddefault}{\updefault}{\color[rgb]{0,0,0}$3$}%
}}}}
\put(3375,-3638){\makebox(0,0)[lb]{\smash{{\SetFigFont{5}{6.0}{\familydefault}{\mddefault}{\updefault}{\color[rgb]{0,0,0}$11$}%
}}}}
\put(3847,-3578){\makebox(0,0)[lb]{\smash{{\SetFigFont{5}{6.0}{\familydefault}{\mddefault}{\updefault}{\color[rgb]{0,0,0}$2$}%
}}}}
\end{picture}%
}
\end{center}
\label{FIG_ERROR}
\caption{\hbox{  }}
\end{figure}

\vspace{.1in} In combination with the computer-implemented Whitehead move
described in the previous section, the sequence of Whitehead moves given in the
proof of Lemma \ref{ROEDER} gives us the path that we will use for our homotopy
method when constructing simple polyhedra.  We provide an outline of the proof
here that is sufficient to describe the sequence of Whitehead moves.  Those who
wish to see a complete proof may refer to \cite{ROE2,ROE}.

\noindent {\bf Outline of the Proof of Lemma \ref{ROEDER}:} 

We assume that $C \neq Pr_N$ is a simple abstract polyhedron with $N>7$ faces.
We will construct a sequence of Whitehead moves that change $C$ to $D_N$, so
that no intermediate complex has a prismatic 3-circuit.

Find a vertex $v_\infty$ of $C^*$ which is connected to the greatest
number of other vertices.  We will call the link of $v_\infty$, a cycle of $k$
vertices and $k$ edges, the {\it outer-polygon}.  Most of the work is to show
that we can do Whitehead moves to increase $k$ to $N-3$ without
introducing any prismatic 3-circuits during the process.  Once this is
completed, it is easy to change the resulting complex to $D^*_N$ by
additional Whitehead moves.

Let us set up some notation.  Draw the dual complex $C^*$ in the
plane with the vertex $v_\infty$ at infinity, and the outer polygon $P$
surrounding the remaining vertices and triangles.  We call the
vertices inside of $P$ {\it interior vertices}.  All of the
edges inside of $P$ which do not have an endpoint on $P$ are called {\it interior
edges}.  

Note that the graph of interior vertices and edges is connected, since
$C^*$ is simple.  An interior vertex which is connected to only one
other interior vertex will be called an {\it endpoint}.

\begin{center}
\begin{picture}(0,0)%
\epsfig{file=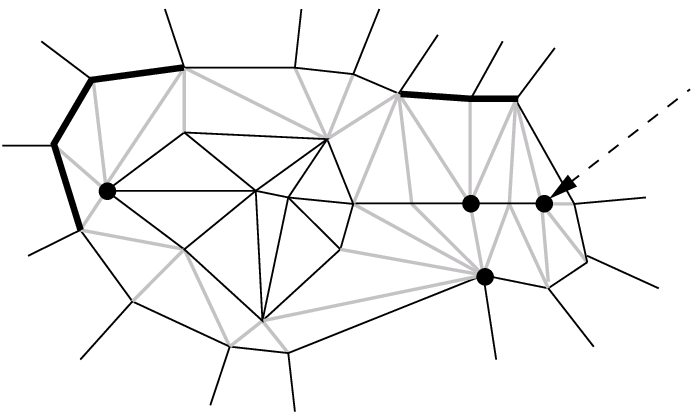}%
\end{picture}%
\setlength{\unitlength}{3947sp}%
\begingroup\makeatletter\ifx\SetFigFont\undefined%
\gdef\SetFigFont#1#2#3#4#5{%
  \reset@font\fontsize{#1}{#2pt}%
  \fontfamily{#3}\fontseries{#4}\fontshape{#5}%
  \selectfont}%
\fi\endgroup%
\begin{picture}(3387,1956)(1189,-2220)
\put(3578,-1771){\makebox(0,0)[lb]{\smash{{\SetFigFont{12}{14.4}{\familydefault}{\mddefault}{\updefault}{\color[rgb]{0,0,0}\small{$F^2_v$}}%
}}}}
\put(1223,-751){\makebox(0,0)[lb]{\smash{{\SetFigFont{12}{14.4}{\familydefault}{\mddefault}{\updefault}{\color[rgb]{0,0,0}\small{$F^1_w$}}%
}}}}
\put(1651,-1329){\makebox(0,0)[lb]{\smash{{\SetFigFont{12}{14.4}{\familydefault}{\mddefault}{\updefault}{\color[rgb]{0,0,0}\small{$w$}}%
}}}}
\put(3270,-638){\makebox(0,0)[lb]{\smash{{\SetFigFont{12}{14.4}{\familydefault}{\mddefault}{\updefault}{\color[rgb]{0,0,0}\small{$F^1_v$}}%
}}}}
\put(3256,-1179){\makebox(0,0)[lb]{\smash{{\SetFigFont{12}{14.4}{\familydefault}{\mddefault}{\updefault}{\color[rgb]{0,0,0}\small{$v$}}%
}}}}
\put(4576,-661){\makebox(0,0)[lb]{\smash{{\SetFigFont{12}{14.4}{\familydefault}{\mddefault}{\updefault}{\color[rgb]{0,0,0}\small{endpoint}}%
}}}}
\end{picture}%

\end{center}

Throughout this proof we will draw $P$ and the interior edges and
vertices of $C^*$ in black.  The connections between $P$ and the
interior vertices will be grey.  Connections between $P$ and $v_{\infty}$ will
be black, if shown at all.

The link of an interior vertex $v$ intersects $P$ in a
number of components $F_v^1,\cdots,F_v^n$ (possibly $n = 0$.)  See the
above figure.  We say that $v$ is {\it connected to $P$ in these
components.}  Notice that since $C^*$ is simple, an endpoint is
always connected to $P$ in exactly one such component.

\begin{move} \label{SUBLEMMA1}
Suppose that there is an interior vertex $A$ of $C^*$ which is connected to $P$ in exactly
one component consisting of exactly two consecutive vertices $Q$ and $R$.
The Whitehead move $Wh(QR)$ on $C^*$ increases the length of the
outer polygon by one, and introduces no prismatic $3$-circuit.
\end{move}

\begin{center}
\begin{picture}(0,0)%
\epsfig{file=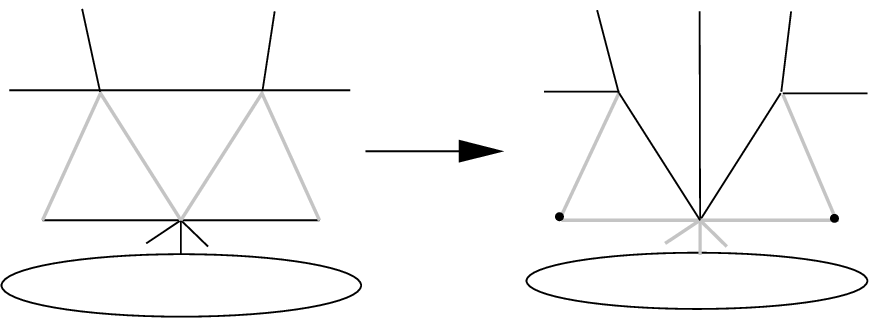}%
\end{picture}%
\setlength{\unitlength}{3947sp}%
\begingroup\makeatletter\ifx\SetFigFont\undefined%
\gdef\SetFigFont#1#2#3#4#5{%
  \reset@font\fontsize{#1}{#2pt}%
  \fontfamily{#3}\fontseries{#4}\fontshape{#5}%
  \selectfont}%
\fi\endgroup%
\begin{picture}(4176,1741)(1537,-2354)
\put(2316,-835){\makebox(0,0)[lb]{\smash{{\SetFigFont{12}{14.4}{\familydefault}{\mddefault}{\updefault}{\color[rgb]{0,0,0}\small{$v_\infty$}}%
}}}}
\put(5542,-2022){\makebox(0,0)[lb]{\smash{{\SetFigFont{8}{9.6}{\familydefault}{\mddefault}{\updefault}{\color[rgb]{0,0,0}\small{$E$}}%
}}}}
\put(4755,-769){\makebox(0,0)[lb]{\smash{{\SetFigFont{12}{14.4}{\familydefault}{\mddefault}{\updefault}{\color[rgb]{0,0,0}\small{$v_\infty$}}%
}}}}
\put(3402,-1387){\makebox(0,0)[lb]{\smash{{\SetFigFont{8}{9.6}{\familydefault}{\mddefault}{\updefault}{\color[rgb]{0,0,0}$Wh(QR)$}%
}}}}
\put(2828,-1223){\makebox(0,0)[lb]{\smash{{\SetFigFont{8}{9.6}{\familydefault}{\mddefault}{\updefault}{\color[rgb]{0,0,0}\small{$R$}}%
}}}}
\put(5315,-1216){\makebox(0,0)[lb]{\smash{{\SetFigFont{8}{9.6}{\familydefault}{\mddefault}{\updefault}{\color[rgb]{0,0,0}\small{$R$}}%
}}}}
\put(4509,-2236){\makebox(0,0)[lb]{\smash{{\SetFigFont{8}{9.6}{\familydefault}{\mddefault}{\updefault}{\color[rgb]{0,0,0}\small{interior stuff}}%
}}}}
\put(1704,-2243){\makebox(0,0)[lb]{\smash{{\SetFigFont{8}{9.6}{\familydefault}{\mddefault}{\updefault}{\color[rgb]{0,0,0}\small{other interior stuff}}%
}}}}
\put(2546,-2018){\makebox(0,0)[lb]{\smash{{\SetFigFont{8}{9.6}{\familydefault}{\mddefault}{\updefault}{\color[rgb]{0,0,0}\small{$A$}}%
}}}}
\put(3106,-2041){\makebox(0,0)[lb]{\smash{{\SetFigFont{8}{9.6}{\familydefault}{\mddefault}{\updefault}{\color[rgb]{0,0,0}\small{$E$}}%
}}}}
\put(1621,-2033){\makebox(0,0)[lb]{\smash{{\SetFigFont{8}{9.6}{\familydefault}{\mddefault}{\updefault}{\color[rgb]{0,0,0}\small{$D$}}%
}}}}
\put(1789,-1207){\makebox(0,0)[lb]{\smash{{\SetFigFont{8}{9.6}{\familydefault}{\mddefault}{\updefault}{\color[rgb]{0,0,0}\small{$Q$}}%
}}}}
\put(4267,-1214){\makebox(0,0)[lb]{\smash{{\SetFigFont{8}{9.6}{\familydefault}{\mddefault}{\updefault}{\color[rgb]{0,0,0}\small{$Q$}}%
}}}}
\put(5057,-2019){\makebox(0,0)[lb]{\smash{{\SetFigFont{8}{9.6}{\familydefault}{\mddefault}{\updefault}{\color[rgb]{0,0,0}\small{$A$}}%
}}}}
\put(4114,-2034){\makebox(0,0)[lb]{\smash{{\SetFigFont{8}{9.6}{\familydefault}{\mddefault}{\updefault}{\color[rgb]{0,0,0}\small{$D$}}%
}}}}
\end{picture}%

\end{center}
\vspace{.05in}

\begin{move} \label{SUBLEMMA2}
Suppose that there is an interior vertex $A$ that is connected to $P$ in a
component consisting of $M$ consecutive vertices $Q_1,\cdots,Q_M$ of $P$ (and
possibly other components).

\noindent (a)
If $A$ is not an endpoint and $M > 2$, the sequence of Whitehead moves
$Wh(AQ_M),\ldots,Wh(AQ_3)$ results in a complex in which $A$ is connected
to the same component of $P$ in only $Q_1$ and $Q_2$.  These moves leave $P$
unchanged, and introduce no prismatic 3-circuit.

\vspace{.07in}
\begin{center}
\begin{picture}(0,0)%
\epsfig{file=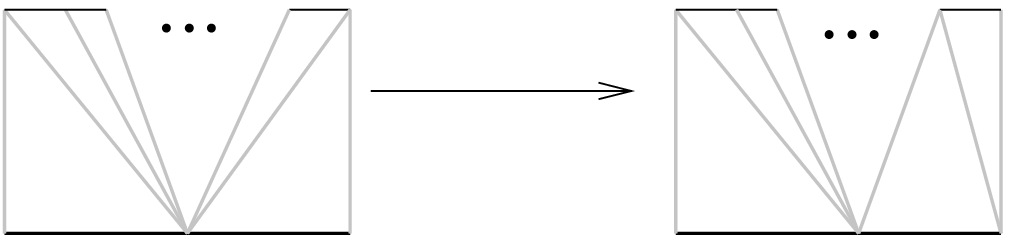}%
\end{picture}%
\setlength{\unitlength}{3947sp}%
\begingroup\makeatletter\ifx\SetFigFont\undefined%
\gdef\SetFigFont#1#2#3#4#5{%
  \reset@font\fontsize{#1}{#2pt}%
  \fontfamily{#3}\fontseries{#4}\fontshape{#5}%
  \selectfont}%
\fi\endgroup%
\begin{picture}(5193,1438)(1104,-1372)
\put(1902,-1341){\makebox(0,0)[lb]{\smash{{\SetFigFont{9}{10.8}{\familydefault}{\mddefault}{\updefault}{\color[rgb]{0,0,0}$A$}%
}}}}
\put(2319,-30){\makebox(0,0)[lb]{\smash{{\SetFigFont{9}{10.8}{\familydefault}{\mddefault}{\updefault}{\color[rgb]{0,0,0}$Q_{M-1}$}%
}}}}
\put(2797,-30){\makebox(0,0)[lb]{\smash{{\SetFigFont{9}{10.8}{\familydefault}{\mddefault}{\updefault}{\color[rgb]{0,0,0}$Q_M$}%
}}}}
\put(4348,-1341){\makebox(0,0)[lb]{\smash{{\SetFigFont{9}{10.8}{\familydefault}{\mddefault}{\updefault}{\color[rgb]{0,0,0}$D$}%
}}}}
\put(5959,-1341){\makebox(0,0)[lb]{\smash{{\SetFigFont{9}{10.8}{\familydefault}{\mddefault}{\updefault}{\color[rgb]{0,0,0}$E$}%
}}}}
\put(5124,-1341){\makebox(0,0)[lb]{\smash{{\SetFigFont{9}{10.8}{\familydefault}{\mddefault}{\updefault}{\color[rgb]{0,0,0}$A$}%
}}}}
\put(4348,-30){\makebox(0,0)[lb]{\smash{{\SetFigFont{9}{10.8}{\familydefault}{\mddefault}{\updefault}{\color[rgb]{0,0,0}$Q_1 Q_2 Q_3$}%
}}}}
\put(5481,-30){\makebox(0,0)[lb]{\smash{{\SetFigFont{9}{10.8}{\familydefault}{\mddefault}{\updefault}{\color[rgb]{0,0,0}$Q_{M-1}$}%
}}}}
\put(5898,-30){\makebox(0,0)[lb]{\smash{{\SetFigFont{9}{10.8}{\familydefault}{\mddefault}{\updefault}{\color[rgb]{0,0,0}$Q_M$}%
}}}}
\put(3035,-387){\makebox(0,0)[lb]{\smash{{\SetFigFont{9}{10.8}{\familydefault}{\mddefault}{\updefault}{\color[rgb]{0,0,0}$Wh(AQ_M)$}%
}}}}
\put(1126,-1341){\makebox(0,0)[lb]{\smash{{\SetFigFont{9}{10.8}{\familydefault}{\mddefault}{\updefault}{\color[rgb]{0,0,0}$D$}%
}}}}
\put(2737,-1341){\makebox(0,0)[lb]{\smash{{\SetFigFont{9}{10.8}{\familydefault}{\mddefault}{\updefault}{\color[rgb]{0,0,0}$E$}%
}}}}
\put(1126,-30){\makebox(0,0)[lb]{\smash{{\SetFigFont{9}{10.8}{\familydefault}{\mddefault}{\updefault}{\color[rgb]{0,0,0}$Q_1 Q_2 Q_3$}%
}}}}
\end{picture}%

\end{center}
\vspace{.07in}

\noindent (b)
If $A$ is an endpoint and $M > 3$, the sequence of Whitehead moves
\newline $Wh(AQ_M),\ldots,Wh(AQ_4)$ results in a complex in which $A$ is connected
to the same component of $P$ in only $Q_1,Q_2$, and $Q_3$.  These moves leave $P$
unchanged and introduce no prismatic 3-circuits.

\vspace{.07in}
\begin{center}
\begin{picture}(0,0)%
\epsfig{file=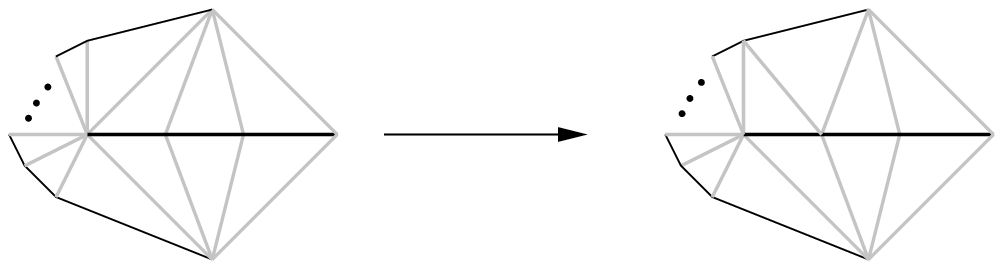}%
\end{picture}%
\setlength{\unitlength}{3947sp}%
\begingroup\makeatletter\ifx\SetFigFont\undefined%
\gdef\SetFigFont#1#2#3#4#5{%
  \reset@font\fontsize{#1}{#2pt}%
  \fontfamily{#3}\fontseries{#4}\fontshape{#5}%
  \selectfont}%
\fi\endgroup%
\begin{picture}(5002,1664)(571,-1179)
\put(4918,-1121){\makebox(0,0)[lb]{\smash{{\SetFigFont{12}{14.4}{\familydefault}{\mddefault}{\updefault}{\color[rgb]{0,0,0}\small{$Q_1$}}%
}}}}
\put(4801,-511){\makebox(0,0)[lb]{\smash{{\SetFigFont{12}{14.4}{\familydefault}{\mddefault}{\updefault}{\color[rgb]{0,0,0}\small{$E$}}%
}}}}
\put(4471,-495){\makebox(0,0)[lb]{\smash{{\SetFigFont{12}{14.4}{\familydefault}{\mddefault}{\updefault}{\color[rgb]{0,0,0}\small{$A$}}%
}}}}
\put(4171,171){\makebox(0,0)[lb]{\smash{{\SetFigFont{12}{14.4}{\familydefault}{\mddefault}{\updefault}{\color[rgb]{0,0,0}\small{$Q_{M-1}$}}%
}}}}
\put(4884,329){\makebox(0,0)[lb]{\smash{{\SetFigFont{12}{14.4}{\familydefault}{\mddefault}{\updefault}{\color[rgb]{0,0,0}\small{$Q_M$}}%
}}}}
\put(4051,-811){\makebox(0,0)[lb]{\smash{{\SetFigFont{12}{14.4}{\familydefault}{\mddefault}{\updefault}{\color[rgb]{0,0,0}\small{$Q_2$}}%
}}}}
\put(3721,-422){\makebox(0,0)[lb]{\smash{{\SetFigFont{12}{14.4}{\familydefault}{\mddefault}{\updefault}{\color[rgb]{0,0,0}\small{$Q_4$}}%
}}}}
\put(3796,-615){\makebox(0,0)[lb]{\smash{{\SetFigFont{12}{14.4}{\familydefault}{\mddefault}{\updefault}{\color[rgb]{0,0,0}\small{$Q_3$}}%
}}}}
\put(3751, 14){\makebox(0,0)[lb]{\smash{{\SetFigFont{12}{14.4}{\familydefault}{\mddefault}{\updefault}{\color[rgb]{0,0,0}\small{$Q_{M-2}$}}%
}}}}
\put(2626,-286){\makebox(0,0)[lb]{\smash{{\SetFigFont{12}{14.4}{\familydefault}{\mddefault}{\updefault}{\color[rgb]{0,0,0}\small{$Wh(AQ_M)$}}%
}}}}
\put(1651,-511){\makebox(0,0)[lb]{\smash{{\SetFigFont{12}{14.4}{\familydefault}{\mddefault}{\updefault}{\color[rgb]{0,0,0}\small{$E$}}%
}}}}
\put(1321,-495){\makebox(0,0)[lb]{\smash{{\SetFigFont{12}{14.4}{\familydefault}{\mddefault}{\updefault}{\color[rgb]{0,0,0}\small{$A$}}%
}}}}
\put(901,-811){\makebox(0,0)[lb]{\smash{{\SetFigFont{12}{14.4}{\familydefault}{\mddefault}{\updefault}{\color[rgb]{0,0,0}\small{$Q_2$}}%
}}}}
\put(571,-422){\makebox(0,0)[lb]{\smash{{\SetFigFont{12}{14.4}{\familydefault}{\mddefault}{\updefault}{\color[rgb]{0,0,0}\small{$Q_4$}}%
}}}}
\put(1021,171){\makebox(0,0)[lb]{\smash{{\SetFigFont{12}{14.4}{\familydefault}{\mddefault}{\updefault}{\color[rgb]{0,0,0}\small{$Q_{M-1}$}}%
}}}}
\put(1734,329){\makebox(0,0)[lb]{\smash{{\SetFigFont{12}{14.4}{\familydefault}{\mddefault}{\updefault}{\color[rgb]{0,0,0}\small{$Q_M$}}%
}}}}
\put(646,-615){\makebox(0,0)[lb]{\smash{{\SetFigFont{12}{14.4}{\familydefault}{\mddefault}{\updefault}{\color[rgb]{0,0,0}\small{$Q_3$}}%
}}}}
\put(601, 14){\makebox(0,0)[lb]{\smash{{\SetFigFont{12}{14.4}{\familydefault}{\mddefault}{\updefault}{\color[rgb]{0,0,0}\small{$Q_{M-2}$}}%
}}}}
\put(1767,-1112){\makebox(0,0)[lb]{\smash{{\SetFigFont{12}{14.4}{\familydefault}{\mddefault}{\updefault}{\color[rgb]{0,0,0}\small{$Q_1$}}%
}}}}
\end{picture}%

\end{center}
\vspace{.07in}
\end{move}

Note: In both parts (1) and (2), each of the Whitehead moves $Wh(AQ_M)$
transfers the connection between $A$ and $Q_M$ to a connection between the
neighboring interior vertex $E$ and $Q_{M-1}$.  This is helpful in
case 2 later.

\begin{move}\label{SUBLEMMA3}
Suppose that there is an interior vertex $A$ whose link contains two distinct
vertices $X$ and $Y$ of $P$.  Then there are Whitehead moves which eliminate
any component in which $A$ is connected to $P$, if that component does not
contain $X$ or $Y$.  $P$ is unchanged, and no prismatic $3$-circuits will be
introduced.
\end{move}

\begin{center}
\vspace{.125in}
\begin{center}
\begin{picture}(0,0)%
\includegraphics{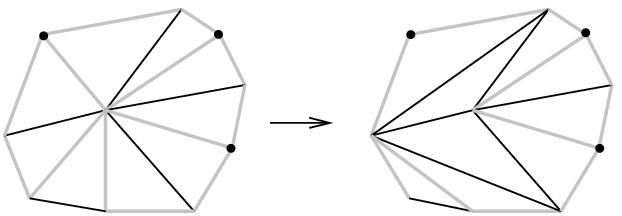}%
\end{picture}%
\setlength{\unitlength}{3947sp}%
\begingroup\makeatletter\ifx\SetFigFont\undefined%
\gdef\SetFigFont#1#2#3#4#5{%
  \reset@font\fontsize{#1}{#2pt}%
  \fontfamily{#3}\fontseries{#4}\fontshape{#5}%
  \selectfont}%
\fi\endgroup%
\begin{picture}(2959,1013)(879,-1352)
\put(1126,-811){\makebox(0,0)[lb]{\smash{\SetFigFont{12}{14.4}{\familydefault}{\mddefault}{\updefault}{\color[rgb]{0,0,0}\small{$A$}}%
}}}
\put(2994,-849){\makebox(0,0)[lb]{\smash{\SetFigFont{12}{14.4}{\familydefault}{\mddefault}{\updefault}{\color[rgb]{0,0,0}\small{$A$}}%
}}}
\put(2026,-511){\makebox(0,0)[lb]{\smash{\SetFigFont{6}{7.2}{\familydefault}{\mddefault}{\updefault}{\color[rgb]{0,0,0}\small{$X$}}%
}}}
\put(3826,-511){\makebox(0,0)[lb]{\smash{\SetFigFont{6}{7.2}{\familydefault}{\mddefault}{\updefault}{\color[rgb]{0,0,0}\small{$X$}}%
}}}
\put(2044,-1089){\makebox(0,0)[lb]{\smash{\SetFigFont{6}{7.2}{\familydefault}{\mddefault}{\updefault}{\color[rgb]{0,0,0}\small{$Y$}}%
}}}
\put(3809,-1094){\makebox(0,0)[lb]{\smash{\SetFigFont{6}{7.2}{\familydefault}{\mddefault}{\updefault}{\color[rgb]{0,0,0}\small{$Y$}}%
}}}
\end{picture}

\end{center}
\vspace{.125in}
\end{center}

Here $A$ is connected to $P$ in four components containing six vertices.  We can
eliminate connections of $A$ to all of the components except for the
single-point components $X$ and $Y$.

The proof that this move does not introduce any new prismatic 3-circuit is rather technical and depends essentially on the fact
that $A$ is connected to $P$ in at least two other vertices $X$ and $Y$.
Andreev describes a nearly identical process to move \ref{SUBLEMMA3} in his
paper \cite{AND} on pages 333-334.  However, he merely assumes that $A$ is
connected to $P$ in at least one component in addition to the components being
eliminated.  He does not require that $A$ is connected to $P$ in
at least {\it two vertices} outside of the components being eliminated. Andreev
then asserts: ``It is readily seen that all of the polyhedra obtained in this way
are simple...'' In fact, the Whitehead move
demonstrated in Figure \ref{FIG_ERROR} creates
a prismatic 3-circuit.

Having assumed this stronger (and incorrect) version of move 3, the remainder
of Andreev's proof is relatively easy.  Unfortunately, the situation pictured
in Figure \ref{FIG_ERROR} is not uncommon (as we will see in Case 3 below!)
Restricted to the weaker hypotheses of Move \ref{SUBLEMMA3} we will have to
work a little bit harder.

Using Move 1, Move 2, and Move 3 we check that if the length of $P$ is less
than $N-3$, then there is a sequence of Whitehead moves that increases the
length of $P$ by one without introducing any prismatic $3$-circuits.

\vspace{.1in}
{\bf Case 1:} An interior vertex that is not an endpoint connects to $P$ in
a component with two or more vertices, and possibly in other components.
\vspace{.05in}

Apply Move \ref{SUBLEMMA2} decreasing this component to two vertices.  We can
then apply Move \ref{SUBLEMMA3}, eliminating any other components since this component
contains two vertices. Finally, apply Move \ref{SUBLEMMA1} to increase the length of
the outer polygon by 1.

\vspace{.1in}
{\bf Case 2:}
An interior vertex that is an endpoint is connected to more than three
vertices of $P$.
\vspace{.05in}

We assume that each of the interior vertices that are not endpoints are
connected to $P$ in components consisting of single vertices, otherwise we
are in Case 1.

Let $A$ be the endpoint which is connected to more than three vertices of $P$.
By Move \ref{SUBLEMMA2}, part (2), there is a Whitehead move that transfers
one of these connections to the interior vertex $E$ that is next to $A$.  Now,
one of the components in which $E$ is connected to $P$ has exactly two
vertices.  The vertex $E$ is not an endpoint since $k < N-3$ implies that there
are at least three interior vertices.  Once this is done, we can apply Case 1.

\vspace{.1in}
{\bf Case 3:}  Each interior vertex that is an endpoint is connected to
exactly three vertices of $P$ and each interior vertex which is not an endpoint
is connected to $P$ in components each consisting of a single vertex.

\vspace{.05in}
First, notice that if the interior vertices and edges form a line, the
restriction on how interior vertices are connected to $P$ results in the prism, 
contrary to the assumption that $C$ is not the prism.  However, there are
many complexes satisfying the hypotheses of this case which have interior
vertices and edges forming a graph more complicated than a line:

\vspace{.125in}
\begin{center}
\begin{picture}(0,0)%
\includegraphics{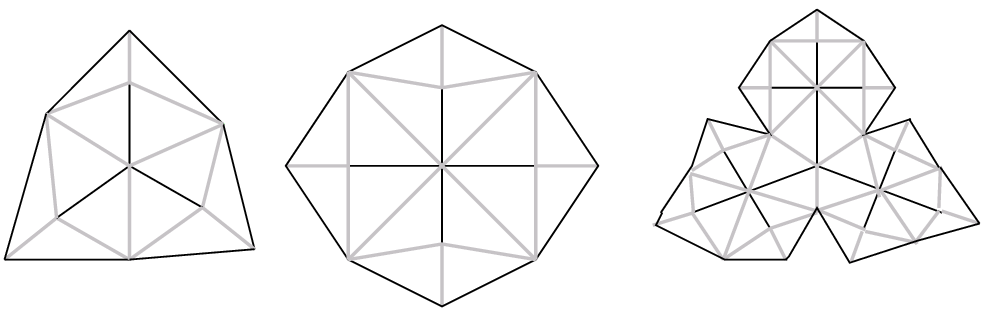}%
\end{picture}%
\setlength{\unitlength}{3947sp}%
\begingroup\makeatletter\ifx\SetFigFont\undefined%
\gdef\SetFigFont#1#2#3#4#5{%
  \reset@font\fontsize{#1}{#2pt}%
  \fontfamily{#3}\fontseries{#4}\fontshape{#5}%
  \selectfont}%
\fi\endgroup%
\begin{picture}(4724,1469)(879,-1283)
\end{picture}

\end{center}
\vspace{.125in}

For such complexes we need a very special sequence of Whitehead moves to
increase the length of $P$.

Pick an interior vertex which is an endpoint and label it $I_1$.  Denote by
$P_1$, $P_2$, and $P_3$
the three vertices of $P$ to which $I_1$ connects.  $I_1$
will be connected to a sequence of interior vertices $I_2, I_3, \cdots
I_m, m \ge 2$, with $I_m$ the first interior vertex in the sequence that is connected to
more than two other interior vertices.  Vertex $I_m$ must exist by the assumption
that the interior vertices don't form a line segment, the configuration that we
ruled out above.  By hypothesis, $I_2,\cdots,I_m$ can only connect to $P$ in
components which each consist of a vertex, hence each must be connected to $P_1$
and to $P_3$.  Similarly, there is an interior vertex (call it $X$) which
connects both to $I_m$ and to $P_1$ and another vertex $Y$ which connects to
$I_m$ and $P_3$.
Vertex $I_m$ may connect to other vertices of $P$ and other interior vertices, as
shown on the left side of the following diagram.

\vspace{.125in} \begin{center} 
\begin{picture}(0,0)%
\epsfig{file=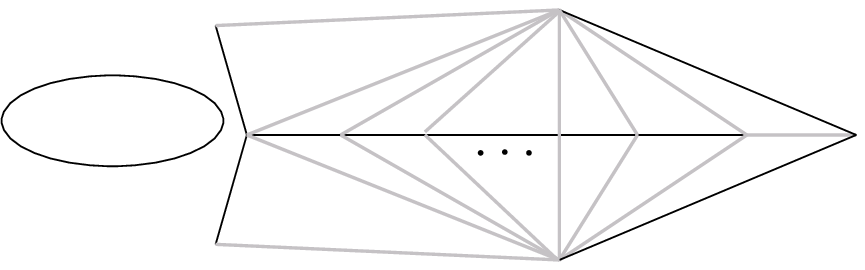}%
\end{picture}%
\setlength{\unitlength}{3947sp}%
\begingroup\makeatletter\ifx\SetFigFont\undefined%
\gdef\SetFigFont#1#2#3#4#5{%
  \reset@font\fontsize{#1}{#2pt}%
  \fontfamily{#3}\fontseries{#4}\fontshape{#5}%
  \selectfont}%
\fi\endgroup%
\begin{picture}(4185,1639)(16,-1844)
\put(1861,-915){\makebox(0,0)[lb]{\smash{{\SetFigFont{5}{6.0}{\familydefault}{\mddefault}{\updefault}{\color[rgb]{0,0,0}\small{$I_{m-1}$}}%
}}}}
\put(2851,-361){\makebox(0,0)[lb]{\smash{{\SetFigFont{12}{14.4}{\familydefault}{\mddefault}{\updefault}{\color[rgb]{0,0,0}\small{$P_1$}}%
}}}}
\put(4201,-961){\makebox(0,0)[lb]{\smash{{\SetFigFont{12}{14.4}{\familydefault}{\mddefault}{\updefault}{\color[rgb]{0,0,0}\small{$P_2$}}%
}}}}
\put(2701,-1786){\makebox(0,0)[lb]{\smash{{\SetFigFont{12}{14.4}{\familydefault}{\mddefault}{\updefault}{\color[rgb]{0,0,0}\small{$P_3$}}%
}}}}
\put(1051,-361){\makebox(0,0)[lb]{\smash{{\SetFigFont{12}{14.4}{\familydefault}{\mddefault}{\updefault}{\color[rgb]{0,0,0}\small{$X$}}%
}}}}
\put(1051,-1711){\makebox(0,0)[lb]{\smash{{\SetFigFont{12}{14.4}{\familydefault}{\mddefault}{\updefault}{\color[rgb]{0,0,0}\small{$Y$}}%
}}}}
\put(136,-961){\makebox(0,0)[lb]{\smash{{\SetFigFont{12}{14.4}{\familydefault}{\mddefault}{\updefault}{\color[rgb]{0,0,0}\tiny{other vertices}}%
}}}}
\put(1201,-833){\makebox(0,0)[lb]{\smash{{\SetFigFont{12}{14.4}{\familydefault}{\mddefault}{\updefault}{\color[rgb]{0,0,0}\small{$I_m$}}%
}}}}
\put(3601,-916){\makebox(0,0)[lb]{\smash{{\SetFigFont{12}{14.4}{\familydefault}{\mddefault}{\updefault}{\color[rgb]{0,0,0}\small{$I_1$}}%
}}}}
\put(2784,-909){\makebox(0,0)[lb]{\smash{{\SetFigFont{12}{14.4}{\familydefault}{\mddefault}{\updefault}{\color[rgb]{0,0,0}\small{$I_3$}}%
}}}}
\put(3084,-908){\makebox(0,0)[lb]{\smash{{\SetFigFont{12}{14.4}{\familydefault}{\mddefault}{\updefault}{\color[rgb]{0,0,0}\small{$I_2$}}%
}}}}
\put(2194,-914){\makebox(0,0)[lb]{\smash{{\SetFigFont{5}{6.0}{\familydefault}{\mddefault}{\updefault}{\color[rgb]{0,0,0}\small{$I_{m-2}$}}%
}}}}
\end{picture}%

\end{center}
\vspace{.125in}

Now we describe a sequence of Whitehead moves that can be used to connect $I_m$
to $P$ in only $P_1$ and $P_2$. This will allow us to use Move
\ref{SUBLEMMA1} to increase the length of $P$ by one.

First, using  Move \ref{SUBLEMMA3} we can eliminate
all possible connections of $I_m$ to $P$ in places other than $P_1$ and $P_3$.
Next, we do the move $Wh(I_mP_3)$ so that $I_m$ connects to $P$ only in
$P_1$.

\vspace{.125in}
\begin{center}
\begin{picture}(0,0)%
\epsfig{file=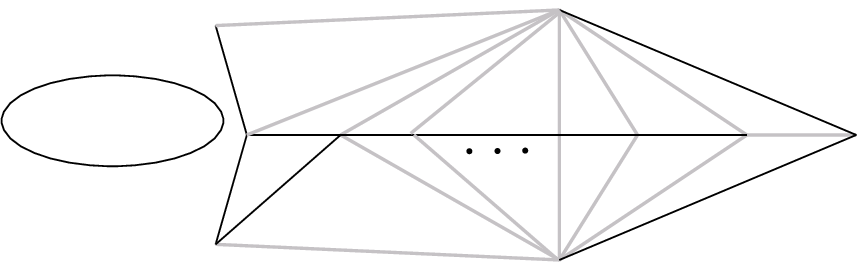}%
\end{picture}%
\setlength{\unitlength}{3947sp}%
\begingroup\makeatletter\ifx\SetFigFont\undefined%
\gdef\SetFigFont#1#2#3#4#5{%
  \reset@font\fontsize{#1}{#2pt}%
  \fontfamily{#3}\fontseries{#4}\fontshape{#5}%
  \selectfont}%
\fi\endgroup%
\begin{picture}(4185,1639)(16,-1844)
\put(1605,-1179){\makebox(0,0)[lb]{\smash{{\SetFigFont{12}{14.4}{\familydefault}{\mddefault}{\updefault}{\color[rgb]{0,0,0}\small{$I_{m-1}$}}%
}}}}
\put(2851,-361){\makebox(0,0)[lb]{\smash{{\SetFigFont{12}{14.4}{\familydefault}{\mddefault}{\updefault}{\color[rgb]{0,0,0}\small{$P_1$}}%
}}}}
\put(4201,-961){\makebox(0,0)[lb]{\smash{{\SetFigFont{12}{14.4}{\familydefault}{\mddefault}{\updefault}{\color[rgb]{0,0,0}\small{$P_2$}}%
}}}}
\put(2701,-1786){\makebox(0,0)[lb]{\smash{{\SetFigFont{12}{14.4}{\familydefault}{\mddefault}{\updefault}{\color[rgb]{0,0,0}\small{$P_3$}}%
}}}}
\put(3601,-886){\makebox(0,0)[lb]{\smash{{\SetFigFont{12}{14.4}{\familydefault}{\mddefault}{\updefault}{\color[rgb]{0,0,0}\small{$I_1$}}%
}}}}
\put(3076,-886){\makebox(0,0)[lb]{\smash{{\SetFigFont{12}{14.4}{\familydefault}{\mddefault}{\updefault}{\color[rgb]{0,0,0}\small{$I_2$}}%
}}}}
\put(2776,-886){\makebox(0,0)[lb]{\smash{{\SetFigFont{12}{14.4}{\familydefault}{\mddefault}{\updefault}{\color[rgb]{0,0,0}\small{$I_3$}}%
}}}}
\put(1201,-886){\makebox(0,0)[lb]{\smash{{\SetFigFont{12}{14.4}{\familydefault}{\mddefault}{\updefault}{\color[rgb]{0,0,0}\small{$I_m$}}%
}}}}
\put(1051,-361){\makebox(0,0)[lb]{\smash{{\SetFigFont{12}{14.4}{\familydefault}{\mddefault}{\updefault}{\color[rgb]{0,0,0}\small{$X$}}%
}}}}
\put(1051,-1711){\makebox(0,0)[lb]{\smash{{\SetFigFont{12}{14.4}{\familydefault}{\mddefault}{\updefault}{\color[rgb]{0,0,0}\small{$Y$}}%
}}}}
\put(2228,-901){\makebox(0,0)[lb]{\smash{{\SetFigFont{12}{14.4}{\familydefault}{\mddefault}{\updefault}{\color[rgb]{0,0,0}\small{$I_{m-2}$}}%
}}}}
\put(151,-961){\makebox(0,0)[lb]{\smash{{\SetFigFont{12}{14.4}{\familydefault}{\mddefault}{\updefault}{\color[rgb]{0,0,0}\tiny{other vertices}}%
}}}}
\end{picture}%

\end{center}
\vspace{.125in}

Next, we must do the moves $Wh(I_{m-1}P_1)$,...,$Wh(I_1P_1)$, in that order
(see the figure below.)

\vspace{.125in}
\begin{center}
\begin{picture}(0,0)%
\epsfig{file=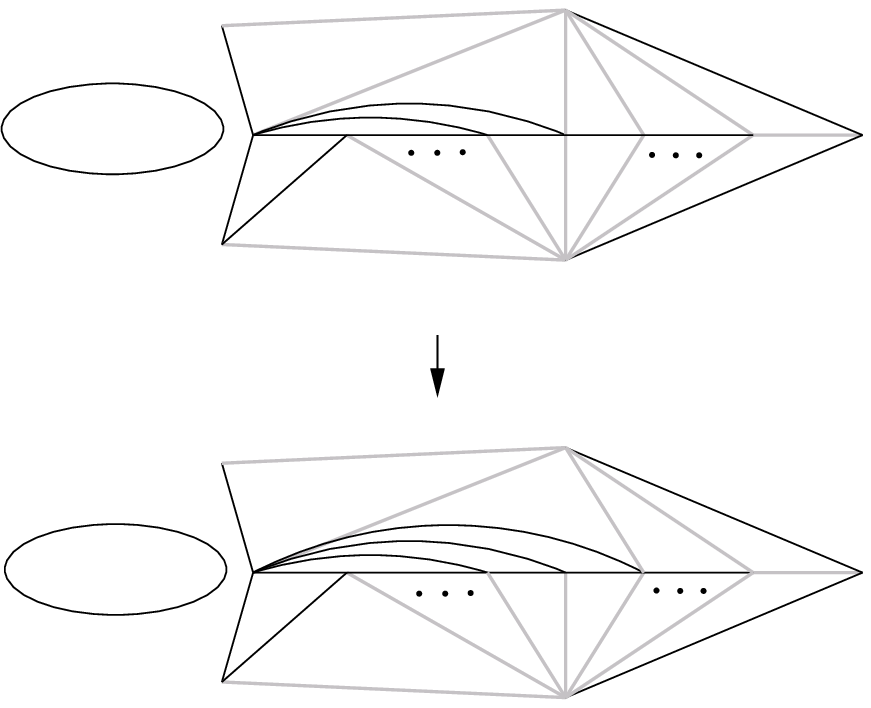}%
\end{picture}%
\setlength{\unitlength}{3947sp}%
\begingroup\makeatletter\ifx\SetFigFont\undefined%
\gdef\SetFigFont#1#2#3#4#5{%
  \reset@font\fontsize{#1}{#2pt}%
  \fontfamily{#3}\fontseries{#4}\fontshape{#5}%
  \selectfont}%
\fi\endgroup%
\begin{picture}(4215,3607)(-14,-3608)
\put(2851,-2161){\makebox(0,0)[lb]{\smash{{\SetFigFont{6}{7.2}{\familydefault}{\mddefault}{\updefault}{\color[rgb]{0,0,0}\small{$P_1$}}%
}}}}
\put(2723,-2919){\makebox(0,0)[lb]{\smash{{\SetFigFont{6}{7.2}{\familydefault}{\mddefault}{\updefault}{\color[rgb]{0,0,0}\small{$I_k$}}%
}}}}
\put(1576,-2986){\makebox(0,0)[lb]{\smash{{\SetFigFont{6}{7.2}{\familydefault}{\mddefault}{\updefault}{\color[rgb]{0,0,0}\small{$I_{m-1}$}}%
}}}}
\put(151,-2761){\makebox(0,0)[lb]{\smash{{\SetFigFont{6}{7.2}{\familydefault}{\mddefault}{\updefault}{\color[rgb]{0,0,0}\tiny{other vertices}}%
}}}}
\put(1051,-3511){\makebox(0,0)[lb]{\smash{{\SetFigFont{6}{7.2}{\familydefault}{\mddefault}{\updefault}{\color[rgb]{0,0,0}\small{$Y$}}%
}}}}
\put(1051,-2161){\makebox(0,0)[lb]{\smash{{\SetFigFont{6}{7.2}{\familydefault}{\mddefault}{\updefault}{\color[rgb]{0,0,0}\small{$X$}}%
}}}}
\put(1201,-2611){\makebox(0,0)[lb]{\smash{{\SetFigFont{6}{7.2}{\familydefault}{\mddefault}{\updefault}{\color[rgb]{0,0,0}\small{$I_m$}}%
}}}}
\put(4201,-2761){\makebox(0,0)[lb]{\smash{{\SetFigFont{6}{7.2}{\familydefault}{\mddefault}{\updefault}{\color[rgb]{0,0,0}\small{$P_2$}}%
}}}}
\put(2191,-1816){\makebox(0,0)[lb]{\smash{{\SetFigFont{6}{7.2}{\familydefault}{\mddefault}{\updefault}{\color[rgb]{0,0,0}\small{$Wh(I_kP_1)$}}%
}}}}
\put(2851,-61){\makebox(0,0)[lb]{\smash{{\SetFigFont{6}{7.2}{\familydefault}{\mddefault}{\updefault}{\color[rgb]{0,0,0}\small{$P_1$}}%
}}}}
\put(2723,-819){\makebox(0,0)[lb]{\smash{{\SetFigFont{6}{7.2}{\familydefault}{\mddefault}{\updefault}{\color[rgb]{0,0,0}\small{$I_k$}}%
}}}}
\put(1576,-886){\makebox(0,0)[lb]{\smash{{\SetFigFont{6}{7.2}{\familydefault}{\mddefault}{\updefault}{\color[rgb]{0,0,0}\small{$I_{m-1}$}}%
}}}}
\put(151,-661){\makebox(0,0)[lb]{\smash{{\SetFigFont{6}{7.2}{\familydefault}{\mddefault}{\updefault}{\color[rgb]{0,0,0}\tiny{other vertices}}%
}}}}
\put(1051,-1411){\makebox(0,0)[lb]{\smash{{\SetFigFont{6}{7.2}{\familydefault}{\mddefault}{\updefault}{\color[rgb]{0,0,0}\small{$Y$}}%
}}}}
\put(1051,-61){\makebox(0,0)[lb]{\smash{{\SetFigFont{6}{7.2}{\familydefault}{\mddefault}{\updefault}{\color[rgb]{0,0,0}\small{$X$}}%
}}}}
\put(1201,-511){\makebox(0,0)[lb]{\smash{{\SetFigFont{6}{7.2}{\familydefault}{\mddefault}{\updefault}{\color[rgb]{0,0,0}\small{$I_m$}}%
}}}}
\put(4201,-661){\makebox(0,0)[lb]{\smash{{\SetFigFont{6}{7.2}{\familydefault}{\mddefault}{\updefault}{\color[rgb]{0,0,0}\small{$P_2$}}%
}}}}
\put(2701,-1486){\makebox(0,0)[lb]{\smash{{\SetFigFont{6}{7.2}{\familydefault}{\mddefault}{\updefault}{\color[rgb]{0,0,0}\small{$P_3$}}%
}}}}
\put(3076,-586){\makebox(0,0)[lb]{\smash{{\SetFigFont{6}{7.2}{\familydefault}{\mddefault}{\updefault}{\color[rgb]{0,0,0}\small{$I_{k-1}$}}%
}}}}
\put(3076,-2686){\makebox(0,0)[lb]{\smash{{\SetFigFont{6}{7.2}{\familydefault}{\mddefault}{\updefault}{\color[rgb]{0,0,0}\small{$I_{k-1}$}}%
}}}}
\put(2414,-819){\makebox(0,0)[lb]{\smash{{\SetFigFont{6}{7.2}{\familydefault}{\mddefault}{\updefault}{\color[rgb]{0,0,0}\small{$I_{k+1}$}}%
}}}}
\put(2701,-3586){\makebox(0,0)[lb]{\smash{{\SetFigFont{6}{7.2}{\familydefault}{\mddefault}{\updefault}{\color[rgb]{0,0,0}\small{$P_3$}}%
}}}}
\put(2414,-2919){\makebox(0,0)[lb]{\smash{{\SetFigFont{6}{7.2}{\familydefault}{\mddefault}{\updefault}{\color[rgb]{0,0,0}\small{$I_{k+1}$}}%
}}}}
\put(3587,-798){\makebox(0,0)[lb]{\smash{{\SetFigFont{6}{7.2}{\familydefault}{\mddefault}{\updefault}{\color[rgb]{0,0,0}\small{$I_1$}}%
}}}}
\put(3592,-2893){\makebox(0,0)[lb]{\smash{{\SetFigFont{6}{7.2}{\familydefault}{\mddefault}{\updefault}{\color[rgb]{0,0,0}\small{$I_1$}}%
}}}}
\end{picture}%

\end{center}
\vspace{.125in}

After this sequence of Whitehead moves we obtain the first diagram below,
with $I_m$ connected to $P$ exactly at $P_1$ and $P_2$, so that we can apply
Move \ref{SUBLEMMA1} to increase the length of $P$ by the move $Wh(P_1P_2)$,
below.

\vspace{.125in}
\begin{center}
\begin{picture}(0,0)%
\epsfig{file=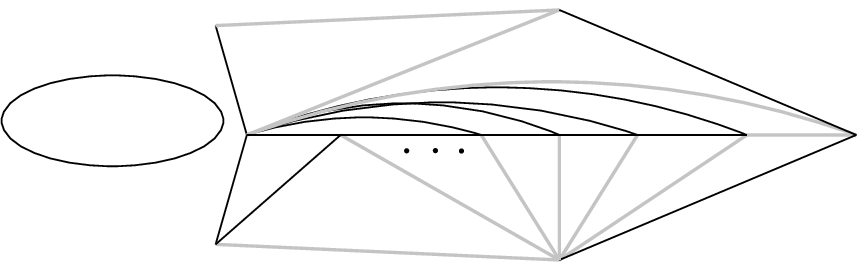}%
\end{picture}%
\setlength{\unitlength}{3947sp}%
\begingroup\makeatletter\ifx\SetFigFont\undefined%
\gdef\SetFigFont#1#2#3#4#5{%
  \reset@font\fontsize{#1}{#2pt}%
  \fontfamily{#3}\fontseries{#4}\fontshape{#5}%
  \selectfont}%
\fi\endgroup%
\begin{picture}(4185,1639)(16,-1844)
\put(3603,-1092){\makebox(0,0)[lb]{\smash{{\SetFigFont{12}{14.4}{\familydefault}{\mddefault}{\updefault}{\color[rgb]{0,0,0}\small{$I_1$}}%
}}}}
\put(2851,-361){\makebox(0,0)[lb]{\smash{{\SetFigFont{12}{14.4}{\familydefault}{\mddefault}{\updefault}{\color[rgb]{0,0,0}\small{$P_1$}}%
}}}}
\put(4201,-961){\makebox(0,0)[lb]{\smash{{\SetFigFont{12}{14.4}{\familydefault}{\mddefault}{\updefault}{\color[rgb]{0,0,0}\small{$P_2$}}%
}}}}
\put(2701,-1786){\makebox(0,0)[lb]{\smash{{\SetFigFont{12}{14.4}{\familydefault}{\mddefault}{\updefault}{\color[rgb]{0,0,0}\small{$P_3$}}%
}}}}
\put(1051,-361){\makebox(0,0)[lb]{\smash{{\SetFigFont{12}{14.4}{\familydefault}{\mddefault}{\updefault}{\color[rgb]{0,0,0}\small{$X$}}%
}}}}
\put(1051,-1711){\makebox(0,0)[lb]{\smash{{\SetFigFont{12}{14.4}{\familydefault}{\mddefault}{\updefault}{\color[rgb]{0,0,0}\small{$Y$}}%
}}}}
\put(151,-961){\makebox(0,0)[lb]{\smash{{\SetFigFont{12}{14.4}{\familydefault}{\mddefault}{\updefault}{\color[rgb]{0,0,0}\tiny{other vertices}}%
}}}}
\put(3076,-1118){\makebox(0,0)[lb]{\smash{{\SetFigFont{12}{14.4}{\familydefault}{\mddefault}{\updefault}{\color[rgb]{0,0,0}\small{$I_2$}}%
}}}}
\put(2723,-1119){\makebox(0,0)[lb]{\smash{{\SetFigFont{12}{14.4}{\familydefault}{\mddefault}{\updefault}{\color[rgb]{0,0,0}\small{$I_3$}}%
}}}}
\put(2415,-1126){\makebox(0,0)[lb]{\smash{{\SetFigFont{12}{14.4}{\familydefault}{\mddefault}{\updefault}{\color[rgb]{0,0,0}\small{$I_4$}}%
}}}}
\put(1201,-811){\makebox(0,0)[lb]{\smash{{\SetFigFont{12}{14.4}{\familydefault}{\mddefault}{\updefault}{\color[rgb]{0,0,0}\small{$I_m$}}%
}}}}
\put(1576,-1186){\makebox(0,0)[lb]{\smash{{\SetFigFont{12}{14.4}{\familydefault}{\mddefault}{\updefault}{\color[rgb]{0,0,0}\small{$I_{m-1}$}}%
}}}}
\end{picture}%

\end{center}
\vspace{.125in}

\vspace{.125in}
\begin{center}
\begin{picture}(0,0)%
\epsfig{file=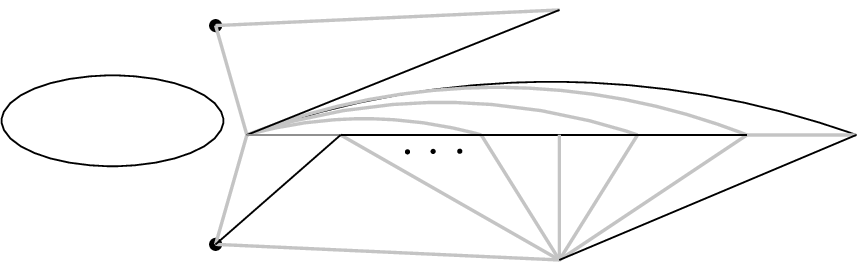}%
\end{picture}%
\setlength{\unitlength}{3947sp}%
\begingroup\makeatletter\ifx\SetFigFont\undefined%
\gdef\SetFigFont#1#2#3#4#5{%
  \reset@font\fontsize{#1}{#2pt}%
  \fontfamily{#3}\fontseries{#4}\fontshape{#5}%
  \selectfont}%
\fi\endgroup%
\begin{picture}(4185,1639)(16,-1844)
\put(3596,-1094){\makebox(0,0)[lb]{\smash{{\SetFigFont{12}{14.4}{\familydefault}{\mddefault}{\updefault}{\color[rgb]{0,0,0}\small{$I_1$}}%
}}}}
\put(2851,-361){\makebox(0,0)[lb]{\smash{{\SetFigFont{12}{14.4}{\familydefault}{\mddefault}{\updefault}{\color[rgb]{0,0,0}\small{$P_1$}}%
}}}}
\put(2701,-1786){\makebox(0,0)[lb]{\smash{{\SetFigFont{12}{14.4}{\familydefault}{\mddefault}{\updefault}{\color[rgb]{0,0,0}\small{$P_3$}}%
}}}}
\put(1051,-361){\makebox(0,0)[lb]{\smash{{\SetFigFont{12}{14.4}{\familydefault}{\mddefault}{\updefault}{\color[rgb]{0,0,0}\small{$X$}}%
}}}}
\put(1051,-1711){\makebox(0,0)[lb]{\smash{{\SetFigFont{12}{14.4}{\familydefault}{\mddefault}{\updefault}{\color[rgb]{0,0,0}\small{$Y$}}%
}}}}
\put(151,-961){\makebox(0,0)[lb]{\smash{{\SetFigFont{12}{14.4}{\familydefault}{\mddefault}{\updefault}{\color[rgb]{0,0,0}\tiny{other vertices}}%
}}}}
\put(1576,-1186){\makebox(0,0)[lb]{\smash{{\SetFigFont{12}{14.4}{\familydefault}{\mddefault}{\updefault}{\color[rgb]{0,0,0}\small{$I_{m-1}$}}%
}}}}
\put(3076,-1118){\makebox(0,0)[lb]{\smash{{\SetFigFont{12}{14.4}{\familydefault}{\mddefault}{\updefault}{\color[rgb]{0,0,0}\small{$I_2$}}%
}}}}
\put(2723,-1119){\makebox(0,0)[lb]{\smash{{\SetFigFont{12}{14.4}{\familydefault}{\mddefault}{\updefault}{\color[rgb]{0,0,0}\small{$I_3$}}%
}}}}
\put(2415,-1126){\makebox(0,0)[lb]{\smash{{\SetFigFont{12}{14.4}{\familydefault}{\mddefault}{\updefault}{\color[rgb]{0,0,0}\small{$I_4$}}%
}}}}
\put(1201,-811){\makebox(0,0)[lb]{\smash{{\SetFigFont{12}{14.4}{\familydefault}{\mddefault}{\updefault}{\color[rgb]{0,0,0}\small{$I_m$}}%
}}}}
\put(4201,-961){\makebox(0,0)[lb]{\smash{{\SetFigFont{12}{14.4}{\familydefault}{\mddefault}{\updefault}{\color[rgb]{0,0,0}\small{$P_2$}}%
}}}}
\end{picture}%

\end{center}
\vspace{.125in}

This concludes Case 3.

\vspace{.1in}

Since $C^*$ must belong to one of these cases, we have seen that if the
length of $P$ is less than $N-3$, we can do Whitehead moves to increase it
to $N-3$ without creating prismatic $3$-circuits.  Hence we can
reduce to the case of two interior vertices, both of which must be
endpoints.  Then we can apply Move \ref{SUBLEMMA2} part (b) to decrease the
number of connections between one of these two interior vertices and $P$
to exactly $3$.  The result is the complex $D_N$, as shown to the right below.

\vspace{.125in}
\begin{center}
%

\end{center}
\vspace{.125in}

\Endproof
\vspace{.2in}

\subsection{Construction of a ``difficult'' simple polyhedron}\label{SEC_DIFFICULT}
We illustrate the algorithm described in the previous section by constructing a
hyperbolic polyhedron realizing an abstract polyhedron for which Case 3 from the proof of Lemma \ref{ROEDER} is
necessary, and hence is particularly difficult.

\begin{center}
\scalebox{1.1}{

%
%
}
\end{center}

\vspace{.1in}
Following the Whitehead moves backwards from $D_{18}^*$ back to $R_{18}^*$
we find the following sequence of Whitehead moves. 

\vspace{.1in}
\noindent
{\footnotesize
$Wh(6,13)$, $Wh(6,9)$, $Wh(3,6)$, $Wh(9,12)$, $Wh(6,15)$, $Wh(6,16)$, $Wh(6,17)$,
$Wh(6,7)$, $Wh(6,8)$, $Wh(6,4)$, $Wh(8,18)$, $Wh(7,9)$, $Wh(9,14)$, $Wh(9,15)$, $Wh(3,18)$,
$Wh(7,8)$, $Wh(4,18)$, $Wh(3,8)$, $Wh(8,9)$, $Wh(5,18)$, $Wh(4,9)$, $Wh(6,9)$, $Wh(2,18)$,
$Wh(5,6)$, $Wh(1,18)$, $Wh(1,13)$, $Wh(1,3)$, $Wh(3,4)$, $Wh(4,5)$, $Wh(2,5)$.
}

\begin{figure}[p]\label{FIG_COMB2}
\begin{center}
\includegraphics[scale=1.2]{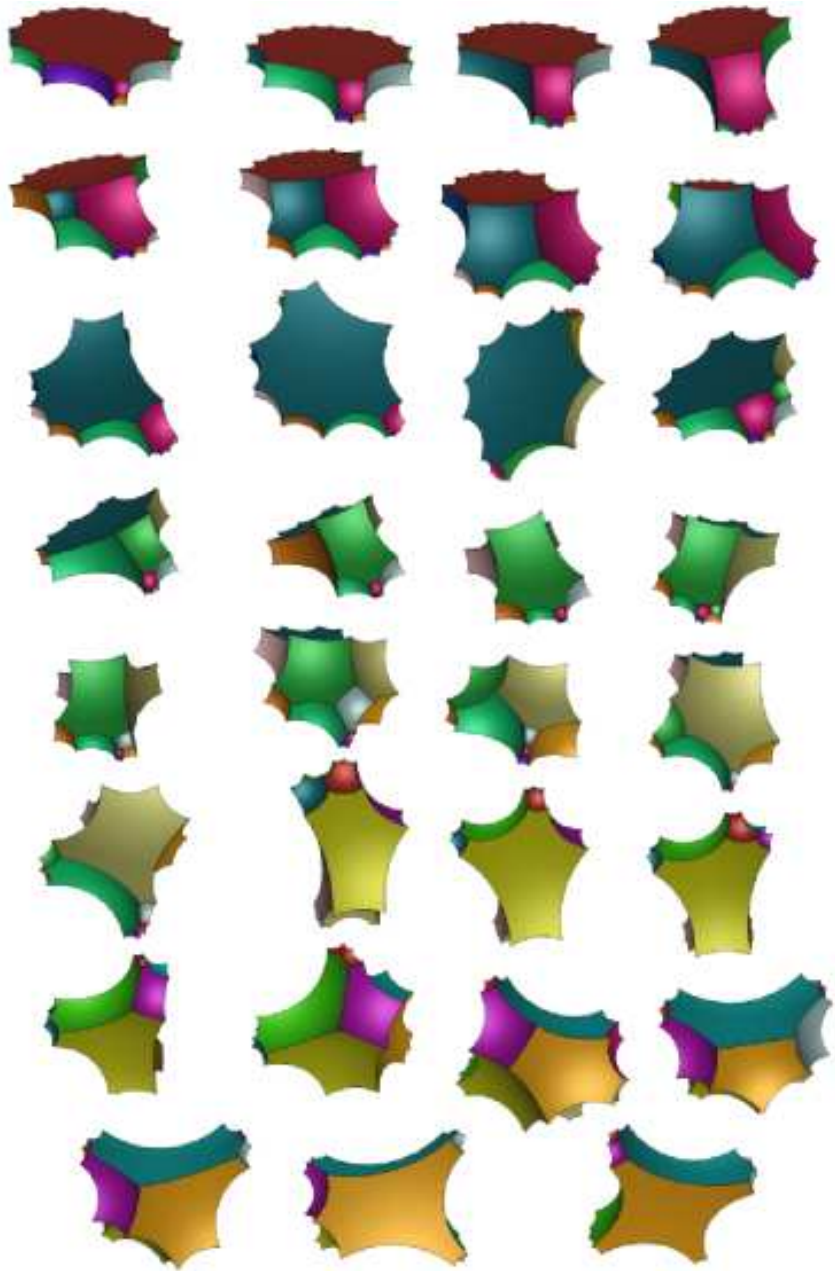}
\end{center}
\caption{Construction of $R_{18}$ from $D_{18}$ using $30$ Whitehead moves.}
\end{figure}

We did this sequence of Whitehead moves geometrically, using Newton's Method.
The result, starting with $D_{18}$, and realizing $R_{18}$ is shown in Figure
\ref{FIG_COMB2}.  Each polyhedron is displayed in the conformal ball model. 

\subsection{Truncation of vertices}\label{SEC_TRUNCATION}

We have seen an outline of how to construct simple polyhedra. 
We now show how to construct all truncated polyhedra, except for the
triangular prism, which we have already constructed in Section \ref{SEC_WH}.

\begin{lem} \label{PI3} If $A_C \neq \emptyset$, then there are points in
$A_C$ arbitrarily close to \\ $(\pi/3,\pi/3,\ldots,\pi/3)$.  \end{lem}

\noindent {\bf Proof:} 
Simply check that if ${\bf a} \in A_C$, then the entire straight line path
to $(\pi/3,\pi/3,\ldots,\pi/3)$, excluding the final point is in $A_C$.
\Endproof

Thus we can assume that ${\bf a}$ is arbitrarily close to
$(\pi/3,\pi/3,\cdots,\pi/3)$, because once we have a polyhedron realizing $C$
with non-obtuse dihedral angles, we can deform it to have any dihedral angles
in $A_C$, as described in Section \ref{SEC_DEFORM}.
Specifically, choose some $0 < \delta < \frac{\pi}{18}$ and assume that each
component of ${\bf a}$ is within $\delta$ of $\frac{\pi}{3}$.

Let $\widetilde{C}$ be the $C$ with each of the triangular faces $f_i^{\rm tr}$
replaced by a single vertex $v_i^{\rm tr}$.  (Or, if $C$ is the truncated triangular prism, let $\widetilde{C}$
be the prism.)
Let $\hat{{\bf a}}$ be the angles
from ${\bf a}$ corresponding to the edges from $C$ that are in $\widetilde{C}$
and let $\beta = (\hat{{\bf a}}_1+2\delta,\hat{{\bf a}}_2+2\delta,\cdots)$.  (If $\widetilde{C}$
is the prism, re-number the edges so that the three edges forming the prismatic cycle
are the first three, and choose $\beta = (\hat{{\bf a}}_1,\hat{{\bf a}}_2,\hat{{\bf a}}_1,\hat{{\bf a}}_5+2\delta,\hat{{\bf a}}_5+2\delta,\cdots)$.)

Note that $\delta$ was chosen so that $\beta \in A_{\widetilde{C}}$.  Then, the straight
line path ${\bf a}(t)$ joining $\beta$ to $\hat{\bf a}$ (parameterized by $t
\in (0,1)$) will remain in $A_{\widetilde{C}}$ except that the sum of the dihedral angles
of edges meeting at each of the vertices $v_i^{\rm tr}$ will decrease past $\pi$ at some time
$t_i \in (0,1)$.

In \cite{ROE2,ROE} the authors use the path ${\bf a}(t)$ to construct a sequence of
polyhedra $\widetilde{P} = P_0,P_1,\cdots,P_{N-1} = P$ where $\widetilde{P}$ realizes
$C$, and $P_i$ is obtained from $P_{i-1}$ truncating the vertices that become ideal
when $t=t_i$.  (The truncation is done using a small perturbation.)  Realizing
$P$ proves that ${\cal P}_C^0 \neq \emptyset$, as needed for the proof of
Andreev's Theorem.

Because the proof from \cite{ROE2,ROE} gives us a priori knowledge that compact
polyhedra exist realizing each of the intermediate combinatorial structures exist, we
can use Newton's method to deform the planes forming $\widetilde{P}$ to realize the
angles in the entire path ${\bf a}(t)$ without truncating each vertex once it
meets $\partial_\infty \mathbb{H}^3$.  We can then solve independently for the
planes corresponding to the triangles in $C$ so that each intersects the three
appropriate planes at the three appropriate angles.  

We illustrate this construction below for the following pair $(C,{\bf a})$, which has
three truncations labeled $f_1^{\rm tr}$, $f_2^{\rm tr}$, and $f_3^{\rm tr}$.  

\begin{center}

\begin{picture}(0,0)%
\epsfig{file=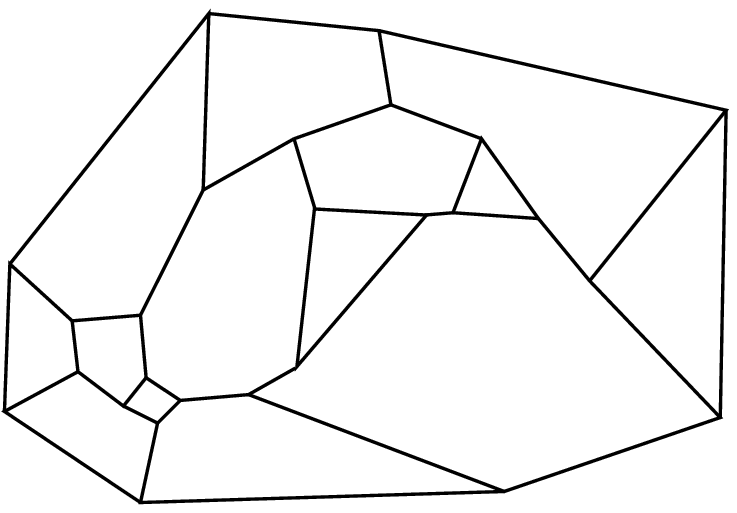}%
\end{picture}%
\setlength{\unitlength}{3947sp}%
\begingroup\makeatletter\ifx\SetFigFont\undefined%
\gdef\SetFigFont#1#2#3#4#5{%
  \reset@font\fontsize{#1}{#2pt}%
  \fontfamily{#3}\fontseries{#4}\fontshape{#5}%
  \selectfont}%
\fi\endgroup%
\begin{picture}(3532,2581)(137,-1763)
\put(3258,-581){\makebox(0,0)[lb]{\smash{{\SetFigFont{8}{9.6}{\familydefault}{\mddefault}{\updefault}{\color[rgb]{0,0,0}$f_3^{\rm tr}$}%
}}}}
\put(356,-617){\makebox(0,0)[lb]{\smash{{\SetFigFont{5}{6.0}{\familydefault}{\mddefault}{\updefault}{\color[rgb]{0,0,0}$\frac{2\pi}{5}$}%
}}}}
\put(609,-725){\makebox(0,0)[lb]{\smash{{\SetFigFont{5}{6.0}{\familydefault}{\mddefault}{\updefault}{\color[rgb]{0,0,0}$\frac{2\pi}{5}$}%
}}}}
\put(1003,-460){\makebox(0,0)[lb]{\smash{{\SetFigFont{5}{6.0}{\familydefault}{\mddefault}{\updefault}{\color[rgb]{0,0,0}$\frac{2\pi}{5}$}%
}}}}
\put(1163,319){\makebox(0,0)[lb]{\smash{{\SetFigFont{5}{6.0}{\familydefault}{\mddefault}{\updefault}{\color[rgb]{0,0,0}$\frac{\pi}{2}$}%
}}}}
\put(2615,-111){\makebox(0,0)[lb]{\smash{{\SetFigFont{5}{6.0}{\familydefault}{\mddefault}{\updefault}{\color[rgb]{0,0,0}$\frac{\pi}{2}$}%
}}}}
\put(2262,-89){\makebox(0,0)[lb]{\smash{{\SetFigFont{5}{6.0}{\familydefault}{\mddefault}{\updefault}{\color[rgb]{0,0,0}$\frac{\pi}{2}$}%
}}}}
\put(1745, 77){\makebox(0,0)[lb]{\smash{{\SetFigFont{5}{6.0}{\familydefault}{\mddefault}{\updefault}{\color[rgb]{0,0,0}$\frac{\pi}{2}$}%
}}}}
\put(3045,-1500){\makebox(0,0)[lb]{\smash{{\SetFigFont{5}{6.0}{\familydefault}{\mddefault}{\updefault}{\color[rgb]{0,0,0}$\frac{\pi}{4}$}%
}}}}
\put(2010,-1360){\makebox(0,0)[lb]{\smash{{\SetFigFont{5}{6.0}{\familydefault}{\mddefault}{\updefault}{\color[rgb]{0,0,0}$\frac{2\pi}{5}$}%
}}}}
\put(958,-1261){\makebox(0,0)[lb]{\smash{{\SetFigFont{5}{6.0}{\familydefault}{\mddefault}{\updefault}{\color[rgb]{0,0,0}$\frac{2\pi}{5}$}%
}}}}
\put(796,-1163){\makebox(0,0)[lb]{\smash{{\SetFigFont{5}{6.0}{\familydefault}{\mddefault}{\updefault}{\color[rgb]{0,0,0}$\frac{2\pi}{5}$}%
}}}}
\put(856,-902){\makebox(0,0)[lb]{\smash{{\SetFigFont{5}{6.0}{\familydefault}{\mddefault}{\updefault}{\color[rgb]{0,0,0}$\frac{2\pi}{5}$}%
}}}}
\put(522,-909){\makebox(0,0)[lb]{\smash{{\SetFigFont{5}{6.0}{\familydefault}{\mddefault}{\updefault}{\color[rgb]{0,0,0}$\frac{2\pi}{5}$}%
}}}}
\put(2818,-390){\makebox(0,0)[lb]{\smash{{\SetFigFont{5}{6.0}{\familydefault}{\mddefault}{\updefault}{\color[rgb]{0,0,0}$\frac{\pi}{4}$}%
}}}}
\put(1319,-112){\makebox(0,0)[lb]{\smash{{\SetFigFont{5}{6.0}{\familydefault}{\mddefault}{\updefault}{\color[rgb]{0,0,0}$\frac{\pi}{2}$}%
}}}}
\put(1630,-98){\makebox(0,0)[lb]{\smash{{\SetFigFont{5}{6.0}{\familydefault}{\mddefault}{\updefault}{\color[rgb]{0,0,0}$\frac{\pi}{4}$}%
}}}}
\put(1513,-1120){\makebox(0,0)[lb]{\smash{{\SetFigFont{5}{6.0}{\familydefault}{\mddefault}{\updefault}{\color[rgb]{0,0,0}$\frac{\pi}{3}$}%
}}}}
\put(1881,-729){\makebox(0,0)[lb]{\smash{{\SetFigFont{5}{6.0}{\familydefault}{\mddefault}{\updefault}{\color[rgb]{0,0,0}$\frac{\pi}{2}$}%
}}}}
\put(1139,-1242){\makebox(0,0)[lb]{\smash{{\SetFigFont{5}{6.0}{\familydefault}{\mddefault}{\updefault}{\color[rgb]{0,0,0}$\frac{2\pi}{5}$}%
}}}}
\put(2471,-361){\makebox(0,0)[lb]{\smash{{\SetFigFont{5}{6.0}{\familydefault}{\mddefault}{\updefault}{\color[rgb]{0,0,0}$\frac{\pi}{2}$}%
}}}}
\put(1501,-1736){\makebox(0,0)[lb]{\smash{{\SetFigFont{5}{6.0}{\familydefault}{\mddefault}{\updefault}{\color[rgb]{0,0,0}$\frac{2\pi}{5}$}%
}}}}
\put(326,-1511){\makebox(0,0)[lb]{\smash{{\SetFigFont{5}{6.0}{\familydefault}{\mddefault}{\updefault}{\color[rgb]{0,0,0}$\frac{2\pi}{5}$}%
}}}}
\put(882,-1482){\makebox(0,0)[lb]{\smash{{\SetFigFont{5}{6.0}{\familydefault}{\mddefault}{\updefault}{\color[rgb]{0,0,0}$\frac{2\pi}{5}$}%
}}}}
\put(695,-1321){\makebox(0,0)[lb]{\smash{{\SetFigFont{5}{6.0}{\familydefault}{\mddefault}{\updefault}{\color[rgb]{0,0,0}$\frac{2\pi}{5}$}%
}}}}
\put(336,-1217){\makebox(0,0)[lb]{\smash{{\SetFigFont{5}{6.0}{\familydefault}{\mddefault}{\updefault}{\color[rgb]{0,0,0}$\frac{2\pi}{5}$}%
}}}}
\put(476, 68){\makebox(0,0)[lb]{\smash{{\SetFigFont{5}{6.0}{\familydefault}{\mddefault}{\updefault}{\color[rgb]{0,0,0}$\frac{2\pi}{5}$}%
}}}}
\put(2237,-352){\makebox(0,0)[lb]{\smash{{\SetFigFont{5}{6.0}{\familydefault}{\mddefault}{\updefault}{\color[rgb]{0,0,0}$\frac{\pi}{6}$}%
}}}}
\put(1490,-559){\makebox(0,0)[lb]{\smash{{\SetFigFont{5}{6.0}{\familydefault}{\mddefault}{\updefault}{\color[rgb]{0,0,0}$\frac{\pi}{2}$}%
}}}}
\put(190,-912){\makebox(0,0)[lb]{\smash{{\SetFigFont{5}{6.0}{\familydefault}{\mddefault}{\updefault}{\color[rgb]{0,0,0}$\frac{2\pi}{5}$}%
}}}}
\put(1458,746){\makebox(0,0)[lb]{\smash{{\SetFigFont{5}{6.0}{\familydefault}{\mddefault}{\updefault}{\color[rgb]{0,0,0}$\frac{2\pi}{5}$}%
}}}}
\put(2711,525){\makebox(0,0)[lb]{\smash{{\SetFigFont{5}{6.0}{\familydefault}{\mddefault}{\updefault}{\color[rgb]{0,0,0}$\frac{\pi}{4}$}%
}}}}
\put(1856,-177){\makebox(0,0)[lb]{\smash{{\SetFigFont{5}{6.0}{\familydefault}{\mddefault}{\updefault}{\color[rgb]{0,0,0}$\frac{\pi}{2}$}%
}}}}
\put(598,-1060){\makebox(0,0)[lb]{\smash{{\SetFigFont{5}{6.0}{\familydefault}{\mddefault}{\updefault}{\color[rgb]{0,0,0}$\frac{2\pi}{5}$}%
}}}}
\put(921,-1070){\makebox(0,0)[lb]{\smash{{\SetFigFont{5}{6.0}{\familydefault}{\mddefault}{\updefault}{\color[rgb]{0,0,0}$\frac{2\pi}{5}$}%
}}}}
\put(2234,251){\makebox(0,0)[lb]{\smash{{\SetFigFont{5}{6.0}{\familydefault}{\mddefault}{\updefault}{\color[rgb]{0,0,0}$\frac{\pi}{4}$}%
}}}}
\put(2401,-211){\makebox(0,0)[lb]{\smash{{\SetFigFont{8}{9.6}{\familydefault}{\mddefault}{\updefault}{\color[rgb]{0,0,0}$f_2^{\rm tr}$}%
}}}}
\put(1726,-511){\makebox(0,0)[lb]{\smash{{\SetFigFont{8}{9.6}{\familydefault}{\mddefault}{\updefault}{\color[rgb]{0,0,0}$f_1^{\rm tr}$}%
}}}}
\put(3119,-951){\makebox(0,0)[lb]{\smash{{\SetFigFont{5}{6.0}{\familydefault}{\mddefault}{\updefault}{\color[rgb]{0,0,0}$\frac{\pi}{2}$}%
}}}}
\put(3128,-165){\makebox(0,0)[lb]{\smash{{\SetFigFont{5}{6.0}{\familydefault}{\mddefault}{\updefault}{\color[rgb]{0,0,0}$\frac{\pi}{2}$}%
}}}}
\put(3669,-496){\makebox(0,0)[lb]{\smash{{\SetFigFont{5}{6.0}{\familydefault}{\mddefault}{\updefault}{\color[rgb]{0,0,0}$\frac{\pi}{2}$}%
}}}}
\put(2036,420){\makebox(0,0)[lb]{\smash{{\SetFigFont{5}{6.0}{\familydefault}{\mddefault}{\updefault}{\color[rgb]{0,0,0}$\frac{2\pi}{5}$}%
}}}}
\end{picture}%
\end{center}

The path ${\bf a}(t)$ described above works with any truncated $C$.  For many
$C$, such as the current one, a much easier path ${\bf a}(t)$ can be found
satisfying conditions (1,3,4) and (5) from Andreev's Theorem for $\widetilde{C}$,
but for which the sum of the dihedral angles of the edges meeting at each
$v_i^{\rm tr}$ decreases past $\pi$ at some $t_i \in (0,1)$.  Such a path is sufficient for
our construction.  For the current construction, $\widetilde{C}$ is shown below with edges labeled
according to an appropriate path ${\bf a}(t)$.  (Notation: ${\bf a}(t) = \beta$ for the edges labeled with a single angle
$\beta$, whereas ${\bf a}(t) = \eta (1-t) + \gamma t$ for the edges  $[\eta, \gamma]$.)

\begin{center}

\begin{picture}(0,0)%
\epsfig{file=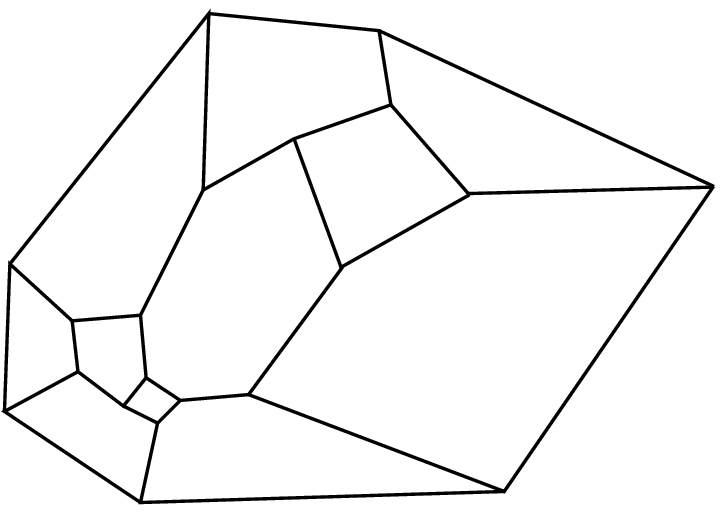}%
\end{picture}%
\setlength{\unitlength}{3947sp}%
\begingroup\makeatletter\ifx\SetFigFont\undefined%
\gdef\SetFigFont#1#2#3#4#5{%
  \reset@font\fontsize{#1}{#2pt}%
  \fontfamily{#3}\fontseries{#4}\fontshape{#5}%
  \selectfont}%
\fi\endgroup%
\begin{picture}(3452,2581)(137,-1763)
\put(2769,-296){\makebox(0,0)[lb]{\smash{{\SetFigFont{5}{6.0}{\familydefault}{\mddefault}{\updefault}{\color[rgb]{0,0,0}$\left[\frac{2\pi}{5},\frac{\pi}{4}\right]$}%
}}}}
\put(356,-617){\makebox(0,0)[lb]{\smash{{\SetFigFont{5}{6.0}{\familydefault}{\mddefault}{\updefault}{\color[rgb]{0,0,0}$\frac{2\pi}{5}$}%
}}}}
\put(609,-725){\makebox(0,0)[lb]{\smash{{\SetFigFont{5}{6.0}{\familydefault}{\mddefault}{\updefault}{\color[rgb]{0,0,0}$\frac{2\pi}{5}$}%
}}}}
\put(1003,-460){\makebox(0,0)[lb]{\smash{{\SetFigFont{5}{6.0}{\familydefault}{\mddefault}{\updefault}{\color[rgb]{0,0,0}$\frac{2\pi}{5}$}%
}}}}
\put(1163,319){\makebox(0,0)[lb]{\smash{{\SetFigFont{5}{6.0}{\familydefault}{\mddefault}{\updefault}{\color[rgb]{0,0,0}$\frac{\pi}{2}$}%
}}}}
\put(1745, 77){\makebox(0,0)[lb]{\smash{{\SetFigFont{5}{6.0}{\familydefault}{\mddefault}{\updefault}{\color[rgb]{0,0,0}$\frac{\pi}{2}$}%
}}}}
\put(958,-1261){\makebox(0,0)[lb]{\smash{{\SetFigFont{5}{6.0}{\familydefault}{\mddefault}{\updefault}{\color[rgb]{0,0,0}$\frac{2\pi}{5}$}%
}}}}
\put(796,-1163){\makebox(0,0)[lb]{\smash{{\SetFigFont{5}{6.0}{\familydefault}{\mddefault}{\updefault}{\color[rgb]{0,0,0}$\frac{2\pi}{5}$}%
}}}}
\put(856,-902){\makebox(0,0)[lb]{\smash{{\SetFigFont{5}{6.0}{\familydefault}{\mddefault}{\updefault}{\color[rgb]{0,0,0}$\frac{2\pi}{5}$}%
}}}}
\put(522,-909){\makebox(0,0)[lb]{\smash{{\SetFigFont{5}{6.0}{\familydefault}{\mddefault}{\updefault}{\color[rgb]{0,0,0}$\frac{2\pi}{5}$}%
}}}}
\put(1319,-112){\makebox(0,0)[lb]{\smash{{\SetFigFont{5}{6.0}{\familydefault}{\mddefault}{\updefault}{\color[rgb]{0,0,0}$\frac{\pi}{2}$}%
}}}}
\put(1139,-1242){\makebox(0,0)[lb]{\smash{{\SetFigFont{5}{6.0}{\familydefault}{\mddefault}{\updefault}{\color[rgb]{0,0,0}$\frac{2\pi}{5}$}%
}}}}
\put(326,-1511){\makebox(0,0)[lb]{\smash{{\SetFigFont{5}{6.0}{\familydefault}{\mddefault}{\updefault}{\color[rgb]{0,0,0}$\frac{2\pi}{5}$}%
}}}}
\put(882,-1482){\makebox(0,0)[lb]{\smash{{\SetFigFont{5}{6.0}{\familydefault}{\mddefault}{\updefault}{\color[rgb]{0,0,0}$\frac{2\pi}{5}$}%
}}}}
\put(695,-1321){\makebox(0,0)[lb]{\smash{{\SetFigFont{5}{6.0}{\familydefault}{\mddefault}{\updefault}{\color[rgb]{0,0,0}$\frac{2\pi}{5}$}%
}}}}
\put(336,-1217){\makebox(0,0)[lb]{\smash{{\SetFigFont{5}{6.0}{\familydefault}{\mddefault}{\updefault}{\color[rgb]{0,0,0}$\frac{2\pi}{5}$}%
}}}}
\put(476, 68){\makebox(0,0)[lb]{\smash{{\SetFigFont{5}{6.0}{\familydefault}{\mddefault}{\updefault}{\color[rgb]{0,0,0}$\frac{2\pi}{5}$}%
}}}}
\put(190,-912){\makebox(0,0)[lb]{\smash{{\SetFigFont{5}{6.0}{\familydefault}{\mddefault}{\updefault}{\color[rgb]{0,0,0}$\frac{2\pi}{5}$}%
}}}}
\put(1458,746){\makebox(0,0)[lb]{\smash{{\SetFigFont{5}{6.0}{\familydefault}{\mddefault}{\updefault}{\color[rgb]{0,0,0}$\frac{2\pi}{5}$}%
}}}}
\put(598,-1060){\makebox(0,0)[lb]{\smash{{\SetFigFont{5}{6.0}{\familydefault}{\mddefault}{\updefault}{\color[rgb]{0,0,0}$\frac{2\pi}{5}$}%
}}}}
\put(921,-1070){\makebox(0,0)[lb]{\smash{{\SetFigFont{5}{6.0}{\familydefault}{\mddefault}{\updefault}{\color[rgb]{0,0,0}$\frac{2\pi}{5}$}%
}}}}
\put(1764,398){\makebox(0,0)[lb]{\smash{{\SetFigFont{5}{6.0}{\familydefault}{\mddefault}{\updefault}{\color[rgb]{0,0,0}$\frac{2\pi}{5}$}%
}}}}
\put(1922,-1317){\makebox(0,0)[lb]{\smash{{\SetFigFont{5}{6.0}{\familydefault}{\mddefault}{\updefault}{\color[rgb]{0,0,0}$\frac{2\pi}{5}$}%
}}}}
\put(1501,-1736){\makebox(0,0)[lb]{\smash{{\SetFigFont{5}{6.0}{\familydefault}{\mddefault}{\updefault}{\color[rgb]{0,0,0}$\frac{2\pi}{5}$}%
}}}}
\put(1583,-911){\makebox(0,0)[lb]{\smash{{\SetFigFont{5}{6.0}{\familydefault}{\mddefault}{\updefault}{\color[rgb]{0,0,0}$\left[\frac{2\pi}{5},\frac{\pi}{3}\right]$}%
}}}}
\put(3589,-142){\makebox(0,0)[lb]{\smash{{\SetFigFont{8}{9.6}{\familydefault}{\mddefault}{\updefault}{\color[rgb]{0,0,0}$v_3^{\rm tr}$}%
}}}}
\put(2027,-482){\makebox(0,0)[lb]{\smash{{\SetFigFont{5}{6.0}{\familydefault}{\mddefault}{\updefault}{\color[rgb]{0,0,0}$\left[\frac{2\pi}{5},\frac{\pi}{6}\right]$}%
}}}}
\put(1500,-563){\makebox(0,0)[lb]{\smash{{\SetFigFont{8}{9.6}{\familydefault}{\mddefault}{\updefault}{\color[rgb]{0,0,0}$v_1^{\rm tr}$}%
}}}}
\put(2391,-306){\makebox(0,0)[lb]{\smash{{\SetFigFont{8}{9.6}{\familydefault}{\mddefault}{\updefault}{\color[rgb]{0,0,0}$v_2^{\rm tr}$}%
}}}}
\put(3093,-982){\makebox(0,0)[lb]{\smash{{\SetFigFont{5}{6.0}{\familydefault}{\mddefault}{\updefault}{\color[rgb]{0,0,0}$\left[\frac{2\pi}{5},\frac{\pi}{4}\right]$}%
}}}}
\put(2944,261){\makebox(0,0)[lb]{\smash{{\SetFigFont{5}{6.0}{\familydefault}{\mddefault}{\updefault}{\color[rgb]{0,0,0}$\left[\frac{2\pi}{5},\frac{\pi}{4}\right]$}%
}}}}
\put(1698,-179){\makebox(0,0)[lb]{\smash{{\SetFigFont{5}{6.0}{\familydefault}{\mddefault}{\updefault}{\color[rgb]{0,0,0}$\left[\frac{2\pi}{5},\frac{\pi}{4}\right]$}%
}}}}
\put(2271, 49){\makebox(0,0)[lb]{\smash{{\SetFigFont{5}{6.0}{\familydefault}{\mddefault}{\updefault}{\color[rgb]{0,0,0}$\left[\frac{2\pi}{5},\frac{\pi}{4}\right]$}%
}}}}
\end{picture}%
\end{center}

We constructed a polyhedron $\widetilde{P}$ realizing the pair $(\widetilde{C},{\bf a}(0))$ and used Newton's
method to deform the faces so that the dihedral angles follow the path ${\bf a}(t)$.
After obtaining a non-compact polyhedron $\widetilde{P}_1$ realizing angles ${\bf a}(1)$ we truncated
the vertices $v_1^{\rm tr},v_2^{\rm tr},v_3^{\rm tr}$.  The final result is the polyhedron:

\begin{center}
\includegraphics[scale=1.0]{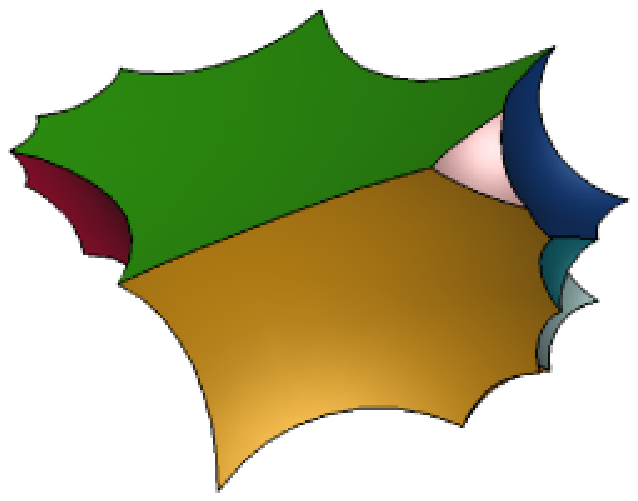}
\end{center}
which will be used in the next section to form part of a compound polyhedron.

\subsection{Constructing compound polyhedra}\label{SEC_COMPOUND}
Any compound polyhedron can be constructed by gluing together a
finite number of truncated polyhedra.  We illustrate this construction
for the polyhedron shown in Section \ref{SEC_EXAMPLE}.

\begin{figure}[h!]
\begin{center}

\begin{picture}(0,0)%
\epsfig{file=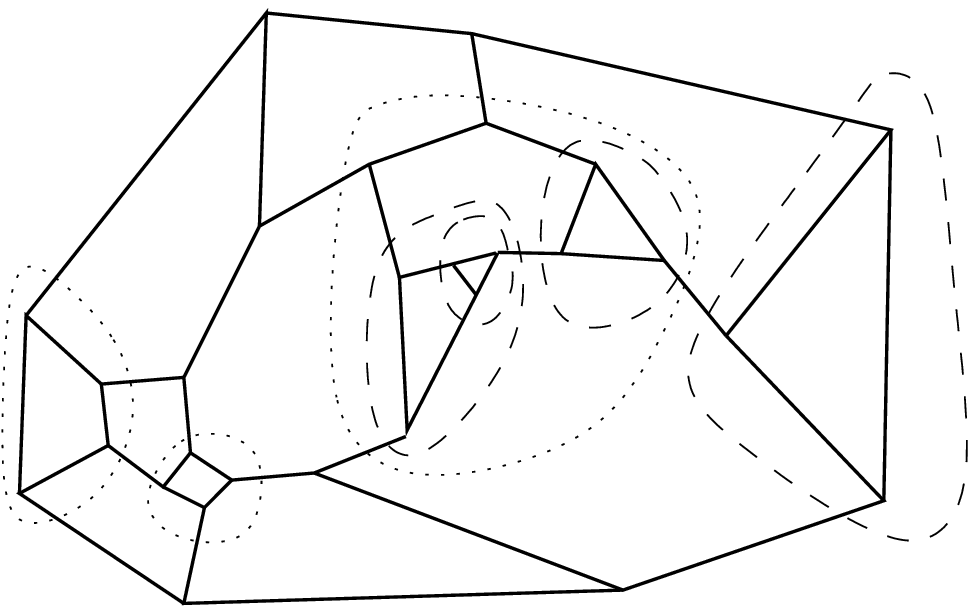}%
\end{picture}%
\setlength{\unitlength}{3947sp}%
\begingroup\makeatletter\ifx\SetFigFont\undefined%
\gdef\SetFigFont#1#2#3#4#5{%
  \reset@font\fontsize{#1}{#2pt}%
  \fontfamily{#3}\fontseries{#4}\fontshape{#5}%
  \selectfont}%
\fi\endgroup%
\begin{picture}(4654,3041)(82,-2275)
\put(2476,-1261){\makebox(0,0)[lb]{\smash{{\SetFigFont{8}{9.6}{\familydefault}{\mddefault}{\updefault}{\color[rgb]{0,0,0}$\gamma$}%
}}}}
\put(581,-98){\makebox(0,0)[lb]{\smash{{\SetFigFont{5}{6.0}{\familydefault}{\mddefault}{\updefault}{\color[rgb]{0,0,0}$\frac{2\pi}{5}$}%
}}}}
\put(1101,-1498){\makebox(0,0)[lb]{\smash{{\SetFigFont{5}{6.0}{\familydefault}{\mddefault}{\updefault}{\color[rgb]{0,0,0}$\frac{2\pi}{5}$}%
}}}}
\put(720,-1049){\makebox(0,0)[lb]{\smash{{\SetFigFont{5}{6.0}{\familydefault}{\mddefault}{\updefault}{\color[rgb]{0,0,0}$\frac{2\pi}{5}$}%
}}}}
\put(1195,-729){\makebox(0,0)[lb]{\smash{{\SetFigFont{5}{6.0}{\familydefault}{\mddefault}{\updefault}{\color[rgb]{0,0,0}$\frac{2\pi}{5}$}%
}}}}
\put(1389,211){\makebox(0,0)[lb]{\smash{{\SetFigFont{5}{6.0}{\familydefault}{\mddefault}{\updefault}{\color[rgb]{0,0,0}$\frac{\pi}{2}$}%
}}}}
\put(3142,-308){\makebox(0,0)[lb]{\smash{{\SetFigFont{5}{6.0}{\familydefault}{\mddefault}{\updefault}{\color[rgb]{0,0,0}$\frac{\pi}{2}$}%
}}}}
\put(2716,-281){\makebox(0,0)[lb]{\smash{{\SetFigFont{5}{6.0}{\familydefault}{\mddefault}{\updefault}{\color[rgb]{0,0,0}$\frac{\pi}{2}$}%
}}}}
\put(3253,416){\makebox(0,0)[lb]{\smash{{\SetFigFont{5}{6.0}{\familydefault}{\mddefault}{\updefault}{\color[rgb]{0,0,0}$\frac{\pi}{4}$}%
}}}}
\put(2091,-81){\makebox(0,0)[lb]{\smash{{\SetFigFont{5}{6.0}{\familydefault}{\mddefault}{\updefault}{\color[rgb]{0,0,0}$\frac{2\pi}{5}$}%
}}}}
\put(3892,-560){\makebox(0,0)[lb]{\smash{{\SetFigFont{5}{6.0}{\familydefault}{\mddefault}{\updefault}{\color[rgb]{0,0,0}$\frac{\pi}{2}$}%
}}}}
\put(4381,-779){\makebox(0,0)[lb]{\smash{{\SetFigFont{5}{6.0}{\familydefault}{\mddefault}{\updefault}{\color[rgb]{0,0,0}$\frac{\pi}{2}$}%
}}}}
\put(3897,-1181){\makebox(0,0)[lb]{\smash{{\SetFigFont{5}{6.0}{\familydefault}{\mddefault}{\updefault}{\color[rgb]{0,0,0}$\frac{\pi}{2}$}%
}}}}
\put(3661,-1985){\makebox(0,0)[lb]{\smash{{\SetFigFont{5}{6.0}{\familydefault}{\mddefault}{\updefault}{\color[rgb]{0,0,0}$\frac{\pi}{4}$}%
}}}}
\put(2411,-1816){\makebox(0,0)[lb]{\smash{{\SetFigFont{5}{6.0}{\familydefault}{\mddefault}{\updefault}{\color[rgb]{0,0,0}$\frac{2\pi}{5}$}%
}}}}
\put(1786,-2248){\makebox(0,0)[lb]{\smash{{\SetFigFont{5}{6.0}{\familydefault}{\mddefault}{\updefault}{\color[rgb]{0,0,0}$\frac{2\pi}{5}$}%
}}}}
\put(1060,-1963){\makebox(0,0)[lb]{\smash{{\SetFigFont{5}{6.0}{\familydefault}{\mddefault}{\updefault}{\color[rgb]{0,0,0}$\frac{2\pi}{5}$}%
}}}}
\put(1141,-1696){\makebox(0,0)[lb]{\smash{{\SetFigFont{5}{6.0}{\familydefault}{\mddefault}{\updefault}{\color[rgb]{0,0,0}$\frac{2\pi}{5}$}%
}}}}
\put(945,-1578){\makebox(0,0)[lb]{\smash{{\SetFigFont{5}{6.0}{\familydefault}{\mddefault}{\updefault}{\color[rgb]{0,0,0}$\frac{2\pi}{5}$}%
}}}}
\put(1018,-1263){\makebox(0,0)[lb]{\smash{{\SetFigFont{5}{6.0}{\familydefault}{\mddefault}{\updefault}{\color[rgb]{0,0,0}$\frac{2\pi}{5}$}%
}}}}
\put(706,-1454){\makebox(0,0)[lb]{\smash{{\SetFigFont{5}{6.0}{\familydefault}{\mddefault}{\updefault}{\color[rgb]{0,0,0}$\frac{2\pi}{5}$}%
}}}}
\put(615,-1271){\makebox(0,0)[lb]{\smash{{\SetFigFont{5}{6.0}{\familydefault}{\mddefault}{\updefault}{\color[rgb]{0,0,0}$\frac{2\pi}{5}$}%
}}}}
\put(400,-1976){\makebox(0,0)[lb]{\smash{{\SetFigFont{5}{6.0}{\familydefault}{\mddefault}{\updefault}{\color[rgb]{0,0,0}$\frac{2\pi}{5}$}%
}}}}
\put(3387,-645){\makebox(0,0)[lb]{\smash{{\SetFigFont{5}{6.0}{\familydefault}{\mddefault}{\updefault}{\color[rgb]{0,0,0}$\frac{\pi}{4}$}%
}}}}
\put(401,-1604){\makebox(0,0)[lb]{\smash{{\SetFigFont{5}{6.0}{\familydefault}{\mddefault}{\updefault}{\color[rgb]{0,0,0}$\frac{2\pi}{5}$}%
}}}}
\put(1577,-309){\makebox(0,0)[lb]{\smash{{\SetFigFont{5}{6.0}{\familydefault}{\mddefault}{\updefault}{\color[rgb]{0,0,0}$\frac{2\pi}{5}$}%
}}}}
\put(1953,-292){\makebox(0,0)[lb]{\smash{{\SetFigFont{5}{6.0}{\familydefault}{\mddefault}{\updefault}{\color[rgb]{0,0,0}$\frac{\pi}{4}$}%
}}}}
\put(2436,-622){\makebox(0,0)[lb]{\smash{{\SetFigFont{5}{6.0}{\familydefault}{\mddefault}{\updefault}{\color[rgb]{0,0,0}$\frac{\pi}{2}$}%
}}}}
\put(2041,-862){\makebox(0,0)[lb]{\smash{{\SetFigFont{5}{6.0}{\familydefault}{\mddefault}{\updefault}{\color[rgb]{0,0,0}$\frac{\pi}{2}$}%
}}}}
\put(1374,-1660){\makebox(0,0)[lb]{\smash{{\SetFigFont{5}{6.0}{\familydefault}{\mddefault}{\updefault}{\color[rgb]{0,0,0}$\frac{2\pi}{5}$}%
}}}}
\put(225,-1275){\makebox(0,0)[lb]{\smash{{\SetFigFont{5}{6.0}{\familydefault}{\mddefault}{\updefault}{\color[rgb]{0,0,0}$\frac{2\pi}{5}$}%
}}}}
\put(1697,694){\makebox(0,0)[lb]{\smash{{\SetFigFont{5}{6.0}{\familydefault}{\mddefault}{\updefault}{\color[rgb]{0,0,0}$\frac{2\pi}{5}$}%
}}}}
\put(2661, 90){\makebox(0,0)[lb]{\smash{{\SetFigFont{5}{6.0}{\familydefault}{\mddefault}{\updefault}{\color[rgb]{0,0,0}$\frac{\pi}{4}$}%
}}}}
\put(2407,360){\makebox(0,0)[lb]{\smash{{\SetFigFont{5}{6.0}{\familydefault}{\mddefault}{\updefault}{\color[rgb]{0,0,0}$\frac{2\pi}{5}$}%
}}}}
\put(402,-932){\makebox(0,0)[lb]{\smash{{\SetFigFont{5}{6.0}{\familydefault}{\mddefault}{\updefault}{\color[rgb]{0,0,0}$\frac{2\pi}{5}$}%
}}}}
\put(2241,-1003){\makebox(0,0)[lb]{\smash{{\SetFigFont{5}{6.0}{\familydefault}{\mddefault}{\updefault}{\color[rgb]{0,0,0}$\frac{\pi}{3}$}%
}}}}
\put(2086,-504){\makebox(0,0)[lb]{\smash{{\SetFigFont{5}{6.0}{\familydefault}{\mddefault}{\updefault}{\color[rgb]{0,0,0}$\frac{\pi}{3}$}%
}}}}
\put(2972,-599){\makebox(0,0)[lb]{\smash{{\SetFigFont{5}{6.0}{\familydefault}{\mddefault}{\updefault}{\color[rgb]{0,0,0}$\frac{\pi}{2}$}%
}}}}
\put(1653,-1430){\makebox(0,0)[lb]{\smash{{\SetFigFont{5}{6.0}{\familydefault}{\mddefault}{\updefault}{\color[rgb]{0,0,0}$\frac{\pi}{3}$}%
}}}}
\put(2219,-670){\makebox(0,0)[lb]{\smash{{\SetFigFont{5}{6.0}{\familydefault}{\mddefault}{\updefault}{\color[rgb]{0,0,0}$\frac{\pi}{2}$}%
}}}}
\put(2276,-451){\makebox(0,0)[lb]{\smash{{\SetFigFont{5}{6.0}{\familydefault}{\mddefault}{\updefault}{\color[rgb]{0,0,0}$\frac{\pi}{2}$}%
}}}}
\put(2569,-406){\makebox(0,0)[lb]{\smash{{\SetFigFont{5}{6.0}{\familydefault}{\mddefault}{\updefault}{\color[rgb]{0,0,0}$\frac{\pi}{6}$}%
}}}}
\put(878,-1751){\makebox(0,0)[lb]{\smash{{\SetFigFont{5}{6.0}{\familydefault}{\mddefault}{\updefault}{\color[rgb]{0,0,0}$\frac{2\pi}{5}$}%
}}}}
\end{picture}%
\end{center}
\caption{\hbox{ }}
\end{figure}

In general, one cuts along every prismatic $3$-circuit which does not
correspond to a triangular face.  Here there is one such circuit which is
labeled $\gamma$.  We cut along $\gamma$ obtaining two combinatorial polyhedra
for which every prismatic $3$-circuit corresponds to a triangular face:

\begin{center}

\begin{picture}(0,0)%
\epsfig{file=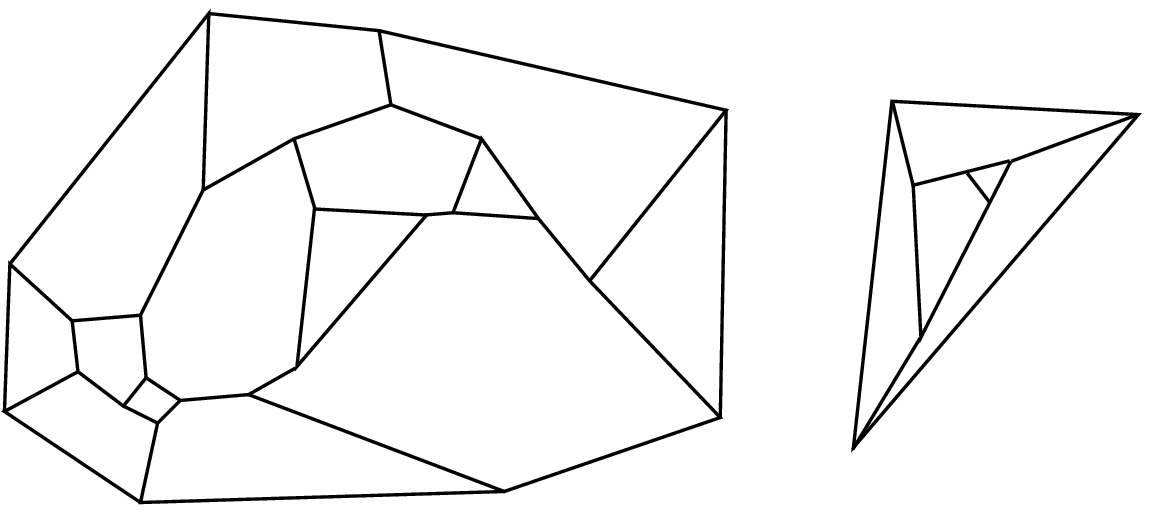}%
\end{picture}%
\setlength{\unitlength}{3947sp}%
\begingroup\makeatletter\ifx\SetFigFont\undefined%
\gdef\SetFigFont#1#2#3#4#5{%
  \reset@font\fontsize{#1}{#2pt}%
  \fontfamily{#3}\fontseries{#4}\fontshape{#5}%
  \selectfont}%
\fi\endgroup%
\begin{picture}(5487,2581)(137,-1763)
\put(2234,251){\makebox(0,0)[lb]{\smash{{\SetFigFont{5}{6.0}{\familydefault}{\mddefault}{\updefault}{\color[rgb]{0,0,0}$\frac{\pi}{4}$}%
}}}}
\put(356,-617){\makebox(0,0)[lb]{\smash{{\SetFigFont{5}{6.0}{\familydefault}{\mddefault}{\updefault}{\color[rgb]{0,0,0}$\frac{2\pi}{5}$}%
}}}}
\put(609,-725){\makebox(0,0)[lb]{\smash{{\SetFigFont{5}{6.0}{\familydefault}{\mddefault}{\updefault}{\color[rgb]{0,0,0}$\frac{2\pi}{5}$}%
}}}}
\put(1003,-460){\makebox(0,0)[lb]{\smash{{\SetFigFont{5}{6.0}{\familydefault}{\mddefault}{\updefault}{\color[rgb]{0,0,0}$\frac{2\pi}{5}$}%
}}}}
\put(1163,319){\makebox(0,0)[lb]{\smash{{\SetFigFont{5}{6.0}{\familydefault}{\mddefault}{\updefault}{\color[rgb]{0,0,0}$\frac{\pi}{2}$}%
}}}}
\put(2615,-111){\makebox(0,0)[lb]{\smash{{\SetFigFont{5}{6.0}{\familydefault}{\mddefault}{\updefault}{\color[rgb]{0,0,0}$\frac{\pi}{2}$}%
}}}}
\put(2262,-89){\makebox(0,0)[lb]{\smash{{\SetFigFont{5}{6.0}{\familydefault}{\mddefault}{\updefault}{\color[rgb]{0,0,0}$\frac{\pi}{2}$}%
}}}}
\put(1745, 77){\makebox(0,0)[lb]{\smash{{\SetFigFont{5}{6.0}{\familydefault}{\mddefault}{\updefault}{\color[rgb]{0,0,0}$\frac{\pi}{2}$}%
}}}}
\put(3236,-320){\makebox(0,0)[lb]{\smash{{\SetFigFont{5}{6.0}{\familydefault}{\mddefault}{\updefault}{\color[rgb]{0,0,0}$\frac{\pi}{2}$}%
}}}}
\put(3641,-501){\makebox(0,0)[lb]{\smash{{\SetFigFont{5}{6.0}{\familydefault}{\mddefault}{\updefault}{\color[rgb]{0,0,0}$\frac{\pi}{2}$}%
}}}}
\put(3240,-834){\makebox(0,0)[lb]{\smash{{\SetFigFont{5}{6.0}{\familydefault}{\mddefault}{\updefault}{\color[rgb]{0,0,0}$\frac{\pi}{2}$}%
}}}}
\put(3045,-1500){\makebox(0,0)[lb]{\smash{{\SetFigFont{5}{6.0}{\familydefault}{\mddefault}{\updefault}{\color[rgb]{0,0,0}$\frac{\pi}{4}$}%
}}}}
\put(2010,-1360){\makebox(0,0)[lb]{\smash{{\SetFigFont{5}{6.0}{\familydefault}{\mddefault}{\updefault}{\color[rgb]{0,0,0}$\frac{2\pi}{5}$}%
}}}}
\put(958,-1261){\makebox(0,0)[lb]{\smash{{\SetFigFont{5}{6.0}{\familydefault}{\mddefault}{\updefault}{\color[rgb]{0,0,0}$\frac{2\pi}{5}$}%
}}}}
\put(796,-1163){\makebox(0,0)[lb]{\smash{{\SetFigFont{5}{6.0}{\familydefault}{\mddefault}{\updefault}{\color[rgb]{0,0,0}$\frac{2\pi}{5}$}%
}}}}
\put(856,-902){\makebox(0,0)[lb]{\smash{{\SetFigFont{5}{6.0}{\familydefault}{\mddefault}{\updefault}{\color[rgb]{0,0,0}$\frac{2\pi}{5}$}%
}}}}
\put(522,-909){\makebox(0,0)[lb]{\smash{{\SetFigFont{5}{6.0}{\familydefault}{\mddefault}{\updefault}{\color[rgb]{0,0,0}$\frac{2\pi}{5}$}%
}}}}
\put(2818,-390){\makebox(0,0)[lb]{\smash{{\SetFigFont{5}{6.0}{\familydefault}{\mddefault}{\updefault}{\color[rgb]{0,0,0}$\frac{\pi}{4}$}%
}}}}
\put(1319,-112){\makebox(0,0)[lb]{\smash{{\SetFigFont{5}{6.0}{\familydefault}{\mddefault}{\updefault}{\color[rgb]{0,0,0}$\frac{\pi}{2}$}%
}}}}
\put(1630,-98){\makebox(0,0)[lb]{\smash{{\SetFigFont{5}{6.0}{\familydefault}{\mddefault}{\updefault}{\color[rgb]{0,0,0}$\frac{\pi}{4}$}%
}}}}
\put(1513,-1120){\makebox(0,0)[lb]{\smash{{\SetFigFont{5}{6.0}{\familydefault}{\mddefault}{\updefault}{\color[rgb]{0,0,0}$\frac{\pi}{3}$}%
}}}}
\put(1881,-729){\makebox(0,0)[lb]{\smash{{\SetFigFont{5}{6.0}{\familydefault}{\mddefault}{\updefault}{\color[rgb]{0,0,0}$\frac{\pi}{2}$}%
}}}}
\put(1139,-1242){\makebox(0,0)[lb]{\smash{{\SetFigFont{5}{6.0}{\familydefault}{\mddefault}{\updefault}{\color[rgb]{0,0,0}$\frac{2\pi}{5}$}%
}}}}
\put(1758,-515){\makebox(0,0)[lb]{\smash{{\SetFigFont{8}{9.6}{\familydefault}{\mddefault}{\updefault}{\color[rgb]{0,0,0}$F$}%
}}}}
\put(4774,-496){\makebox(0,0)[lb]{\smash{{\SetFigFont{5}{6.0}{\familydefault}{\mddefault}{\updefault}{\color[rgb]{0,0,0}$\frac{\pi}{3}$}%
}}}}
\put(4955,-148){\makebox(0,0)[lb]{\smash{{\SetFigFont{5}{6.0}{\familydefault}{\mddefault}{\updefault}{\color[rgb]{0,0,0}$\frac{\pi}{2}$}%
}}}}
\put(4780,337){\makebox(0,0)[lb]{\smash{{\SetFigFont{5}{6.0}{\familydefault}{\mddefault}{\updefault}{\color[rgb]{0,0,0}$\frac{\pi}{2}$}%
}}}}
\put(4330,-971){\makebox(0,0)[lb]{\smash{{\SetFigFont{5}{6.0}{\familydefault}{\mddefault}{\updefault}{\color[rgb]{0,0,0}$\frac{\pi}{3}$}%
}}}}
\put(4595,-192){\makebox(0,0)[lb]{\smash{{\SetFigFont{5}{6.0}{\familydefault}{\mddefault}{\updefault}{\color[rgb]{0,0,0}$\frac{\pi}{3}$}%
}}}}
\put(4514, 66){\makebox(0,0)[lb]{\smash{{\SetFigFont{5}{6.0}{\familydefault}{\mddefault}{\updefault}{\color[rgb]{0,0,0}$\frac{\pi}{3}$}%
}}}}
\put(4800,-979){\makebox(0,0)[lb]{\smash{{\SetFigFont{8}{9.6}{\familydefault}{\mddefault}{\updefault}{\color[rgb]{0,0,0}$\tilde{F}$}%
}}}}
\put(4705,-1116){\makebox(0,0)[lb]{\smash{{\SetFigFont{8}{9.6}{\familydefault}{\mddefault}{\updefault}{\color[rgb]{0,0,0}{\small on the outside}}%
}}}}
\put(4894,-747){\makebox(0,0)[lb]{\smash{{\SetFigFont{5}{6.0}{\familydefault}{\mddefault}{\updefault}{\color[rgb]{0,0,0}$\frac{\pi}{2}$}%
}}}}
\put(4846, 26){\makebox(0,0)[lb]{\smash{{\SetFigFont{5}{6.0}{\familydefault}{\mddefault}{\updefault}{\color[rgb]{0,0,0}$\frac{\pi}{2}$}%
}}}}
\put(5192,-27){\makebox(0,0)[lb]{\smash{{\SetFigFont{5}{6.0}{\familydefault}{\mddefault}{\updefault}{\color[rgb]{0,0,0}$\frac{\pi}{6}$}%
}}}}
\put(4144,-497){\makebox(0,0)[lb]{\smash{{\SetFigFont{5}{6.0}{\familydefault}{\mddefault}{\updefault}{\color[rgb]{0,0,0}$\frac{\pi}{2}$}%
}}}}
\put(4558,-441){\makebox(0,0)[lb]{\smash{{\SetFigFont{5}{6.0}{\familydefault}{\mddefault}{\updefault}{\color[rgb]{0,0,0}$\frac{\pi}{2}$}%
}}}}
\put(4727,-199){\makebox(0,0)[lb]{\smash{{\SetFigFont{5}{6.0}{\familydefault}{\mddefault}{\updefault}{\color[rgb]{0,0,0}$\frac{\pi}{2}$}%
}}}}
\put(2471,-361){\makebox(0,0)[lb]{\smash{{\SetFigFont{5}{6.0}{\familydefault}{\mddefault}{\updefault}{\color[rgb]{0,0,0}$\frac{\pi}{2}$}%
}}}}
\put(1501,-1736){\makebox(0,0)[lb]{\smash{{\SetFigFont{5}{6.0}{\familydefault}{\mddefault}{\updefault}{\color[rgb]{0,0,0}$\frac{2\pi}{5}$}%
}}}}
\put(326,-1511){\makebox(0,0)[lb]{\smash{{\SetFigFont{5}{6.0}{\familydefault}{\mddefault}{\updefault}{\color[rgb]{0,0,0}$\frac{2\pi}{5}$}%
}}}}
\put(882,-1482){\makebox(0,0)[lb]{\smash{{\SetFigFont{5}{6.0}{\familydefault}{\mddefault}{\updefault}{\color[rgb]{0,0,0}$\frac{2\pi}{5}$}%
}}}}
\put(695,-1321){\makebox(0,0)[lb]{\smash{{\SetFigFont{5}{6.0}{\familydefault}{\mddefault}{\updefault}{\color[rgb]{0,0,0}$\frac{2\pi}{5}$}%
}}}}
\put(336,-1217){\makebox(0,0)[lb]{\smash{{\SetFigFont{5}{6.0}{\familydefault}{\mddefault}{\updefault}{\color[rgb]{0,0,0}$\frac{2\pi}{5}$}%
}}}}
\put(476, 68){\makebox(0,0)[lb]{\smash{{\SetFigFont{5}{6.0}{\familydefault}{\mddefault}{\updefault}{\color[rgb]{0,0,0}$\frac{2\pi}{5}$}%
}}}}
\put(2237,-352){\makebox(0,0)[lb]{\smash{{\SetFigFont{5}{6.0}{\familydefault}{\mddefault}{\updefault}{\color[rgb]{0,0,0}$\frac{\pi}{6}$}%
}}}}
\put(1490,-559){\makebox(0,0)[lb]{\smash{{\SetFigFont{5}{6.0}{\familydefault}{\mddefault}{\updefault}{\color[rgb]{0,0,0}$\frac{\pi}{2}$}%
}}}}
\put(190,-912){\makebox(0,0)[lb]{\smash{{\SetFigFont{5}{6.0}{\familydefault}{\mddefault}{\updefault}{\color[rgb]{0,0,0}$\frac{2\pi}{5}$}%
}}}}
\put(1458,746){\makebox(0,0)[lb]{\smash{{\SetFigFont{5}{6.0}{\familydefault}{\mddefault}{\updefault}{\color[rgb]{0,0,0}$\frac{2\pi}{5}$}%
}}}}
\put(2711,525){\makebox(0,0)[lb]{\smash{{\SetFigFont{5}{6.0}{\familydefault}{\mddefault}{\updefault}{\color[rgb]{0,0,0}$\frac{\pi}{4}$}%
}}}}
\put(1856,-177){\makebox(0,0)[lb]{\smash{{\SetFigFont{5}{6.0}{\familydefault}{\mddefault}{\updefault}{\color[rgb]{0,0,0}$\frac{\pi}{2}$}%
}}}}
\put(598,-1060){\makebox(0,0)[lb]{\smash{{\SetFigFont{5}{6.0}{\familydefault}{\mddefault}{\updefault}{\color[rgb]{0,0,0}$\frac{2\pi}{5}$}%
}}}}
\put(921,-1070){\makebox(0,0)[lb]{\smash{{\SetFigFont{5}{6.0}{\familydefault}{\mddefault}{\updefault}{\color[rgb]{0,0,0}$\frac{2\pi}{5}$}%
}}}}
\end{picture}%
\end{center}

In this case, the diagram on the
left is that for the polyhedron that
we constructed in the previous section.  The diagram on the right is that of 
the truncated triangular prism, which can also be easily constructed.

We require that the new triangular faces $F$ and $\widetilde{F}$ obtained by
cutting along $\gamma$ be perpendicular to each of the other faces that they
intersect.  Then each face angle equals the dihedral angle outside of $F$, or
$\widetilde{F}$, that leads to that vertex.  Because we obtained the two diagrams by
cutting the original diagram along $\gamma$, the dihedral angles on the edges leading to $F$ and
$\widetilde{F}$ are the same and we naturally obtain that $F$ and $\widetilde{F}$
have the same face angles, but are mirror images of each other.

\begin{center}

\begin{picture}(0,0)%
\epsfig{file=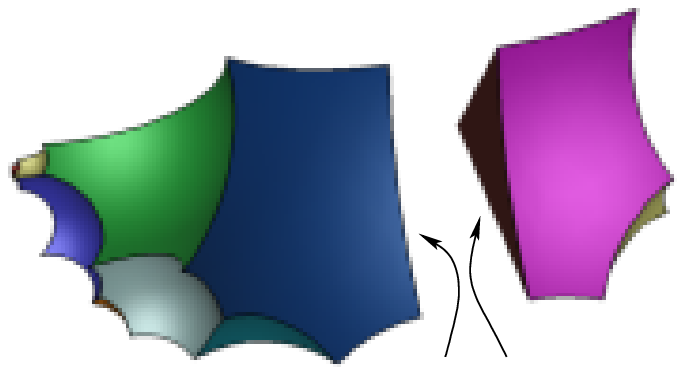}%
\end{picture}%
\setlength{\unitlength}{3947sp}%
\begingroup\makeatletter\ifx\SetFigFont\undefined%
\gdef\SetFigFont#1#2#3#4#5{%
  \reset@font\fontsize{#1}{#2pt}%
  \fontfamily{#3}\fontseries{#4}\fontshape{#5}%
  \selectfont}%
\fi\endgroup%
\begin{picture}(3375,2008)(1201,-2369)
\put(3601,-2311){\makebox(0,0)[lb]{\smash{{\SetFigFont{12}{14.4}{\familydefault}{\mddefault}{\updefault}{\color[rgb]{0,0,0}$\tilde{F}$}%
}}}}
\put(3226,-2311){\makebox(0,0)[lb]{\smash{{\SetFigFont{12}{14.4}{\familydefault}{\mddefault}{\updefault}{\color[rgb]{0,0,0}$F$}%
}}}}
\end{picture}%
\end{center}

These two polyhedra glue perfectly together to form a polyhedron realizing
$(C,{\bf a})$ as shown in the following figure.

\begin{center}
\includegraphics[scale=1.0]{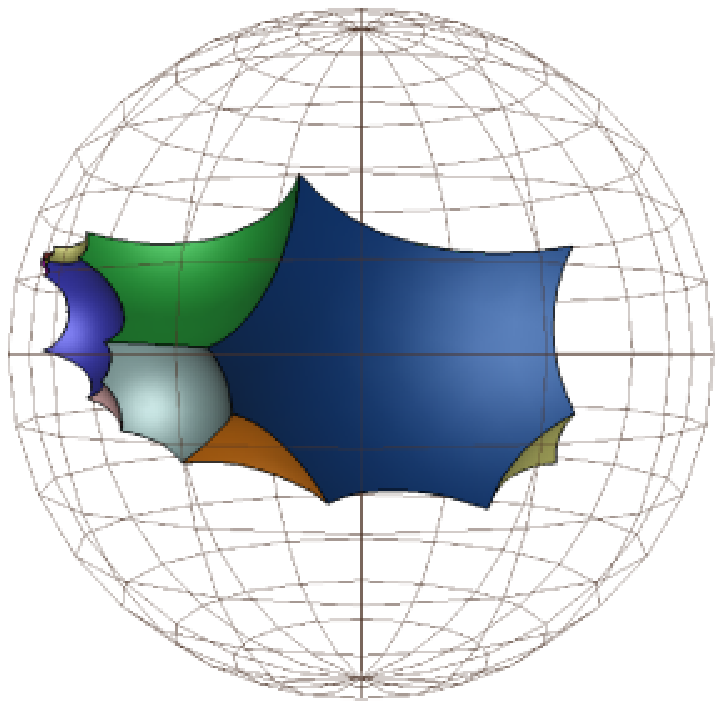}
\end{center}

\section{Applications to discrete groups and polyhedral orbifolds}\label{SEC_APPLICATIONS}
Let $P$ be a finite volume hyperbolic polyhedron having dihedral angles each of
which is a proper integer sub-multiple of $\pi$.  It is a well-known
application of the Poincar\'e Polyhedron Theorem \cite{POINCARE} that the group
generated by reflections in the faces $P$ forms a discrete subgroup $\Gamma_P$ of
$Isom(\mathbb{H}^3)$.  Such groups have been extensively studied, see
\cite{VINREFL}, and the references therein.  

Given such a discrete reflection group $\Gamma_P$, we denote the corresponding
orbifold by $O_P = \mathbb{H}^3/\Gamma_P$.  We will use the term {\em
polyhedral orbifolds} to describe orbifolds obtained in this way.  (Note: often
in the literature, the term ``polyhedral orbifold'' is used to describe the
oriented double cover $\mathbb{H}^3/\Gamma^+_P$, where $\Gamma^+_P$ 
is index two subgroup consisting of orientation preserving elements of
$\Gamma_P$.)  See Thurston \cite[Chapter 13]{TH} and Reni \cite{RENI} for more
details on polyhedral orbifolds.

We use the computer program described in this paper to construct
examples from three classes of polyhedral orbifolds, the Lambert
cubes \cite{KELLERHALS}, the L\"obell orbifolds \cite{LOBELL,VESNIN_LOB,HYPER_ELL}, and
an  mysterious orbifold described by Mednykh and Vesnin \cite{HYPER_ELL} whose
16-fold cover is a``hyperelliptic'' compact hyperbolic manifold.

We output the generators each reflection group as elements of $SO(3,1)$ into
SnapPea \cite{SNAPPEA}, computing volumes and length spectra of these
orbifolds.  For details on how SnapPea calculates the length spectrum refer to
\cite{WEEKS_LENGTHS}.

\subsection{Construction of Lambert Cubes}
A Lambert cube is a compact polyhedron realizing the combinatorial type of a cube,
with three non-coplanar edges chosen and assigned dihedral angles $\alpha,$ $\beta$, and $\gamma$,
and the remaining edges assigned dihedral angles $\frac{\pi}{2}$.
It is easy to verify that if $0 <
\alpha, \beta, \gamma < \frac{\pi}{2}$ then, such an 
assignment of dihedral angles  satisfies the hypotheses of Andreev's Theorem.  The resulting
polyhedron is called the $(\alpha,\beta,\gamma)$ Lambert Cube, which we will
denote by $P_{\alpha,\beta,\gamma}$.  Thus, there are discrete reflection groups
generated in the faces of a Lambert Cube when $\alpha = \frac{\pi}{p}$, $\beta = \frac{\pi}{q}$, and $\gamma =
\frac{\pi}{r}$ for integers $p,q,r > 2$.  We denote the corresponding orbifold $O_{\rm Lambert}(p,q,r)$  In the following table, we
present volumes and the lengths of the shortest geodesics for a sampling of
Lambert Cubes for small $p,q$, and $r$:

\vspace{.1in}
\begin{tabular}{|l|l|l|} 
\hline
$O_{\rm Lambert}(3,3,3)$ & Computed Volume: $0.324423$ & Theoretical Volume: $0.3244234492$ \\
\hline
Short Geodesics & 6 mI $1.087070$ & 3 mI  $1.087070  +  i \cot 2\pi/3$  \\
 & 6 mI $1.400257$  & 6 mI $1.400257 + i \cdot \pi$  \\
 & 3 mI $1.601733  +  i \cdot 2.765750$ & 3 mI  $1.790480  +  i \cdot  0.762413$ \\
 & 6 mI  $1.864162$ & 3 mI  $2.138622$  \\
 & 4 mI  $2.174140$ &  6 mI  $2.199243  +  i \cdot  2.436822$ \\
\hline
\hline
$O_{\rm Lambert}(3,4,5)$ & Computed Volume: $0.479079$ & Theoretical Volume: $0.4790790206$ \\
\hline
Short Geodesics & 2 mI  $0.622685$ & 1 mI  $0.622685  +  i \cdot  1.256637$      \\
 & 1 mI  $0.622685  +  i \cdot 2.513274$ & 3 mI  $0.883748$  \\
 &  1 mI  $0.883748  +  i \cdot 1.570797$ & 1 mI  $0.883748  +  i \cdot 3.141592$ \\
 & 1 mI  $1.123387$   &   1 mI  $1.123387  +  i \cdot 3.141593$ \\
 & 1 mI  $1.245371$     & \\
\hline
\hline
$O_{\rm Lambert}(4,4,4)$ & Computed Volume: $0.554152$ & Theoretical Volume: $0.5382759501$ \\
\hline
Short Geodesics & 2 mI  $0.175240$  & 1 mI  $0.175240  +  i \cdot  0.369599$ \\
 & 1 mI  $0.175240  +  i \cdot 1.108797$ & 1 mI  $0.175240  +  i \cdot 0.739198$\\
 & 1 mI  $0.175240  +  i \cdot 0.739198$ &  1 mI  $0.175240  +  i \cdot 1.478396$ \\
 & 1 mI  $0.175240  +  i \cdot  1.847996$ &  1 mI  $0.175240  +  i \cdot 2.217595$   \\
 & 1 mI  $0.175240  +  i \cdot  2.587194$ & 1 mI  $0.175240  +  i \cdot 2.956793$ \\
 & 1 mI  $0.350479$ &  \\  
\hline
\hline
$O_{\rm Lambert}(5,8,12)$ & Computed Volume: $0.768801$ & Theoretical Volume: $0.7688005863$ \\
\hline
Short Geodesics & 3 mI  $0.407809$ & 1 mI  $0.407809 +  i \cdot 0.523599$ \\
 & 1 mI  $0.407809  +  i \cdot 1.047198$ & 1 mI  $0.407809  +  i \cdot  1.570797$ \\
 &  1 mI  $0.407809  +  i \cdot 2.094396$ & 1 mI  $0.407809  +  i \cdot 2.617995$ \\
 & 1 mI  $0.407809  +  i \cdot 3.141592$  & 2 mI $0.643110$ \\
 & 1 mI $0.643110 + i \cdot 0.785398$ & 1 mI $0.643110 + i \cdot 1.570796$ \\
 & 1 mI $0.643110 + i \cdot 2.356194$ & 1 mI $0.643110 + i \cdot 3.141593$ \\
\hline
\end{tabular} 

\vspace{.1in} The format of the lists of Geodesic lengths presented in this and
in the following tables, is the same as that presented by SnapPea.  The first
entry is the multiplicity of distinct geodesics having the same complex length.
The second entry is either ``mI'' to indicate that the geodesic has the
topological type of a mirrored interval or is empty, if the geodesic has the
topological type of a circle.  The third entry is the complex length.  Nearly all
of the short geodesics that we present in these tables are mirrored intervals
because our orbifolds are mirrored polyhedra and because we have only listed
rather short geodesics.

Also notice, that while SnapPea provides many more digits of precision for the
geodesic length, we have rounded to the first 6 decimal places in order to
group geodesics that are likely to correspond to the same class, but weren't
listed that way due to numerical imprecision.  \vspace{.1in}

The volumes of Lambert cubes have been explicitly calculated by R. Kellerhals \cite{KELLERHALS}.
If we write $\Delta(\eta,\xi) = \Lambda(\eta+\xi) - \Lambda(\eta-\xi)$, where $\Lambda$ is the well-known Lobachevskii function $\Lambda(x) = -\int_0^x \log \vert 2\sin(t) \vert dt$, then

\begin{eqnarray}\label{LAM_VOL}
{\rm Vol}(P_{\alpha,\beta,\gamma}) = \frac{1}{4}\left(\Delta(\alpha,\theta)+\Delta(\beta,\theta)+\Delta(\gamma,\theta)-2 \cdot \Delta\left(\frac{\pi}{2},\theta\right) -\Delta(0,\theta)\right).
\end{eqnarray}

where $\theta$, with $0 < \theta < \frac{\pi}{2}$ is the parameter defined by:
\begin{eqnarray*}
{\rm tan}^2(\theta) = p + \sqrt{p^2+L^2M^2N^2},\\
p = \frac{L^2+M^2+N^2+1}{2} \mbox{,  and  }\\
L = {\rm tan}\alpha, M = {\rm tan}\beta, N = {\rm tan}\gamma.
\end{eqnarray*}

The column in the above table labeled ``approximate volume'' gives the volume
of $O_{\rm Lambert}(p,q,r)$ as computed using SnapPea, while the column labeled
``actual volume'' gives the volume of $O_{\rm Lambert}(p,q,r)$ computed using
Equation \ref{LAM_VOL}.

\subsection{Construction of L\"obell Orbifolds}

For each $n > 5$, there is an radially symmetric combinatorial polyhedron having two $n$-sided faces
and having $2n$ faces with $5$ sides which provides a natural generalization of
the dodecahedron.  This combinatorial polyhedron is depicted below for $n=8$.

\begin{center}

\begin{picture}(0,0)%
\epsfig{file=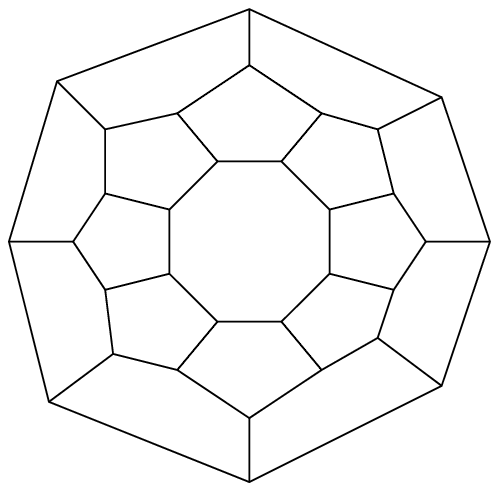}%
\end{picture}%
\setlength{\unitlength}{3947sp}%
\begingroup\makeatletter\ifx\SetFigFont\undefined%
\gdef\SetFigFont#1#2#3#4#5{%
  \reset@font\fontsize{#1}{#2pt}%
  \fontfamily{#3}\fontseries{#4}\fontshape{#5}%
  \selectfont}%
\fi\endgroup%
\begin{picture}(2333,2294)(3889,-4808)
\end{picture}%
\end{center}

Andreev's Theorem provides the existence of a compact right angled polyhedron
$R_n$ realizing this abstract polyhedron because it contains no prismatic
$3$-circuits or prismatic $4$-circuits.  (In fact, the work of L\"obell
predates Andreev by many years, and one can also verify the existence of $R_n$
as an appropriate truncation and gluing of tetrahedra.) We refer to the group
generated by reflections in the faces of $R_n$ by $\Gamma_n$ and the corresponding orbifold
$O_{\rm L\ddot{o}bell}(n) = \mathbb{H}^3/\Gamma_n$.

\vspace{0.05in}
\noindent
{\bf Historical Note:}

\noindent
While we restrict our attention to the orbifold $O_{\rm L\ddot{o}bell}(n)$ in this
paper, the reader may wish to notice that the first example of a closed
hyperbolic manifolds was constructed by L\"obell \cite{LOBELL} in 1931 by an
appropriate gluing of $8$ copies of $R_5$.  Generalizing this notion, Vesnin
\cite{VESNIN_LOB} has described a convenient algebraic method to construct a
torsion free subgroup $\Gamma_n' \subset \Gamma_n$ of index 8.  This, in the
{\em $n$-th L\"obell Manifold} is the compact, orientable, hyperbolic manifold
$M_{\rm L\ddot{o}bell}(n) := \mathbb{H}^3/\Gamma'_n$.  Naturally, $M_{\rm
L\ddot{o}bell}(n)$ is an $8$-fold (orbifold) cover of $O_{\rm L\ddot{o}bell}(n)$.  We
refer the reader to the nice exposition in \cite{VESNIN_LOB,HYPER_ELL} for the
details.  The delightful paper by Reni \cite{RENI} provides further details on the
construction of hyperbolic manifolds and orbifolds and finite covers of right-angled polyhedra.
\vspace{0.05in}

We include the following table of data computed using SnapPea for the $n=5,\cdots,8$ L\"obell orbifolds.

\vspace{.1in}
\begin{tabular}{|l|l|l|} 
\hline
$O_{\rm L\ddot{o}bell}(5)$ (Dodecahedron) & Computed Volume: $4.306208$ & Theoretical Volume: $4.3062076007$ \\
\hline
Short Geodesics & 60 mI  $2.122550$ & 60 mI  $2.122550 +  i \cdot\pi $   \\
  & 60 mI  $2.938703$ & 60 mI   $2.938703 + i \cdot \pi$      \\
 &  126 mI $3.233843$    & 60 mI $3.579641$    \\
 &  60 mI  $3.783112 + i \cdot 1.376928$ & 12 $3.835986$   \\
 &  12  $3.835986 + i \cdot \pi$ & 60 mI  $3.835986 +  i \cdot \pi$  \\
 &  60 mI  $3.966774$ & 60 mI  $4.0270318 + i \cdot 2.264758$ \\
\hline
\hline
$O_{\rm L\ddot{o}bell}(6)$ & Computed Volume: $6.023046$ & Theoretical Volume: $6.0230460200$ \\
\hline
Short Geodesics & 36 mI  $1.762747$ & 12 mI  $1.762747  +  i \cdot \pi$   \\
  & 37 mI  $2.292431$ & 12 mI  $2.292431 + i \cdot \pi$      \\
 &  48 mI  $2.633916$    & 24 mI  $2.633916  +  i \cdot \pi$    \\
 &  36 mI  $2.887271$ & 24 mI  $2.887271  +  i \cdot \pi$   \\
 &  48 mI  $3.088970$ & 12 mI  $3.154720  +  i \cdot  1.312496$  \\
 &  24 mI  $3.256614$ & 36 mI  $3.256614 + i \cdot \pi$ \\
\hline
\hline
$O_{\rm L\ddot{o}bell}(7)$ & Computed Volume: $7.563249$ & Theoretical Volume: $7.5632490914$ \\
\hline
Short Geodesics & 42 mI  1.611051 & 14 mI  $1.611051  +  i \cdot \pi$  \\
  & 1 mI  $1.823106$ &  42 mI  $2.388409$    \\
 & 14 mI  $2.388409  +  i \cdot \pi$  & 14 mI  $2.512394$ \\
 &  14 mI  $2.512394  +  i \cdot \pi$    & 14 mI  $2.601666$\\
 &  70 mI  $2.898149$ & 14 mI $2.898149  +  i \cdot 1.280529$ \\
 & 42 mI $2.898149 + i \cdot \pi$ & 14 mi $3.031090 +i \cdot \pi$ \\
\hline
\hline
$O_{\rm L\ddot{o}bell}(8)$ & Computed Volume: $9.019053$ & Theoretical Volume: $9.0190527274$\\
\hline
Short Geodesics & 49 mI  $1.528571$ & 16 mI  $1.528571  +  i \cdot \pi$ \\
 & 80 mI  $2.448452$ & 32 mI  $2.448452  +  i \cdot \pi$     \\
 & 16 mI $2.760884 + i \cdot 1.261789$ & 32 mI $2.914035$ \\
 & 48 mI $2.914035 + i \cdot \pi$ & 160 mI $3.057142$ \\
 & 32 mI $3.057142 + i \cdot \pi$ & 16 mI $3.461816 + i \cdot 2.650944$ \\
 & 64 mI $3.553688$ & 32 mI $3.553688 + i \cdot \pi$ \\
\hline
\end{tabular} 

\vspace{.1in}

The column labeled ``Computed Volume'' gives the volume as computed in SnapPea, whereas``Theoretical Volume'' provides the volume of $O_{\rm L\ddot{o}bell}(n)$ using explicit formula from \cite{LOB_VOL}.  (In fact we have divided by $8$ the volume formula presented in \cite{LOB_VOL}, because they
study the volume of the $8$-fold cover $M_{\rm L\ddot{o}bell}(n)$.)
If we let $\theta = \frac{\pi}{2} -\arccos\left(\frac{1}{2\cos(\pi/n)}\right)$, then
\begin{eqnarray}\label{EQ_LOB_VOL}
{\rm Vol}(O_{\rm L\ddot{o}bell}(n)) = \frac{n}{2}\left(2\Lambda(\theta)+\Lambda\left(\theta+\frac{\pi}{n}\right)+\Lambda\left(\theta-\frac{\pi}{n}\right)
-\Lambda\left(2\theta+\frac{\pi}{2}\right)\right).
\end{eqnarray}
where $\Lambda$ is the  Lobachevskii function.

Notice that for each of the L\"obell orbifolds that we computed, the volume computed in SnapPea agrees
perfectly (within the six digits of precision available) with that given by Equation \ref{EQ_LOB_VOL}.

\subsection{An orbifold due to Mednykh and Vesnin}
In a very similar way to the construction of L\"obell manifolds, Mednykh and
Vesnin describe in \cite{HYPER_ELL} a compact three-dimensional hyperbolic
manifold $G$ which forms a $2$-fold branched cover over the sphere
$\mathbb{S}^3$.  They call manifolds with such a covering property over
$\mathbb{S}^3$ ``hyperelliptic,'' generalizing the classical notion of
hyperelliptic Riemann surfaces.  See also \cite{MEDNYKH,HYPER_ELL2,HYPER_ELL3}.

The combinatorial polyhedron considered by Mednykh and Vesnin (and apparently
originally due to Grinbergs) is depicted below.

\begin{center}
\begin{picture}(0,0)%
\epsfig{file=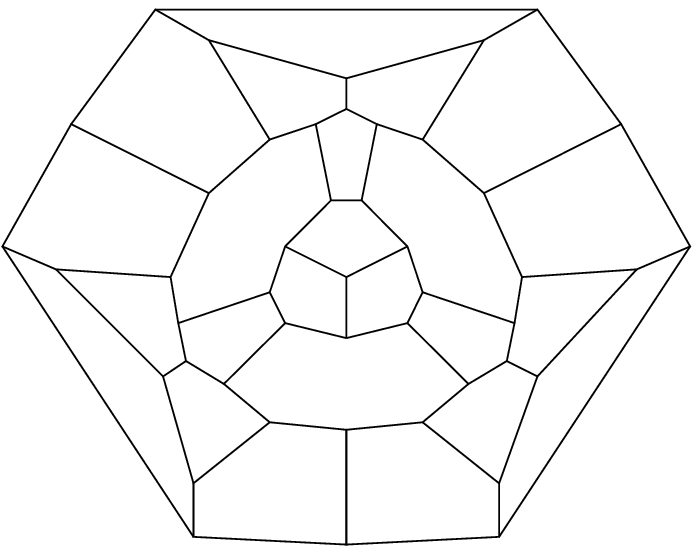}%
\end{picture}%
\setlength{\unitlength}{3947sp}%
\begingroup\makeatletter\ifx\SetFigFont\undefined%
\gdef\SetFigFont#1#2#3#4#5{%
  \reset@font\fontsize{#1}{#2pt}%
  \fontfamily{#3}\fontseries{#4}\fontshape{#5}%
  \selectfont}%
\fi\endgroup%
\begin{picture}(3325,2592)(1564,-3488)
\end{picture}%
\end{center}

This abstract polyhedron has no prismatic $3$-circuits or prismatic
$4$-circuits, so Andreev's Theorem garuntees the existence of a polyhedron
$R_{MV}$ realizing it with $\pi/2$ dihedral angles.  We denote the group generated
by reflections in the faces of $R_{MV}$ by $\Gamma_{MV}$ and the orbifold by $O_{MV}$.
Combinatorial details on the construction of $M_{MV}$ as a $16$-fold cover of $O_{MV}$ can be found
in \cite{HYPER_ELL}.

The following table contains invariants of the orbifold $O_{MV} =
\mathbb{H}^3/\Gamma_{MV}$ obtained by entering an explicit list of generators for
$\Gamma_{MV}$ into SnapPea.

\vspace{.1in}
\begin{tabular}{|l|l|l|}
\hline
$O_{MV}$ & Computed Volume: $6.023046$ & Theoretical Volume: unknown \\
\hline
Short Geodesics & 9 mI  $0.989308$  & 3 mI  $0.989308  +  i \cdot \pi$  \\
& 9 mI  $1.183451$   &  3 mI  $1.183451  +  i \cdot \pi$        \\ 
 & 18 mI  $1.834468$ & 6 mI  $1.834468  +  i \cdot\pi$     \\ 
 & 18 mI $1.859890$ & 6 mI  $1.859890+  i \cdot\pi$      \\ 
 & 27 mI  $1.882318$ &  9 mI  $1.882318  +  i \cdot\pi$     \\ 
 &  6 mI  $1.978616$  & 3 mI  $1.978616  +  i \cdot\pi$     \\ 
 &  9 mI  $2.214787$  & 3 mI  $2.214787 + i \cdot \pi$ \\
 &  18 mI $2.252719$  & 6 mI $2.252719+ i \cdot \pi$ \\
 &  6 mI $2.366902$  & 3 mI $2.366902 + i \cdot \pi$ \\
 &  6 mI  $2.433170$ & 6 mI $2.433170 + i \cdot \pi$ \\
 &  6 mI  $2.446977$ & 6 mI $2.446977 + i \cdot \pi$ \\ 
\hline
\end{tabular}
\vspace{.1in}

As an application, we obtain the estimate ${\rm Vol}(M_{MV}) = 16\cdot 15.608119 = 249.729904$ using that 
$M_{MV}$ is a $16$-fold orbifold cover over $O_{MV}$.

\subsection{Spectral Staircases}

For a given hyperbolic manifold or orbifold $M$, the ``spectral staircase'' is a plot of the number of closed
geodesics of length less than $l$, which we denote $N(l)$, versus $l$.  (In fact, it is much more common
to plot $\log(N(l))$ due to the exponential grown predicted by (\ref{EQ_MARG}) below.) 
The spectral staircase provides both a nice way to graphically display the spectrum of $M$ and a 
illustration of the classical result of Margulis \cite{MARGULIS}, who proves the following universal
formula for the asymptotics of $N(l)$:

\begin{eqnarray}\label{EQ_MARG}
N(l) \sim \frac{exp(\tau l)}{\tau l}  \mbox{   as  } l \rightarrow \infty
\end{eqnarray}

\noindent
where the constant $\tau$ is the topological entropy, which for hyperbolic space $\mathbb{H}^d$ is given
by $\tau = d-1$.  For an exposition and nice experimental works considering spectral staircases, see \cite{INOUE} and
the references within.

\begin{figure}\label{SPEC_STAIR}
\includegraphics[scale = 0.90]{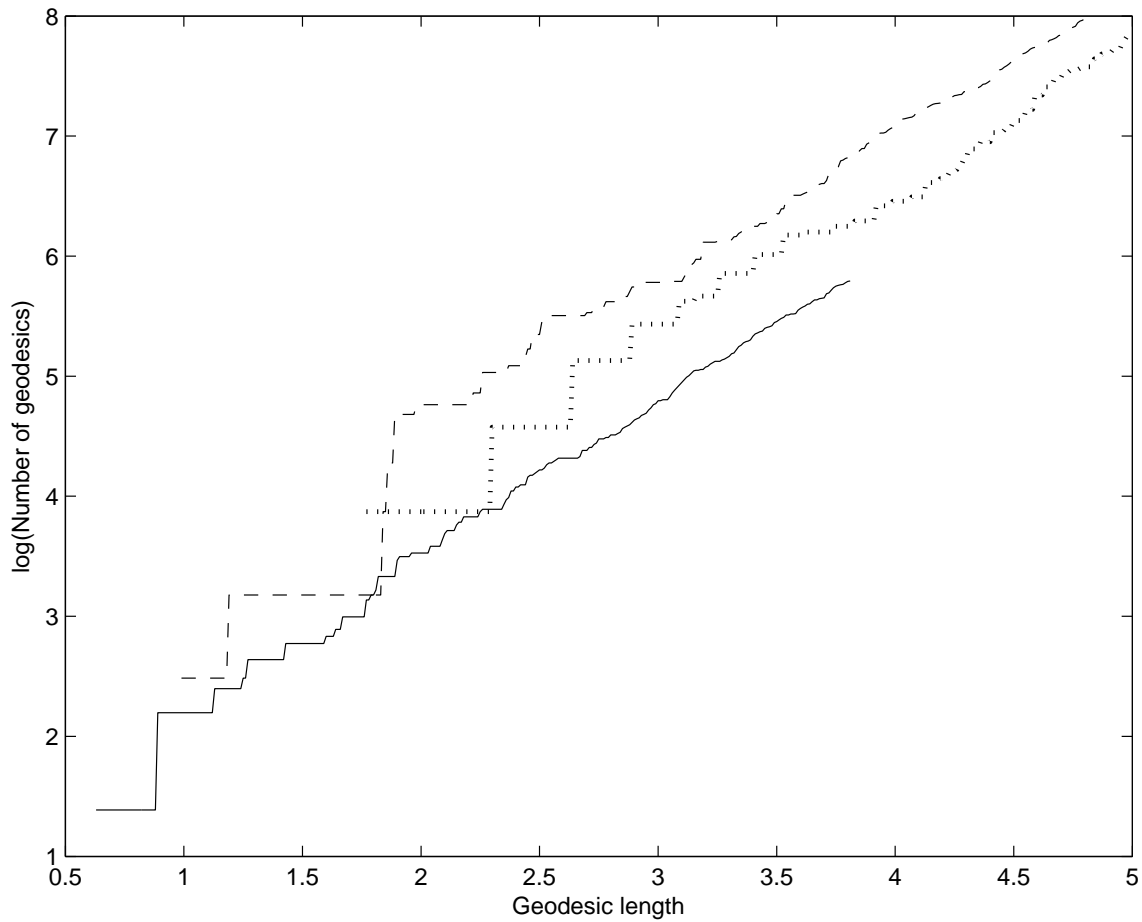}
\caption{Spectral Staircases for the $O_{\rm Lambert}(3,4,5)$ (solid line), $O_{\rm L\ddot{o}bell}(6)$ (dotted line), and $O_{MV}$ (dashed line).}
\end{figure}

We compute these spectral staircases for $O_{\rm Lambert}(3,4,5)$, $O_{\rm L\ddot{o}bell}(6)$, and $O_{MV}$
displaying the results in Figure \ref{SPEC_STAIR}.  (The data for $O_{\rm Lambert}(3,4,5)$ ends at roughly $l = 3.8$.   SnapPea encounters
an error computing at this length, probably due to the comparatively small dihedral angles of $O_{\rm Lambert}(3,4,5)$.)

\section{Questions for further study}
We present a non-comprehensive list of interesting questions for further study:

\begin{itemize}
\item[1] Determine if there is a faster way of computing Andreev Polyhedra.
\item[1b] Related question: determine if CirclePack \cite{CIRCLE_PACK} can be used to construct compact hyperbolic polyhedra and their
reflection groups.  If so, it would provide a faster method of construction.

\item[2] Construct manifold covers of the polyhedral orbifolds that were
considered in Section \ref{SEC_APPLICATIONS}, including the L\"obell Manifolds
\cite{VESNIN_LOB}  the ``Small Covers of the Dodecahedron''
\cite{SMALL_COVERS}, and the Hyperelliptic Manifold \cite{HYPER_ELL}.  Such a
construction, could potentially lead to computations of many additional
interesting invariants of these manifolds using SnapPea, as well as on drilling
and Dehn fillings on them (which would also be possible in SnapPea).

\item[2b] Related question: use the program SNAP \cite{SNAP} to compute arithmetic invariants for these manifolds.

\item[2c] Related question: using SNAP, or the ideas used in SNAP \cite{SNAP_PAPER},  study the arithmetic invariants of polyhedral reflection groups.

\item[3] Perform a study of volumes of hyperbolic polyhedra corresponding to general angles in $A_C$.  
(While SnapPea computes volumes only for polyhedra with discrete reflection groups, the functions from the SnapPea kernel could probably
be used for this more general study.)

\end{itemize}

\appendix
\vspace{.25in}
\noindent
{\Large \bf Appendix:}

\section{Using the program to construct compact polyhedra and  their reflection groups}
The computer program described in this paper 
is a functional but slightly rough collection of Matlab \cite{MAT} (or Octave
\cite{OCTAVE}) scripts.  Using a single command the program produces a
polyhedron realizing a simple abstract polyhedron $C$ with all dihedral angles
$\frac{2\pi}{5}$.  However, one must do a little bit more work to construct a
polyhedron realizing $(C,{\bf a})$ of truncated or compound type.  These later
steps are not as automatic in the program, but one can follow the description in this paper
step by step to do them ``by hand.''  Please see the README file enclosed with the
program for further information.

There are two ways to output a polyhedron that has been constructed using this program:
the {\em Object File Format} (e.g. ``filename.off'') and the {\em Generators File
Format} (e.g. ``filename.txt'').  The {\em Object File Format} output can be
read into Geomview \cite{GEO} and displayed nicely in the Poincar\'e Ball
model there.  The  {\em Generators File Format} is for SnapPea \cite{SNAPPEA}
and (if the polyhedron output this way has dihedral angles that are proper
integer sub-multiples of $\pi$) this file can be loaded into SnapPea. 

\vspace{0.1in}

An analysis of the computational complexity of this method would be quite
involved and is not feasible at this point.  In fact, the
computational complexity of Newton's Method is quite difficult
\cite{SHUB1,SHUB2,SHUB3,SHUB4,BELTRAN}.  

In practice, on a contemporary PC running Linux, the most complicated
construction in this paper took approximately 2 minutes; it is by no means
fast, but certainly usable.  We expect that the program described in this paper
can be used to construct all polyhedra having no more than 30 faces and having
dihedral angles bounded away from $\partial A_C$ (with the exception of the
part of $\partial A_C$ corresponding to dihedral angles $\pi/2$, where there
should be no problem.) Certainly the number of iteration steps in each homotopy
(the parameter $K$  from Section \ref{SEC_DEFORM}) and the parameter $\epsilon$ used in the Whitehead
move (see Section \ref{SEC_WH}) may need to be modified in special
circumstances.  With improvements of the program, perhaps by implementing it in
a faster programming language, and, if necessary, in a higher precision
arithmetic, we expect the same to be possible for up to 100 faces, or possibly
more.

\vspace{.2in}
\noindent
{\bf \Large Acknowledgments:}

The author gratefully thanks John H. Hubbard for suggesting that he write this
program and for providing an outline for using the homotopy method in
combination with Andreev's proof.  It was this suggestion that led to the
discovery of the error in Andreev's proof \cite{ROE2,ROE} and eventually to a
working program.  The author has greatly benefited from many discussions with
William D. Dunbar and he also thanks Rodrigo Perez and Travis Waddington for
providing suggestions about the manuscript and attached computer files.

The computer program SnapPea, written by Jeff Weeks and his collaborators has
allowed for the exciting applications in Section \ref{SEC_APPLICATIONS} of this
paper.  The author thanks Jeff Weeks and his collaborators for writing this
wonderful program, and also thanks Jeff Weeks for helpful personal
communications.

The author also gratefully thanks the referee who offered many suggestions
which have led to a clearer explanation of why we use the technique that we
use.  His or her comments have also provided the inspiration for the
last two sections of this paper, and significant further research in these directions.

\bibliographystyle{plain}
\bibliography{andreev.bib,newton.bib}

\begin{thebibliography}{10}

\bibitem{MAT}
www.mathworks.com{\tt}.

\bibitem{OCTAVE}
www.octave.org.

\bibitem{GEO}
www.geomview.org, Developed by The Geometry Center at the University of
  Minnesota in the late 1990's.

\bibitem{ALLGOWER}
Eugene~L. Allgower and Kurt Georg.
\newblock {\em Numerical Continuation Methods, an introduction}.
\newblock Springer-Verlag, 1990.

\bibitem{AND}
E.~M. Andreev.
\newblock On convex polyhedra in {L}obacevskii spaces ({E}nglish
  {T}ranslation).
\newblock {\em Math. USSR Sbornik}, 10:413--440, 1970.

\bibitem{AND2}
E.~M. Andreev.
\newblock On convex polyhedra in {L}obacevskii spaces (in {R}ussian).
\newblock {\em Mat. Sb.}, 81(123):445--478, 1970.

\bibitem{BLUM}
Lenore Blum, Felipe Cucker, Michael Shub, and Steve Smale.
\newblock {\em Complexity and real computation}.
\newblock Springer-Verlag, New York, 1998.
\newblock With a foreword by Richard M. Karp.

\bibitem{STEVE}
Phil Bowers and Kenneth Stephenson.
\newblock A branched {A}ndreev-{T}hurston theorem for circle packings of the
  sphere.
\newblock {\em Proc. London Math. Soc. (3)}, 73(1):185--215, 1996.

\bibitem{H}
C.~D. Hodgson.
\newblock Deduction of {A}ndreev's theorem from {R}ivin's characterization of
  convex hyperbolic polyhedra.
\newblock In {\em Topology 90}, pages 185--193. de Gruyter, 1992.

\bibitem{HP}
John~H. Hubbard and Peter Papadopol.
\newblock Newton's method applied to two quadratic equations in $\mathbb{C}^2$
  viewed as a global dynamical system.
\newblock To appear in the Memoires of the AMS.

\bibitem{HBH}
John~Hamal Hubbard and Barbara~Burke Hubbard.
\newblock {\em Vector calculus, linear algebra, and differential forms}.
\newblock Prentice Hall Inc., Upper Saddle River, NJ, 1999.
\newblock A unified approach.

\bibitem{ROE2}
Roland K. W. Roeder John~H. Hubbard and William~D. Dunbar.
\newblock Andreev's theorem on hyperbolic polyhedra.
\newblock Submitted to {\it Les Annales de l'Institute Fourier}, Dec 2004.

\bibitem{KANT}
L.~V. Kantorovi{\v{c}}.
\newblock On {N}ewton's method.
\newblock {\em Trudy Mat. Inst. Steklov.}, 28:104--144, 1949.

\bibitem{MR}
A.~Marden and B.~Rodin.
\newblock On {T}hurston's formulation and proof of {A}ndreev's {T}heorem.
\newblock In {\em Computational Methods and Function Theory}, volume 1435 of
  {\em Lecture Notes in Mathematics}, pages 103--115. Springer-Verlag, 1990.

\bibitem{RH}
I.~Rivin and C.~D. Hodgson.
\newblock A characterization of compact convex polyhedra in hyperbolic 3-space.
\newblock {\em Invent. Math.}, 111:77--111, 1993.

\bibitem{ROE_TET}
Roland K.~W. Roeder.
\newblock Compact hyperbolic tetrahedra with non-obtuse dihedral angles.
\newblock to appear.

\bibitem{RNEWT}
Roland K.~W. Roeder.
\newblock A degenerate newton's map in two complex variables: linking with
  currents.
\newblock submitted to J. Geometric Analysis.

\bibitem{ROE}
Roland K.~W. Roeder.
\newblock Le th\'eor\`eme d'andreev sur poly\`edres hyperboliques.
\newblock Doctoral thesis (In English), May 2004.
\newblock Universit\'e de Provence, Aix-Marseille 1.

\bibitem{ROE_WEB}
http://www.fields.utoronto.ca/$\sim$rroeder/andreev\_auxilliaries.tar.gz

\bibitem{SHUB1}
Michael Shub and Steve Smale.
\newblock Complexity of {B}\'ezout's theorem. {I}. {G}eometric aspects.
\newblock {\em J. Amer. Math. Soc.}, 6(2):459--501, 1993.

\bibitem{SHUB2}
M.~Shub and S.~Smale.
\newblock Complexity of {B}ezout's theorem. {II}. {V}olumes and probabilities.
\newblock In {\em Computational algebraic geometry (Nice, 1992)}, volume 109 of
  {\em Progr. Math.}, pages 267--285. Birkh\"auser Boston, Boston, MA, 1993.


\bibitem{SHUB3}
Michael Shub and Steve Smale.
\newblock Complexity of {B}ezout's theorem. {III}. {C}ondition number and
  packing.
\newblock {\em J. Complexity}, 9(1):4--14, 1993.
\newblock Festschrift for Joseph F. Traub, Part I.

\bibitem{SHUB4}
Michael Shub and Steve Smale.
\newblock Complexity of {B}ezout's theorem. {IV}. {P}robability of success;
  extensions.
\newblock {\em SIAM J. Numer. Anal.}, 33(1):128--148, 1996.

\bibitem{TH}
W.~P. Thurston.
\newblock {\em The Geometry and Topology of 3-manifolds}.
\newblock Princeton University Notes, Princeton, New Jersey, 1980.

\end{thebibliography}


\begin{thebibliography}{10}

\bibitem{MAT}
www.mathworks.com{\tt}.

\bibitem{OCTAVE}
www.octave.org.

\bibitem{GEO}
www.geomview.org, Developed by The Geometry Center at the University of
  Minnesota in the late 1990's.

\bibitem{WEEKS}
Colin Adams, Martin Hildebrand, and Jeffrey Weeks.
\newblock Hyperbolic invariants of knots and links.
\newblock {\em Trans. Amer. Math. Soc.}, 326(1):1--56, 1991.

\bibitem{AVS}
D.~V. Alekseevskij, {\`E}.~B. Vinberg, and A.~S. Solodovnikov.
\newblock Geometry of spaces of constant curvature.
\newblock In {\em Geometry, II}, volume~29 of {\em Encyclopaedia Math. Sci.},
  pages 1--138. Springer, Berlin, 1993.

\bibitem{ALLGOWER}
Eugene~L. Allgower and Kurt Georg.
\newblock {\em Numerical Continuation Methods, an introduction}.
\newblock Springer-Verlag, 1990.

\bibitem{BELTRAN}
Carlos~Beltr\'an \'Alvarez.
\newblock Sobre el {P}roblema 17 de {S}male: Teor\'ia de la intersecci\'on y
  geometr\'ia integral.
\newblock Doctoral thesis (In Spanish), September 2006.
\newblock Universidad de Cantabria.

\bibitem{AND}
E.~M. Andreev.
\newblock On convex polyhedra in {L}obacevskii spaces ({E}nglish
  {T}ranslation).
\newblock {\em Math. USSR Sbornik}, 10:413--440, 1970.

\bibitem{AND2}
E.~M. Andreev.
\newblock On convex polyhedra in {L}obacevskii spaces (in {R}ussian).
\newblock {\em Mat. Sb.}, 81(123):445--478, 1970.

\bibitem{BAO}
Xiliang Bao and Francis Bonahon.
\newblock Hyperideal polyhedra in hyperbolic 3-space.
\newblock {\em Bull. Soc. Math. France}, 130(3):457--491, 2002.

\bibitem{BLUM}
Lenore Blum, Felipe Cucker, Michael Shub, and Steve Smale.
\newblock {\em Complexity and real computation}.
\newblock Springer-Verlag, New York, 1998.
\newblock With a foreword by Richard M. Karp.

\bibitem{STEVE}
Phil Bowers and Kenneth Stephenson.
\newblock A branched {A}ndreev-{T}hurston theorem for circle packings of the
  sphere.
\newblock {\em Proc. London Math. Soc. (3)}, 73(1):185--215, 1996.

\bibitem{LUO}
Bennett Chow and Feng Luo.
\newblock Combinatorial {R}icci flows on surfaces.
\newblock {\em J. Differential Geom.}, 63(1):97--129, 2003.

\bibitem{SNAP_PAPER}
David Coulson, Oliver~A. Goodman, Craig~D. Hodgson, and Walter~D. Neumann.
\newblock Computing arithmetic invariants of 3-manifolds.
\newblock {\em Experiment. Math.}, 9(1):127--152, 2000.

\bibitem{DIAZ}
Raquel D{\'{\i}}az.
\newblock Non-convexity of the space of dihedral angles of hyperbolic
  polyhedra.
\newblock {\em C. R. Acad. Sci. Paris S\'er. I Math.}, 325(9):993--998, 1997.

\bibitem{DIAZ_ANDREEV}
Raquel D{\'{\i}}az.
\newblock A generalization of {A}ndreev's theorem.
\newblock {\em J. Math. Soc. Japan}, 58(2):333--349, 2006.

\bibitem{SMALL_COVERS}
Anne Garrison and Richard Scott.
\newblock Small covers of the dodecahedron and the 120-cell.
\newblock {\em Proc. Amer. Math. Soc.}, 131(3):963--971 (electronic), 2003.

\bibitem{SNAP}
O.~A. Goodman, C.~D. Hodgson, and Neumann w.~D.
\newblock ``snap home page''.
\newblock See http://www.ms.unimelb.edu.au/$\sim$snap. Includes source
  distribution and extensive tables of results of Snap computations.

\bibitem{GUE}
Fran{\c{c}}ois Gu{\'e}ritaud.
\newblock On an elementary proof of {R}ivin's characterization of convex ideal
  hyperbolic polyhedra by their dihedral angles.
\newblock {\em Geom. Dedicata}, 108:111--124, 2004.

\bibitem{POINCARE}
Poincar\'e Henri.
\newblock M\'emoire sur les groups klein\'eens.
\newblock {\em Acta Math}, 3:49--92, 1883.

\bibitem{H}
C.~D. Hodgson.
\newblock Deduction of {A}ndreev's theorem from {R}ivin's characterization of
  convex hyperbolic polyhedra.
\newblock In {\em Topology 90}, pages 185--193. de Gruyter, 1992.

\bibitem{WEEKS_LENGTHS}
Craig~D. Hodgson and Jeffrey~R. Weeks.
\newblock Symmetries, isometries and length spectra of closed hyperbolic
  three-manifolds.
\newblock {\em Experiment. Math.}, 3(4):261--274, 1994.

\bibitem{HP}
John~H. Hubbard and Peter Papadopol.
\newblock Newton's method applied to two quadratic equations in $\mathbb{C}^2$
  viewed as a global dynamical system.
\newblock To appear in the Memoires of the AMS.

\bibitem{HBH}
John~Hamal Hubbard and Barbara~Burke Hubbard.
\newblock {\em Vector calculus, linear algebra, and differential forms}.
\newblock Prentice Hall Inc., Upper Saddle River, NJ, 1999.
\newblock A unified approach.

\bibitem{ROE2}
Roland K. W. Roeder John~H. Hubbard and William~D. Dunbar.
\newblock Andreev's theorem on hyperbolic polyhedra.
\newblock To appear {\it Les Annales de l'Institute Fourier}.

\bibitem{INOUE}
Kaiki~Taro Inoue.
\newblock Numerical study of length spectra and low-lying eigenvalue spectra of
  compact hyperbolic 3-manifolds.
\newblock {\em Classical Quantum Gravity}, 18(4):629--652, 2001.

\bibitem{KANT}
L.~V. Kantorovi{\v{c}}.
\newblock On {N}ewton's method.
\newblock {\em Trudy Mat. Inst. Steklov.}, 28:104--144, 1949.

\bibitem{KELLERHALS}
Ruth Kellerhals.
\newblock On the volume of hyperbolic polyhedra.
\newblock {\em Math. Ann.}, 285(4):541--569, 1989.

\bibitem{LOBELL}
F.~L\"obell.
\newblock Beispiele geschlossener dreidimensionaler clifford-kleinische r\"aume
  negativer kr\"ummung.
\newblock {\em Ber. S\"achs. Akad. Wiss.}, 83:168--174, 1931.

\bibitem{MR}
A.~Marden and B.~Rodin.
\newblock On {T}hurston's formulation and proof of {A}ndreev's {T}heorem.
\newblock In {\em Computational Methods and Function Theory}, volume 1435 of
  {\em Lecture Notes in Mathematics}, pages 103--115. Springer-Verlag, 1990.

\bibitem{MARGULIS}
G.~A. Margulis.
\newblock Certain applications of ergodic theory to the investigation of
  manifolds of negative curvature.
\newblock {\em Functional Anal. Appl.}, 3(4):335--336, 1969.

\bibitem{MEDNYKH}
A.~D. Mednykh.
\newblock Three-dimensional hyperelliptic manifolds.
\newblock {\em Ann. Global Anal. Geom.}, 8(1):13--19, 1990.

\bibitem{HYPER_ELL3}
Alexander Mednykh and Marco Reni.
\newblock Twofold unbranched coverings of genus two 3-manifolds are
  hyperelliptic.
\newblock {\em Israel J. Math.}, 123:149--155, 2001.

\bibitem{HYPER_ELL2}
Alexander Mednykh, Marco Reni, Andrei Vesnin, and Bruno Zimmermann.
\newblock Three-fold coverings and hyperelliptic manifolds: a three-dimensional
  version of a result of {A}ccola.
\newblock {\em Rend. Istit. Mat. Univ. Trieste}, 32(suppl. 1):181--191 (2002),
  2001.
\newblock Dedicated to the memory of Marco Reni.

\bibitem{HYPER_ELL}
Alexander Mednykh and Andrei Vesnin.
\newblock Colourings of polyhedra and hyperelliptic 3-manifolds.
\newblock In {\em Recent advances in group theory and low-dimensional topology
  (Pusan, 2000)}, volume~27 of {\em Res. Exp. Math.}, pages 123--131.
  Heldermann, Lemgo, 2003.

\bibitem{RENI}
Marco Reni.
\newblock Dihedral branched coverings of hyperbolic orbifolds.
\newblock {\em Geom. Dedicata}, 67(3):271--283, 1997.

\bibitem{RH}
I.~Rivin and C.~D. Hodgson.
\newblock A characterization of compact convex polyhedra in hyperbolic 3-space.
\newblock {\em Invent. Math.}, 111:77--111, 1993.

\bibitem{RIV_IDEAL1}
Igor Rivin.
\newblock On geometry of convex ideal polyhedra in hyperbolic {$3$}-space.
\newblock {\em Topology}, 32(1):87--92, 1993.

\bibitem{RIV_IDEAL2}
Igor Rivin.
\newblock A characterization of ideal polyhedra in hyperbolic {$3$}-space.
\newblock {\em Ann. of Math. (2)}, 143(1):51--70, 1996.

\bibitem{RNEWT}
Roland K.~W. Roeder.
\newblock A degenerate newton's map in two complex variables: linking with
  currents.
\newblock To appear J. Geometric Analysis.

\bibitem{ROE}
Roland K.~W. Roeder.
\newblock Le th\'eor\`eme d'andreev sur poly\`edres hyperboliques.
\newblock Doctoral thesis (In English), May 2004.
\newblock Universit\'e de Provence, Aix-Marseille 1.

\bibitem{ROE_TET}
Roland K.~W. Roeder.
\newblock Compact hyperbolic tetrahedra with non-obtuse dihedral angles.
\newblock {\em Publicacions Matem\`atiques}, 50(1):211--227, 2006.

\bibitem{SCH2}
J.-M. Schlenker.
\newblock Dihedral angles of convex polyhedra.
\newblock {\em Discrete Comput. Geom.}, 23(3):409--417, 2000.

\bibitem{SCH1}
Jean-Marc Schlenker.
\newblock M\'etriques sur les poly\`edres hyperboliques convexes.
\newblock {\em J. Differential Geom.}, 48(2):323--405, 1998.

\bibitem{SCH3}
Jean-Marc Schlenker.
\newblock Hyperbolic manifolds with convex boundary.
\newblock {\em Invent. Math.}, 163(1):109--169, 2006.

\bibitem{SHUB2}
M.~Shub and S.~Smale.
\newblock Complexity of {B}ezout's theorem. {II}. {V}olumes and probabilities.
\newblock In {\em Computational algebraic geometry (Nice, 1992)}, volume 109 of
  {\em Progr. Math.}, pages 267--285. Birkh\"auser Boston, Boston, MA, 1993.

\bibitem{SHUB1}
Michael Shub and Steve Smale.
\newblock Complexity of {B}\'ezout's theorem. {I}. {G}eometric aspects.
\newblock {\em J. Amer. Math. Soc.}, 6(2):459--501, 1993.

\bibitem{SHUB3}
Michael Shub and Steve Smale.
\newblock Complexity of {B}ezout's theorem. {III}. {C}ondition number and
  packing.
\newblock {\em J. Complexity}, 9(1):4--14, 1993.
\newblock Festschrift for Joseph F. Traub, Part I.

\bibitem{SHUB4}
Michael Shub and Steve Smale.
\newblock Complexity of {B}ezout's theorem. {IV}. {P}robability of success;
  extensions.
\newblock {\em SIAM J. Numer. Anal.}, 33(1):128--148, 1996.

\bibitem{CIRCLE_PACK}
Ken Stephenson.
\newblock `{C}ircle{P}ack''.
\newblock See http://www.math.utk.edu/~kens/CirclePack/.

\bibitem{TH}
W.~P. Thurston.
\newblock {\em The Geometry and Topology of 3-manifolds}.
\newblock Princeton University Notes, Princeton, New Jersey, 1980.

\bibitem{VESNIN_LOB}
Andrei Vesnin.
\newblock Three-dimensional hyperbolic manifolds of {L}\"obell type.
\newblock {\em Siberian Math. J}, 28(5):731--733, 1987.

\bibitem{LOB_VOL}
Andrei Vesnin.
\newblock Volumes of {L}\"obell $3$-manifolds.
\newblock {\em Math. Notes}, 64(1-2):15--19, 1998.

\bibitem{VIN}
{\`E}.~B. Vinberg.
\newblock Discrete groups generated by reflections in {L}oba\v cevski\u\i \
  spaces.
\newblock {\em Mat. Sb. (N.S.)}, 72 (114):471--488; correction, ibid. { 73
  (115) (1967), 303}, 1967.

\bibitem{VINREFL}
{\`E}.~B. Vinberg.
\newblock Hyperbolic groups of reflections.
\newblock {\em Russian Math. Surveys}, 40(1):31--75, 1985.

\bibitem{VINVOL}
{\`E}.~B. Vinberg.
\newblock The volume of polyhedra on a sphere and in {L}obachevsky space.
\newblock In {\em Algebra and analysis (Kemerovo, 1988)}, volume 148 of {\em
  Amer. Math. Soc. Transl. Ser. 2}, pages 15--27. Amer. Math. Soc., Providence,
  RI, 1991.

\bibitem{VS}
{\`E}.~B. Vinberg and O.~V. Shvartsman.
\newblock Discrete groups of motions of spaces of constant curvature.
\newblock In {\em Geometry, II}, volume~29 of {\em Encyclopaedia Math. Sci.},
  pages 139--248. Springer, Berlin, 1993.

\bibitem{SNAPPEA}
J.~R. Weeks.
\newblock ``{S}nap{P}ea''.
\newblock See http://geometrygames.org/SnapPea/index.html.

\end{thebibliography}
\end{document}